\documentclass[12pt]{amsart}
\usepackage[colorlinks, linkcolor=blue, anchorcolor=blue, citecolor=green]{hyperref}
\usepackage[normalem]{ulem}
\usepackage{graphics,epic}
\usepackage{amsmath,amssymb, amsthm,mathrsfs}
\usepackage[all,2cell]{xy}
\usepackage{multicol}
\usepackage{multirow}
\usepackage{caption}
\usepackage{subcaption}
\usepackage{shorttoc}
\usepackage{url}
\usepackage{tikz}
\usepackage{geometry}      
\usepackage{graphicx}

\newtheorem{theorem}{Theorem}[section]
\newtheorem*{theorem*}{Theorem}

\newtheorem{lemma}[theorem]{Lemma}

\newtheorem{proposition}[theorem]{Proposition}
\newtheorem{corollary}[theorem]{Corollary}

\newtheorem*{conjecture*}{Conjecture}

\newtheorem{remark}[theorem]{Remark}
\newtheorem{definition}[theorem]{Definition}

\newtheorem{lem}[theorem]{Lemma}

\numberwithin{equation}{section}
\numberwithin{figure}{subsection}
\begin{document}

\title[Pseudo grading on cluster automorphism group]{Pseudo grading on cluster automorphism group with application to cluster algebras of rank $3$}

\author{Changjian Fu}
\address{Changjian Fu\\Department of Mathematics\\SiChuan University\\610064 Chengdu\\P.R.China}
\email{changjianfu@scu.edu.cn}
\author{Zhanhong Liang}
\address{Zhanhong Liang\\Department of Mathematics\\SiChuan University\\610064 Chengdu\\P.R.China}
\email{1278074468@qq.com}
\thanks{This work is partially supported by the National Natural Science Foundation of China (Grant No. 11971326, 12471037).}

\begin{abstract}
We introduce a pseudo $\mathbb{N}$-grading on the cluster auotmorphism group $\operatorname{Aut}(\mathcal{A})$ with respect to an initial seed of $\mathcal{A}$, which consists of a family of subsets $\{G_i\}_{i\in \mathbb{N}}$ of $\operatorname{Aut}(\mathcal{A})$ such that $\operatorname{Aut}(\mathcal{A})=\bigcup_{i\in \mathbb{N}}G_i$ and $G_k\cdot G_l\subset \bigcup_{i=0}^{k+l}G_i$. We prove that $\operatorname{Aut}(\mathcal{A})$ is generated by $G_0\cup G_1$, leading to an elementary approach for 
calculating cluster automorphism groups of certain cluster algebras. As an application, we completely determined the cluster automorphism groups of cluster algebras of rank $3$ with indecomposable exchange matrices.

\end{abstract}
\maketitle

\tableofcontents

\section{Introduction}
  A cluster algebra \cite{FZ02} is a commutative algebra with an additional combinatorial structure. In particular, it has a countable family of generators (called {\em cluster variables}), which are grouped into overlapping sets (called {\em clusters}) of equal finite size.  All the cluster variables can be obtained from an initial cluster via mutation.
From the ring-theoretic viewpoint, Assem, Schiffler and Shramchenko \cite{ASS2012} proposed the investigation of homomorphisms between cluster algebras and introduced the concept of cluster auotmorphism. By definition, a cluster automorphism is a homomorphism of algebras that is compatible with cluster structure, i.e., it maps a cluster to a cluster and preserves mutations. 
Since their introduction, cluster automorphism and their generalizations \cite{Fr16,CHL24} have attracted significant attention \cite{CZ16b,CS19,CZ20,CLLP19,LL21} and have proven  to be extremely useful in the study of the structure theory of cluster algebras, cf. \cite{ASS14,BM16,BMP20,CL2020}.

The set  of all cluster automorphisms of a cluster algebra forms a group under the composition, referred to as the {\em cluster automorphism group}. 
Since this group describes the symmetry of the cluster algebra with respect to its cluster structure, it is important to explicitly characterize it. However, it seems unlikely that such a characterization is possible for an arbitrary cluster algebra.
 To the best of our knowledge, there is no general method to compute the cluster automorphism group for a general cluster algebra.

Nevertheless, significant work have been done to describe cluster automorphism groups for certain important classes of cluster algebras.  Assem, Schiffler and Shramchenko \cite{ASS2012} studied the cluster automorphism groups for acyclic skew-symmetric cluster algebras via representation theory of the underlying path algebras, computing these groups explicitly for simply-laced Dynkin and Euclidean types. Using folding techniques, Chang and Zhu \cite{CZ16a} determined the cluster automorphism groups explicitly for cluster algebras of finite type. 
Blanc and Dolgachev \cite{BD15} computed the cluster automorphism groups for cluster algebras of rank two using geometric tools.
On the other hand, for a cluster algebra arising from a marked surface, Assem et al \cite{ASS2012} first observed that the cluster automorphism is deeply connected to  the mapping class group of the  marked surface. Gu \cite{Gu11}, Bridgeland and Smith \cite{BS15} further established an explicit relation between cluster automorphism groups and  mapping class groups, except for the 4-punctured sphere, the once-punctured 4-gon and the
 twice-punctured digon. Buliding on the work \cite{Gu11,BS15}, Dong and Li \cite{DL23}  obtain the presentations of cluster autormorphism groups for all cluster algebras arising from  marked surfaces.


Let $\mathcal{A}$ be a cluster algebra of rank $n$ with trivial coefficients, and let $\mathbf{\Sigma}=\{\Sigma_t=(\mathbf{x}_t, B_t)\}_{t\in \mathbb{T}_n}$ be a cluster pattern of $\mathcal{A}$ with root vertex $t_0$. 
By observing that each cluster automorphism $f$ of $\mathcal{A}$ can be expressed as $f=(\Sigma_{t_0},\Sigma_{t},\sigma, \varepsilon)$ for some $t\in \mathbb{T}_n$, $\sigma\in S_n$ and $\varepsilon\in \{\pm\}$, we introduce a weight on the path $p(t_0,t)$ and obtain a pseudo $\mathbb{N}$-grading on the cluster automorphism group $\operatorname{Aut}(\mathcal{A})$.  This grading consists of a family of subsets $\{G_i(t_0)\}_{i\in \mathbb{N}}$ of $\operatorname{Aut}(\mathcal{A})$. 
It is worth mentioning that $G_0(t_0)$ is the group of symmetries of the initial exchange matrix $B_{t_0}$, which is a subgroup of $S_n$ and $\operatorname{Aut}(\mathcal{A})$. By proving a factorization lemma (cf. Lemma \ref{lem:d-vector-initial}) for cluster automorphisms, we show that $\operatorname{Aut}(\mathcal{A})$ is generated by $G_0(t_0)\cup G_1(t_0)$ (cf. Theorem \ref{thm:generator-set}). In most case, the subset $G_1(t_0)$ can be further refined (cf. Proposition \ref{prop:generator-set-assumption}). Once $G_1(t_0)$ is finite, $\operatorname{Aut}(\mathcal{A})$ is finitely generated. This yields a new, elementary method for explicitly computing $\operatorname{Aut}(\mathcal{A})$.
As a first application, we completely determine the cluster automorphism groups for all cluster algebras of rank $3$ with indecomposable exchange matrices. It seems likely this approach can also be applied to study the cluster automorphism groups for cluster algebras of finite mutation type.

The paper is organized as follows. In Section \ref{s:preliminaries}, we recall the basic definitions in cluster algebras and the classification of the mutation classes of skew-symmetrizable $3\times 3$ matrices. In Section \ref{s:cluster-automorphism}, after recalling the definition and properties of cluster automorphism, we introduce a weight for the path associated with a cluster automorphism and obtain the pseudo $\mathbb{N}$-grading of $\operatorname{Aut}(\mathcal{A})$. After establishing the facterization property of cluster automorphism in Lemma \ref{lem:factorization-cluster-automorphism}, we prove that $\operatorname{Aut}(\mathcal{A})$ is generated by the $0$-th  and $1$-st components in Theorem \ref{thm:generator-set} and Proposition \ref{prop:generator-set-assumption}.
As an application of Theorem \ref{thm:generator-set}  and Proposition \ref{prop:generator-set-assumption}, we completely determine the cluster automorphism groups for cluster algebras of rank $3$ with indecomposable exchange matrices in Section \ref{s:cluster-automorphism-rank-3}, see Theorem \ref{thm:cyclic-case} and Theorem \ref{thm:acyclic-case}. The proofs of Theorem \ref{thm:cyclic-case} and Theorem \ref{thm:acyclic-case} are given in Section \ref{s:mutation-cyclic} and \ref{s:mutation-acyclic}, respectively.

\section{Preliminaries}\label{s:preliminaries}
\subsection{Cluster algebras}
In this section, we recall basic definitions in cluster algebras and restrict our attention to cluster algebras with trivial coefficients. We follow \cite{FZ02}. For an integer $b$, we use the notation $[b]_+=\max\{b,0\}$.
Let $n$ be a positive integer. Let $\mathcal F$ be a field isomorphic to the rational function field in $n$ variables with coefficients in the field $\mathbb{Q}$ of rational numbers. This field $\mathcal{F}$ serves as the {\em ambient field} for cluster algebras we will introduce.

\begin{definition}[Seed]
\begin{itemize}
    \item An $n\times n$ integer matrix $B=(b_{ij})^n_{i,j=1}$ is  {\em skew-symmetrizable} if there exists a diagonal matrix $D=\textup{diag}\{d_1,\dots, d_n\}$ whose diagonal entries $d_i$ are positive integers such that $DB$ is skew-symmetric. In this case, $D$ is called a {\em skew-symmetrizer} of $B$.
    \item  A (labeled) seed of rank $n$ is a pair $\Sigma = (\mathbf x,B)$, where $\mathbf{x} = (x_1,\dots,x_n)$ is an $n$-tuple of elements of $\mathcal{F}$ forming a free generating set of $\mathcal{F}$, and $B=(b_{ij})^n_{i,j=1}$ is an $n\times n$ skew-symmetrizable integer matrix. We refer to $\mathbf x$ and $B$ as the {\em cluster} and the {\em exchange matrix} of $\Sigma$, respectively. The variables $x_i$ are called the {\em cluster variables}.
\end{itemize}
\end{definition}

\begin{definition}[Seed mutation]
For any seed $\Sigma = (\mathbf x,B)$ of rank $n$ and $k \in  \{1,\dots,n\}$, the {\em mutation $\mu_k$ in direction $k$} transforms $\Sigma$ into another labeled seed $\mu_k(\Sigma):=(\mathbf{x'},B')$ defined as follows:
\begin{itemize}
    \item The entries of $B':=\mu_k(B):=(b'_{ij})_{i,j=1}^n$ are given by
    \begin{align*}
    b'_{ij}=&
\begin{cases}
-b_{ij}  &\text{$i=k$ or $j=k$}; \\
b_{ij}+b_{ik}[b_{kj}]_+ +[-b_{ik}]_+ b_{kj} &\text{else}.
\end{cases}
\end{align*}
\item The cluster variables $\mathbf{x}'=(x_1',\dots, x_n')$ are given by 
\begin{align*}
x'_i=&
\begin{cases}
x_k^{-1}\left (\prod\limits_{j=1}^n x_j^{[b_{jk}]_+ }+\prod\limits_{j=1}^n x_j^{[-b_{jk}]_+}  \right )   &i=k; \\
x_i  & i\ne k.
\end{cases} 
\end{align*}
\end{itemize}
\end{definition}
It is straightforward to check that $\mu_k$ is an involution, i.e., $\mu_k^2(\Sigma)=\Sigma$. The matrix $\mu_k(B)$ is also called the {\em mutation of $B$ in direction $k$}. Two skew-symmetrizable matrices $B$ and $B'$ are {\em mutation-equivalent}, if $B'$ can be obtained from $B$ by a finite sequence of mutations. A sequence $i_1i_2\cdots i_k$ with $i_j\in \{1,\dots, n\}$ for $1\leq j\leq k$ is called {\em reduced} if $i_l\neq i_{l+1}$ for $1\leq l<k$. For a reduced sequence $i_1i_2\cdots i_k$, we define $\mu_{i_1i_2\cdots i_k}:=\mu_{i_1}\mu_{i_2}\cdots \mu_{i_k}$ as the composition of mutations $\mu_{i_l}$, for $1\leq l\leq k$.

Let $\mathbb{T}_n$ be the $n$-regular tree whose edges are labeled by the numbers $1,\dots, n$ such that the $n$ edges emanating from each vertex carry distinct labels.
For each vertex $t\in \mathbb{T}_n$, there is a unique reduced sequence  $i_k\dots i_1$ such that $t_0$ and $t$ are connected by edges labeled by $i_1,\dots, i_k$, i.e.,
 \[
    \xymatrix{t_0\ar@{-}[r]^{i_1}&\cdot\ar@{-}[r]^{i_2}&\cdot\cdots\ar@{-}[r]^{i_k}&t}. \]
 We also write $t=\mu_{i_k\cdots i_1}(t_0)$ and denote the path from $t_0$ to $t$ by $p(t_0,t)$.

\begin{definition}[Cluster pattern]
 A collection of labeled seeds ${\bf \Sigma }= \{\Sigma_t=(\mathbf x_t, B_t)\}_{t\in \mathbb T_n}$ indexed by the $n$-regular tree $\mathbb T_n$ is called a {\em cluster pattern} or {\em seed pattern} of rank $n$ if for any pair of vertices $t,t' \in \mathbb T_n$ that are $k$-adjacent, the equality $\Sigma_{t'}=\mu_k(\Sigma_t)$ holds.
 \end{definition}
 Fix a vertex $t_0\in \mathbb{T}_0$.
 It is clear that a cluster pattern is uniquely determined by assigning $(\mathbf{x}_{t_0}, B_{t_0})$ to $t_0\in \mathbb{T}_n$, and we refer to $t_0$ as the {\em root vertex} of this cluster pattern. For a given cluster pattern ${\bf \Sigma}=\{\Sigma_t=(\mathbf{x}_t,B_t)\}_{t\in \mathbb{T}_n}$, we always write $\mathbf{x}_t=(x_{1;t},\dots, x_{n;t})$ for the cluster and $B_t=(b_{ij;t})_{i,j=1}^n$ for the exchange matrix. We also refer to $\Sigma_{t_0}$ as the {\em initial seed}, $\mathbf{x}_{t_0}$ as the {\em initial cluster}, and the cluster variables $x_i:=x_{i;t_0}$ as the {\em initial cluster variables}.

 \begin{definition}[Cluster algebra]
 For a given cluster pattern ${\bf \Sigma}=\{\Sigma_t=(\mathbf{x}_t,B_t)\}_{t\in \mathbb{T}_n}$, the {\em cluster algebra} $\mathcal{A} = \mathcal{A}({\bf \Sigma})$ associated with ${\bf \Sigma}$ is the $\mathbb{Z}$-subalgebra of  $\mathcal{F}$ generated by all the cluster variables $\mathcal{X}=\{x_{i;t} ~|~i = 1,\dots,n; t\in \mathbb{T}_n\}$.
\end{definition}
\begin{remark}
    Up to isomorphism, a cluster algebra only depends on the mutation-equivalence class of its exchange matrices.
\end{remark}

One of the fundamental results in cluster algebras is the so-called Laurent phenomenon \cite[Theorem 3.5]{FZ02}, which states that every cluster variable $x_{i;t}\in \mathcal{A}$ can be expressed as a Laurent polynomial in the initial cluster variables. As a consequence, Fomin and Zelevinsky introduced the definition of denominator vector (d-vector for short).

\begin{definition}[d-vector]
    Let ${\bf \Sigma}$ be a cluster pattern of rank $n$. For any cluster variable $x_{i;t}$, denote by  $\mathbf{d}_{i;t}=(d_{1},\dots,d_n)$ the integer vector such that $-d_{k}$ is the lowest degree of $x_{k;t_0}$ in the Laurent polynomial expression of $x_{i;t}$ in $\mathbf{x}_{t_0}$. We call $\mathbf{d}_{i;t}$ the {\em denominator vector} of $x_{i;t}$ with respect to $t_0$.
\end{definition}
As a direct consequence of the definition, we have $\mathbf{d}_{i;t_0}=-e_i$ for $1\leq i\leq n$, where $e_1,\dots, e_n$ is the standard $\mathbb{Z}$-basis of $\mathbb{Z}^n$. In fact, the following statement holds (cf. \cite{CL20}).
\begin{lemma}\label{lem:d-vector-initial}
    A cluster variable $x$ is an initial cluster variable if and only if its denominator vector is $-e_i$ for some $1\leq i\leq n$.
\end{lemma}
The following recursion formula for d-vectors is also useful.
\begin{lemma}\label{lem:mutation-d-vector}
    Let ${\bf \Sigma}$ be a cluster pattern of rank $n$. For any pair of vertices $t,t'\in \mathbb{T}_n$ that are $k$-adjacent, the d-vectors of them satisfy the following relation:
    \[
\mathbf{d}_{k;t'}=-\mathbf{d}_{k;t}+\max\left\{ \sum_{i=1}^n  [b_{ik;t}]_+ \mathbf{d}_{i;t}, \sum_{i=1}^n [-b_{ik;t}]_+\mathbf{d}_{i;t} \right\},
\]
where $\max\{\mathbf{\alpha},\mathbf{\beta}\}=(\max\{a_1,b_1\},\dots,\max\{ a_n,b_n\})$ for $\alpha=(a_1,\dots, a_n)$ and $\beta=(b_1,\dots, b_n)$.
\end{lemma}

\subsection{Diagrams of exchange matrices}
In this subsection, we recall the diagrams of exchange matrices and their mutations, which were first introduced by Fomin and Zelevinsky in \cite{FZ03}.

\begin{definition}
Let $B=(b_{ij})^n_{i,j=1}$ be an $n \times n$ skew-symmetrizable matrix. The {\em diagram} $\Gamma(B)$ of $B$ is a weighted directed graph with vertex set $\{1,\dots, n\}$, where there is a directed edge from $i$ to $j$ if and only if $b_{ij} > 0$, and this edge is assigned the weight $\sqrt {|b_{ij}b_{ji}|}$. We say that $B$ is indecomposable if the diagram $\Gamma(B)$ is connected.
\end{definition}
We remark that our definition of the diagram is slightly different from the one in \cite{FZ03}, where the weight of the edge corresponding to $b_{ij}>0$ is $|b_{ij}b_{ji}|$, rather than $\sqrt{|b_{ij}b_{ji}|}$. The reason for this difference is that if $B$ is skew-symmetric, then $\Gamma(B)$ is the quiver corresponding to $B$ and the diagram mutation,  which will be recalled later, is the quiver mutation. On the other hand, the diagram $\Gamma(B)$ alone does not determine the matrix $B$ unless the skew-symmetrizer $D$ is also given.

\begin{definition}
For any vertex $k$ in a diagram $\Gamma:=\Gamma(B)$ of a skew-symmetrizable matrix $B$, the mutation $\mu_{k}(\Gamma)$ of $\Gamma$ at vertex $k$ is a weighted directed graph obtained from $\Gamma$ as follows: 
\begin{itemize}
\item  The orientations of all edges incident to k are reversed, their weights intact.
\item  For any vertices $i$ and $j$ which are connected in $\Gamma$ via a two-edge oriented path going through $k$ (see Figure \ref{fig:diagram-mutation}), the direction of the edge between $i, j$ in $\mu_k(\Gamma)$ and its weight $a'$ are uniquely determined by the rule
\[
\pm  a \pm  a'=bc,
\]
where the sign before $ a$ (resp. before $a'$) is $``+"$ if $i, j, k$ form an oriented cycle in $\Gamma$ (resp., in $\mu_k(\Gamma)$), and is $``-"$ otherwise. Here either $a$ or $a'$ can be equal to $0$, which means that the corresponding edge is absent.
\item  The rest of the edges and their weights in $\Gamma$ remain unchanged.
\begin{figure}
    \centering
    \begin{align*}
\xymatrix@R=0.5cm@C=0.5cm
{
&k\ar[rd]^{c}&\\
i\ar[ru]^{b}&&j\ar@{-}[ll]^a\\
}
\xymatrix@R=0.5cm@C=0.6cm
{
&\\
\ar@{-}[r]^{\mu_k}&\\
}
\xymatrix@R=0.5cm@C=0.5cm
{
&k\ar[ld]_b&\\
i\ar@{-}[rr]_{a'}&&j\ar[lu]_{c}
}
\end{align*}
    \caption{Diagram mutation.}
    \label{fig:diagram-mutation}
\end{figure}

\end{itemize}
\end{definition}
 It has been proved in \cite[Proposition 8.1]{FZ03} that the diagram mutation is compatible with the mutation of exchange matrices, i.e., $\Gamma(\mu_k(B))=\mu_k(\Gamma(B))$. Two diagrams $\Gamma$ and $\Gamma'$ are {\em mutation-equivalent} if one can be obtained from the other by a finite sequence of mutations.
A diagram $\Gamma$ of a skew-symmetrizable matrix is {\em mutation-acyclic} if it is mutation-equivalent to an acyclic diagram; otherwise, it is referred to as {\em mutation-cyclic}. For a diagram $\Gamma$, we denote by $s(\Gamma)$ the sum of weights in $\Gamma$.



The classification of  mutation-equivalence classes of rank $3$ skew-symmetrizable matrices has been obtained in \cite{ABBS08,Sev13} by their diagrams.
\begin{theorem}\cite[Theorem 1.1\&1.2]{Sev13}\label{thm:classification-rank-3-matrix}
Let $\mathcal{M}$ be a mutation-equivalence class of diagrams with $3$ vertices.  Then there is a diagram $\Gamma_0$ in $\mathcal{M}$ such that $s(\Gamma_0)$ is minimal. Furthermore,
\begin{itemize}
    \item[(1)] If $\mathcal{M}$ is mutation-cyclic, then $\Gamma_0$ can be described as follows up to an enumeration of vertices:
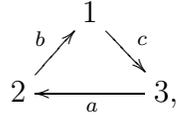
\[
\xymatrix@R=0.5cm@C=0.5cm
{
&1\ar[rd]^c&\\
2\ar[ru]^b&&3,\ar[ll]^a\\
}
\]
where $2\le a\le b\le c$ and $ab\ge 2c$.
\item[(2)] If $\mathcal{M}$ is mutation-acyclic, then  $\Gamma_0$ can be described as follows up to an enumeration of vertices:

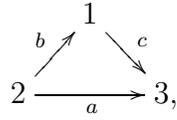
\[
\xymatrix@R=0.5cm@C=0.5cm
{
&1\ar[rd]^c&\\
2\ar[ru]^b\ar[rr]_a&&3,\\
}
\]
where $0\le a\le b\le c$.
\end{itemize}

\end{theorem}
\begin{definition}
Let $\Gamma$ be the cyclic diagram  in Figure \ref{fig:cyclic-diagram}.
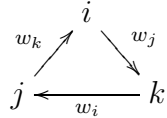
\begin{figure}[ht]
    \centering
   \[ \xymatrix@R=0.5cm@C=0.5cm
{
&i\ar[rd]^{w_j}&\\
j\ar[ru]^{w_k}&&k\ar[ll]^{w_i}\\
}
\]
    \caption{Cyclic weighted diagram.}
    \label{fig:cyclic-diagram}
\end{figure}
We say that the diagram $\Gamma$ has the {\em property $(M)$ at vertex $i$} if one of the following conditions is satisfied:
\begin{itemize}
    \item[(i)] $w_i\geq \max\{w_j,w_k\}\geq \min\{w_j,w_k\}> 2$;
    \item[(ii)] $w_i> \max\{w_j,w_k\}\geq \min\{w_j,w_k\}\geq 2$.
\end{itemize}
\end{definition}

\begin{lemma}\label{lem:weight-increasing}
    Let $\Gamma$ be a cyclic diagram as shown in Figure \ref{fig:cyclic-diagram} . If $\Gamma$ has the property $(M)$ at vertex $i$, then $\mu_j(\Gamma)$ and $\mu_k(\Gamma)$ have the property $(M)$ at $j$ and $k$, respectively. Furthermore, $s(\mu_j(\Gamma))>s(\Gamma)$ and $s(\mu_k(\Gamma))>s(\Gamma)$.
\end{lemma}
\begin{proof}
We provide the proof for $\mu_j(\Gamma)$, and the proof for $\mu_k(\Gamma)$ is similar.
    Note that $w_iw_k-w_j>0$, which implies that $\mu_j(\Gamma)$ is a cyclic weighted diagram, and the weight of the edge not incident to $j$ is $w_iw_k-w_j$.  On the other hand, 
        $w_iw_k-w_j-w_i=w_i(w_k-1)-w_j>0$,
    which implies that  $w_iw_k-w_j> \max\{w_i,w_k\}=w_i$. Hence, $\mu_j(\Gamma)$ has the property $(M)$ at vertex $j$.
\end{proof}

We conclude this subsection by introducing the $S_n$-actions on labeled seeds and diagrams.
Let $S_n$ denote the symmetric group on $n$ letters.
There are $S_n$-actions on labeled seeds and  diagrams of skew-symmetrizable matrices, respectively.
\begin{definition}[$S_n$-action]
Let $\Sigma=(\mathbf x,B)$ be a labeled seed of rank $n$ and $\Gamma(B)$ the diagram of $B$. Let $\sigma\in S_n$.
\begin{itemize}
    \item We define the action of $\sigma$ on $\Sigma$ by
\[
\sigma\Sigma=(\sigma \mathbf x,\sigma B),
\]
where $\sigma \mathbf{x}:=(x_1',\dots,x_n'), \sigma B:=(b_{ij}')_{i,j=1}^n$ are defined by
\[
x'_i=x_{\sigma^{-1}(i)}\quad \text{and}\quad b'_{ij}=b_{\sigma^{-1}(i)\sigma^{-1}(j)}.
\]
\item We define $\sigma(\Gamma(B))$ to be the diagram with vertex set $\{1,\dots,n\}$ such that the direction and weight of the edge between $i,j$ is the same as the direction and weight of the edge between $\sigma^{-1}(i),\sigma^{-1}(j)$ in the diagram $\Gamma(B)$.
\end{itemize}
\end{definition}
The following is straightforward.
\begin{lemma}\label{lem:compatible-mutation-s-action}
    Let $\Sigma=(\mathbf{x},B)$ be a labeled seeds of rank $n$ and $\Gamma(B)$ the diagram of $B$, and $1\leq k\leq n$. The following statements hold:
    \begin{itemize}
    \item[(1)] $\sigma(\Gamma(B))=\Gamma(\sigma(B))$;
        \item[(2)] $\sigma(\mu_k(\Sigma))=\mu_{\sigma(k)}(\sigma(\Sigma))$ and $\sigma(\mu_k(\Gamma(B)))=\mu_{\sigma(k)}(\sigma(\Gamma(B)))$ for any vertex $k$ of $\Gamma(B)$.
    \end{itemize}
\end{lemma}

The $S_n$-action on labeled seeds induces an equivalence relation on labeled seeds. Specifically, 
two labeled seeds $\Sigma$ and $\Sigma'$ are {\em equivalent}, denote by $\Sigma\sim \Sigma'$, if there exists a permutation $\sigma\in S_n$ such that $\sigma(\Sigma)=\Sigma'$. The equivalence class $[\Sigma]$ of $\Sigma$ is referred to as an {\em unlabeled seed}. For skew-symmetrizable matrices $B$ and $B'$, we also write $B\sim B'$ (resp. $\Gamma(B)\sim\Gamma(B')$) if there exists a $\sigma\in S_n$ such that  $B'=\sigma(B)$ (resp. $\Gamma(B')\sim\sigma\Gamma(B)$).


\section{Pseudo grading on cluster automorphism groups} \label{s:cluster-automorphism}

In this section, we recall the definition of cluster automorphism \cite{ASS2012} and introduce a pseudo $\mathbb{N}$-grading  for cluster automorphism groups. As a consequence, we prove that the cluster automorphism group is generated by its $0$-th  and $1$-st components.

Throughout this section, let $\mathcal{A}$ be a cluster algebra of rank $n$ with a fixed cluster pattern $\mathbf{\Sigma}$. Assume that the exchange matrices of $\mathcal{A}$ are indecomposable and fix a skew-symmetrizer $D=\operatorname{diag}\{d_1,\dots, d_n\}$.

\subsection{Cluster automorphism}
\begin{definition}\label{def:cluster-automorphism}
 An automorphism $f:\mathcal A \to \mathcal A$ of $\mathbb Z$-algebra is called a {\em cluster automorphism} if there exist two labeled seeds $\Sigma:=(\mathbf x ,B)$, $\Sigma':=(\mathbf x',B')$ in $\mathbf{\Sigma}$ and a permutation $\sigma \in S_n$, such that the following conditions are satisfied:
\begin{itemize}
    \item[(i)] $\sigma f(\mathbf{x})=\mathbf{x}'$;
    \item[(ii)] $\sigma(f(\mu_k(\mathbf{x})))=\mu_{\sigma(k)}(\mathbf{x}')$ for $ k=1,\dots,n$.
\end{itemize} 
\end{definition}
\begin{remark}
It was conjecture by Chang and Schiffler and proved in \cite{CLLP19} that the condition $(ii)$ is superfluous.
\end{remark}

The following is a direct consequence of Definition \ref{def:cluster-automorphism}, cf. \cite[Lemma 2.3]{ASS2012}.
\begin{corollary}\label{cor:property-cluster-automorpshim}
    Let $\Sigma=(\mathbf{x},B)$ and $\Sigma'=(\mathbf{x}',B')$ be labeled seeds of $\mathcal{A}$ and $\sigma\in S_n$ such that $\sigma f(\mathbf{x})=\mathbf{x}'$. 
    Then $\sigma B=\pm B'$, equivalently, $\sigma \Gamma(B)=\Gamma(B')$ or $\sigma\Gamma(B)=\Gamma(B')^{\text{op}}$ and $d_i=d_{\sigma^{-1}(i)}$ for $1\leq i\leq n$, where $\Gamma(B')^{\text{op}}$ is the opposite direct graph of $\Gamma(B')$.
\end{corollary}

On the other hand, we have the following.
\begin{lemma}\label{lem:matrix-determine-cluster-automorphism}
    Let $\Sigma:=(\mathbf{x},B)$ and $\Sigma':=(\mathbf{x}',B')$ be two labeled seeds of $\mathcal{A}$. Assume that there is a permutation $\sigma\in S_n$ such that $\sigma B=\pm B'$. Then there is a unique cluster automorphism $f:\mathcal{A}\to \mathcal{A}$ such that $\sigma f(\mathbf{x})=\mathbf{x}'$.
\end{lemma}
\begin{proof}
Denote by $\mathbf{x}=(x_1,\dots,x_n)$ and $\mathbf{x}'=(x_1',\dots,x_n')$.
    Note that $\sigma f(\mathbf{x})=\mathbf{x}'$ if and only if $f(x_i)=x_{\sigma(i)}'$ for $1\leq i\leq n$.
    Since each cluster is a free generating set of $\mathcal{F}$, it follows that there is a unique isomorphism of algebras $f:\mathcal{F}\to \mathcal{F}$ such that $f(x_i)=x_{\sigma(i)}'$ for $1\leq i\leq n$. By Lemma \ref{lem:compatible-mutation-s-action}, for any $1\leq k\leq n$, we have \[\sigma \mu_k(B)=\mu_{\sigma(k)}(\sigma B)=\mu_{\sigma(k)}(\pm B')=\pm \mu_{\sigma(k)}(B').\]
    Inductively, one can show that $f$ maps a cluster variable to a cluster variable, and every cluster variable has a preimage. Hence, $f$ restricts to a cluster automorphism $f:\mathcal{A}\to \mathcal{A}$.
\end{proof}

As a consequence of Corollary \ref{cor:property-cluster-automorpshim} and Lemma \ref{lem:matrix-determine-cluster-automorphism}, every cluster automorphism of $\mathcal{A}$ can be written as a quadruple
\[
    f=(\Sigma_1,\Sigma_2,\sigma, \varepsilon),
\]
where $\Sigma_1=(\mathbf{x}_1,B_1)$ and $\Sigma_2=(\mathbf{x}_2,B_2)$ are labeled seeds, $\sigma\in S_n$, and $\varepsilon\in\{\pm\}$ such that $\sigma B_1=\varepsilon B_2$. The automorphism $f$ is uniquely determined by the condition $\sigma f(\mathbf{x}_1)=\mathbf{x}_2$.
Note that the representation of a cluster auotmorphism as a quadruple is not unique. In fact, we have the following result.

\begin{proposition}\label{prop:expression-mutation-composition}
    Let $f=(\Sigma_1,\Sigma_2,\sigma, \varepsilon)$ be a cluster automorphism of $\mathcal{A}$.
    \begin{enumerate}
        \item For any mutation sequence $\mu_{i_k},\dots, \mu_{i_1}$, the following equality holds
        \[
(\Sigma_1,\Sigma_2, \sigma ,  \varepsilon)=(\mu_{i_k}\cdots\mu_{i_1}(\Sigma_1), \mu_{\sigma(i_k)}\cdots\mu_{\sigma(i_1)}(\Sigma_2), \sigma, \varepsilon).
\]
\item Let $g=(\Sigma_2,\Sigma_3,\tau, \varepsilon')$ be another cluster automorphism of $\mathcal{A}$. We have
\[
(\Sigma_2,\Sigma_3, \tau ,  \varepsilon') \circ (\Sigma_1,\Sigma_2, \sigma ,  \varepsilon)=(\Sigma_1,\Sigma_3, \tau\sigma , \varepsilon'\varepsilon).
\]
    \end{enumerate}
\end{proposition}
\begin{proof}
Denote by $\Sigma_1=(\mathbf{x}_1,B_1)$ and $\Sigma_2=(\mathbf{x}_2,B_2)$.
By Lemma \ref{lem:compatible-mutation-s-action}, we have 
\begin{align*}
    \sigma (\mu_{i_k}\dots\mu_{i_1}(B_1))&=  \mu_{\sigma(i_k)}\dots\mu_{\sigma(i_1)}(\sigma(B_1))\\
    &= \mu_{\sigma(i_k)}\dots\mu_{\sigma(i_1)}(\varepsilon B_2)\\
    &=\varepsilon  \mu_{\sigma(i_k)}\dots\mu_{\sigma(i_1)}(B_2).
\end{align*}
In particular, $(\mu_{i_k}\dots\mu_{i_1}(\Sigma_1), \mu_{\sigma(i_k)}\dots\mu_{\sigma(i_1)}(\Sigma_2), \sigma, \varepsilon)$ is a cluster auotmorphism. It remains to show that \[
\sigma f(\mu_{i_k}\cdots \mu_{i_1}(\mathbf{x}_1))=\mu_{\sigma(i_k)}\cdots \mu_{\sigma(i_1)}(\mathbf{x}_2),
\]
 which is a consequence of Definition \ref{def:cluster-automorphism}. This finishes the proof of $(1)$.

For $(2)$, let us denote by $\Sigma_3=(\mathbf{x}_3,B_3)$. We first note that
\[
\tau\sigma(B_1)=\tau(\varepsilon B_2)=\varepsilon \tau(B_2)=\varepsilon\varepsilon' B_3.
\]
This implies that $(\Sigma_1,\Sigma_3, \tau\sigma , \varepsilon'\varepsilon)$ is a cluster automorphism.
On the other hand, 
\[
    \tau\circ \sigma(g\circ f(\mathbf{x}_1))=\tau(g(\sigma f(\mathbf{x}_1)))=\tau g(\mathbf{x}_2)=\mathbf{x}_3.
\]
It follows that $g\circ f=(\Sigma_1,\Sigma_3, \tau\sigma , \varepsilon'\varepsilon)$.
This completes the proof of $(2)$.
\end{proof}
As a consequence of Proposition \ref{prop:expression-mutation-composition}, the set $\operatorname{Aut}(\mathcal{A})$ of all cluster automorphisms of $\mathcal{A}$ forms a group under composition. We refer to $\operatorname{Aut}(\mathcal{A})$ as the {\em cluster automorphism group} of $\mathcal{A}$, cf. \cite[Lemma 2.8]{ASS2012}.

\begin{remark}\label{rem:abberivation-morphism}
   Let $\Sigma_{t_0}=(\mathbf{x}_{t_0},B_{t_0})$ be the initial seed of the cluster pattern $\mathbf{\Sigma}$. 
As a consequence of Proposition \ref{prop:expression-mutation-composition} (1), every cluster automorphism can be expressed as
\[
    (\Sigma_{t_0}, \Sigma_t,\sigma, \varepsilon)
\]
for some labeled seed $\Sigma_t=(\mathbf{x}_t,B_t)$ such that $\sigma B_{t_0}=\varepsilon B_t$. Note that there is a unique reduced sequence $i_k\cdots i_1$ such that $t=\mu_{i_k}\cdots \mu_{i_1}(t_0)$. In the following, we will also denote the cluster automorphism as \[
g_{i_k\cdots i_1}^{\sigma, \varepsilon}:=(\Sigma_{t_0}, \Sigma_t,\sigma, \varepsilon).\] 
If $t=t_0$, we further denote the cluster automorphism by
\[
    \psi_{\sigma}^\varepsilon:=(\Sigma_{t_0}, \Sigma_{t_0},\sigma, \varepsilon).
    \]
\end{remark}


\subsection{Pseudo $\mathbb{N}$-grading on $\operatorname{Aut}(\mathcal{A})$}

Recall that we have fixed a cluster pattern $\mathbf{\Sigma}=\{\Sigma_t=(\mathbf{x}_t,B_t)\}_{t\in \mathbb{T}_n}$ with root vertex $t_0$. 
For any two vertices $t,t'\in \mathbb{T}_n$, let $p(t,t')$ be the unique path in $\mathbb{T}_n$ from $t$ to $t'$. Assume that $p(t,t')$ is as follows:
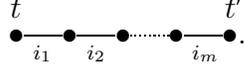
\[
\xymatrix@R=0.1cm@C=0.5cm@M=0cm
{
t &&&  &t' \\
\bullet \ar@{-}[r]_{i_1}&\bullet \ar@{-}[r]_{i_2} &\bullet \ar@{.}[r]&\bullet\ar@{-}[r]_{i_m}&\bullet.
}
\]
We also define $\mu_{p(t,t')}:=\mu_{i_m}\cdots \mu_{i_1}$ and $\mu_{\sigma(p(t,t'))}:=\mu_{\sigma(i_m)}\cdots \mu_{\sigma(i_1)}$ for any $\sigma\in S_n$.
For a fixed vertex $s \in \mathbb{T}_n$, we define
\[
    \mathbf{w}_{s}(t,t'):=\# \{t''\in \mathbb{T}_n~|~\text{the path $p(t,t')$ passes through the vertex $t''$ and $B_{t''}\sim \pm B_{s}$}\},
\]
which is called the {\em weight of $p(t,t')$ with respect to $s$}.
\begin{lemma}\label{lem:factorization-cluster-automorphism}
    Let $s,t\in \mathbb{T}_n$ such that $B_{t_0}\sim \pm B_s$ and $B_{t_0}\sim \pm B_t$. Let $f=(\Sigma_{t_0},\Sigma_t,\sigma,\varepsilon),g=(\Sigma_{t_0},\Sigma_s,\tau,\delta)\in \operatorname{Aut}(\mathcal{A})$. Denote by \[h=(\Sigma_{t_0},\mu_{\tau^{-1}(p(s,t))}(\Sigma_{t_0}), \tau^{-1}\sigma, \varepsilon\delta).\] Then $h$ is a cluster automorphism of $\mathcal{A}$ and $f=g\circ h$. Furthermore, for any vertex $t'\in \mathbb{T}_n$, we have $\mathbf{w}_{t'}(t_0,\mu_{\tau^{-1}(p(s,t))}(t_0))=\mathbf{w}_{t'}(s,t)$.
\end{lemma}
\begin{proof}
    By Proposition \ref{prop:expression-mutation-composition}, we have
    \begin{align*}
        g^{-1}&=(\Sigma_s,\Sigma_{t_0},\tau^{-1},\delta)\\
        &=(\mu_{p(s,t)}(\Sigma_s),\mu_{\tau^{-1}(p(s,t))}(\Sigma_{t_0}), \tau^{-1},\delta)\\
        &=(\Sigma_t,\mu_{\tau^{-1}(p(s,t))}(\Sigma_{t_0}), \tau^{-1},\delta).
    \end{align*}  Hence,
    \begin{align*}
        g^{-1}\circ f&=(\Sigma_t,\mu_{\tau^{-1}(p(s,t))}(\Sigma_{t_0}), \tau^{-1},\delta)\circ (\Sigma_{t_0},\Sigma_t,\sigma, \varepsilon)\\
        &=(\Sigma_{t_0},\mu_{\tau^{-1}(p(s,t))}(\Sigma_{t_0}),\tau^{-1}\sigma,\delta\varepsilon)=h.
    \end{align*}
   In particular, $h\in \operatorname{Aut}(\mathcal{A})$ and $f=g\circ h$. 

Since $\tau B_{t_0}=\pm B_s$, for any $1\leq k\leq n$, we have $\tau(\mu_k(B_{t_0}))=\pm\mu_{\tau(k)}(B_s)$ by Lemma \ref{lem:compatible-mutation-s-action}.
It follows that 
\[
    \mathbf{w}_{t'}(t_0,\mu_{\tau^{-1}(p(s,t))}(t_0))=\mathbf{w}_{t'}(s,\mu_{p(s,t)}(s))=\mathbf{w}_{t'}(s,t).
\]

\end{proof}

For any nonnegative integer $m$, we define 
\[
    \mathcal{P}_m(t_0)=\{p(t_0,t)~\mid~\text{$\mathbf{w}_{t_0}(t_0,t)=m+1$ and $B_t\sim \pm B_{t_0}$}\}.
\]
In particular, every path $p(t_0,t)$ in $\mathcal{P}_m(t_0)$ determines at least one cluster automorphism $f=(\Sigma_{t_0},\Sigma_t,\sigma, \varepsilon)$ for some $\sigma\in S_n$ and $\varepsilon\in \{\pm\}$.
Note that $\mathcal{P}_0(t_0)$ consists only of the trivial path $p(t_0,t_0)$.

Denote by 
\[
    G_m(t_0):=\left\{f\in \operatorname{Aut}(\mathcal{A})~\middle|~\begin{array}{l} \text{$\exists s\in \mathbb{T}_n$, $\sigma\in S_n$ and $ \varepsilon\in \{\pm\}$} \\ \text{such that $\mathbf{w}_{{t_0}}(t_0,s)=m+1$ and}  \\ \text{ $f=(\Sigma_{t_0},\Sigma_s,\sigma, \varepsilon)$} \end{array}
    \right\}\subset \operatorname{Aut}(\mathcal{A}).
\]
As a consequence of the following proposition, we refer to the family $\{G_m\}_{m\in \mathbb{N}}$ as a {\em pseudo $\mathbb{N}$-grading} on $\operatorname{Aut}(\mathcal{A})$.
\begin{proposition}
    The following statements hold:
    \begin{itemize}
        \item[(1)] $\operatorname{Aut}(\mathcal{A})=\bigcup_{m\in \mathbb{N}}G_m(t_0)$;
        \item[(2)] $G_0(t_0)$ is a subgroup of $\operatorname{Aut}(\mathcal{A})$;
        \item[(3)] For any $f\in G_k(t_0),g\in G_l(t_0)$, we have $g\circ f\in \bigcup_{m=0}^{k+l}G_m(t_0)$.
    \end{itemize}
\end{proposition}
\begin{proof}
    The statement $(1)$ follows directly from the definition. The statements $(2)$ and $(3)$ are  consequences of Proposition \ref{prop:expression-mutation-composition}.
\end{proof}

\begin{remark}
	It is clear that $G_0(t_0)$ is isomorphic to the subgroup of $S_n$ given by
	\[\{\sigma\in S_n~|~\sigma(B_{t_0})=\pm B_{t_0}\}.
	\]
	Thus, we refer to $G_0(t_0)$ as the group of  symmetries of the exchange matrix $B_{t_0}$.
\end{remark}

\subsection{Generating sets of $\operatorname{Aut}(\mathcal{A})$}\label{ss:generator-set}
For each path $p(t_0,t)\in \mathcal{P}_m(t_0)$, we fix a permutation $\sigma_t\in S_n$ and $\varepsilon_t\in \{\pm\}$ such that $\sigma_t(B_{t_0})=\varepsilon_tB_t$, and denote by $f_{p(t_0,t)}:=(\Sigma_{t_0},\Sigma_t,\sigma_t,\varepsilon_t)$ the associated cluster automorphism. 
\begin{remark}\label{rem:factorization-trivial-path}
   Let $f=(\Sigma_{t_0},\Sigma_t,\sigma, \varepsilon)$  be a cluster automorphism detemined by the path $p(t_0,t)$. According to Lemma \ref{lem:factorization-cluster-automorphism}, we have the following facterization \[f=f_{p(t_0,t)}\circ h,\] 
   where $h=(\Sigma_{t_0},\Sigma_{t_0},\sigma_t^{-1}\sigma, \varepsilon_t\varepsilon)$. In particular, every cluster automorphism determied by $p(t_0,t)$ is a product of $f_{p(t_0,t)}$ with a cluster automorphism associated with the trivial path $p(t_0,t_0)$.
   
\end{remark}

For any positive integer $m$, define 
\[H_m(t_0):=\{f_{p(t_0,t)}~|~p(t_0,t)\in \mathcal{P}_m(t_0)\}\subset G_m(t_0).\]

\begin{theorem}\label{thm:generator-set}
    The cluster automorphism group $\operatorname{Aut}(\mathcal{A})$ is generated by $G_0(t_0)\cup H_1(t_0)$.
\end{theorem}
\begin{proof}
According to Remark \ref{rem:factorization-trivial-path}, $G_1(t_0)$ is generated by $G_0(t_0)\cup H_1(t_0)$. It remains to prove that $\operatorname{Aut}(\mathcal{A})$ is generated by $G_0(t_0)\cup G_1(t_0)$. Note that $\operatorname{Aut}(\mathcal{A})=\bigcup_{m\geq 0}G_m(t_0)$. It suffices to show that $G_m(t_0)$ is generated by $G_0(t_0)\cup G_1(t_0)$ for any $m\geq 2$.

Suppose that $G_m(t_0)$ is generated by $G_0(t_0)\cup G_1(t_0)$, and let $f\in G_{m+1}(t_0)$. Without loss of generality, assume that $f=(\Sigma_{t_0},\Sigma_t,\sigma,\varepsilon)$ for some $\sigma\in S_n$, $\varepsilon\in \{\pm\}$ and $\mathbf{w}_{t_0}(t_0,t)=m+2$. Let $s$ be the vertex on the path $p(t_0,t)$ such that $B_{t_0}\sim \pm B_{s}$, $\mathbf{w}_{t_0}(t_0,s)=2$ and $\mathbf{w}_{t_0}(s,t)=m+1$. Let $f'$ be a cluster automorphism determined by the path $p(t_0,s)$. It follows that $f'\in G_1(t_0)$.
By Lemma \ref{lem:factorization-cluster-automorphism}, there exists $f''\in G_m(t_0)$ such that $f=f'\circ f''$. Hence, $f$ is generated by $G_0(t_0)\cup G_1(t_0)$. By induction, we conclude that $\operatorname{Aut}(\mathcal{A})$ is generated by $G_0(t_0)\cup G_1(t_0)$.
\end{proof}

The generator set $G_0(t_0)\cup H_1(t_0)$ can be further refined in many cases.
 In the following, we assume that $\mathcal{A}$ satisfies the following assumption:

\noindent{\bf Assumption ($\diamond$):}
there exists a vertex  $t_{\diamond}\in \mathbb{T}_n$ such that $B_{t_{\diamond}}\not\sim \pm B_{t_0}$ and $\mathbf{w}_{t_0}(t_0,t_{\diamond})=\mathbf{w}_{t_{\diamond}}(t_0,t_{\diamond})=1$. 

\begin{remark}
	Let $\mathcal{B}$ be the mutation-equivalence class of $B_{t_0}$. The assumption ($\diamond$) is satisfied if there exists a matrix $B'\in \mathcal{B}$ such that $B'\not\sim \pm B_{t_0}$. Equivalently, there exists a $k$ such that $\mu_k(B_{t_0})\not\sim B_{t_0}$. However, the choice of $t_\diamond$ is not unique in general. 
 \end{remark}
 
 In the following, we fix a vertex $t_\diamond$ satisfying assumption ($\diamond$).

For any $m\geq 0$, we define 
\[
    \mathcal{P}_{1,m}^{t_\diamond}(t_0)=\{p(t_0,t)\in \mathcal{P}_1(t_0)~\mid~\mathbf{w}_{t_\diamond}(t_0,t)=m\}\subset \mathcal{P}_1(t_0),
\]
\[
    H_{1,m}^{t_\diamond}(t_0):=\{f_{p(t_0,t)}~\mid~p(t_0,t)\in \mathcal{P}_{1,m}^{t_\diamond}(t_0)\}\subset H_1(t_0),
\]
\[
     G_{1,m}^{t_\diamond}(t_0)=\left\{f\in \operatorname{Aut}(\mathcal{A})~\middle|~ \begin{array}{l}\text{$\exists\quad p(t_0,t)\in \mathcal{P}_{1,m}^{t_\diamond}(t_0)$, $\sigma\in S_n$ and $\varepsilon\in \{\pm\}$}\\ \text{ such that $f=(\Sigma_{t_0},\Sigma_t,\sigma, \varepsilon)$}\end{array}\right\}\subset G_1(t_0).
 \]
It is clear that \[\mathcal{P}_1(t_0)=\bigcup\limits_{m\geq 0}\mathcal{P}_{1,m}^{t_\diamond}(t_0),\quad H_1(t_0)=\bigcup\limits_{m\geq 0}H_{1,m}^{t_\diamond}(t_0),\quad G_1(t_0)=\bigcup\limits_{m\geq 0}G_{1,m}^{t_\diamond}(t_0).\]
Finally, define 
\[
    K_1^{t_\diamond}(t_0):=\{f_{p(t_0,t)}~\mid~\text{$p(t_0,t)\in \mathcal{P}_{1,2}^{t_\diamond}(t_0)$ and passes through $t_\diamond$}\}\subset H_{1,2}^{t_\diamond}(t_0).
\]
\begin{proposition}\label{prop:generator-set-assumption}
    The cluster automorphism group $\operatorname{Aut}(\mathcal{A})$ is generated by $G_0(t_0)\cup H_{1,0}^{t_\diamond}(t_0)\cup H_{1,1}^{t_\diamond}(t_0)\cup K_1^{t_\diamond}(t_0)$.
\end{proposition}
\begin{proof}
Let us first prove that $\operatorname{Aut}(\mathcal{A})$ is generated by $G_0(t_0)\cup H_{1,0}^{t_\diamond}(t_0)\cup H_{1,1}^{t_\diamond}(t_0)\cup H_{1,2}^{t_\diamond}(t_0)$. According to Theorem \ref{thm:generator-set}, it suffices to show that $H_1(t_0)$ is generated by $G_0(t_0)\cup H_{1,0}^{t_\diamond}(t_0)\cup H_{1,1}^{t_\diamond}(t_0)\cup H_{1,2}^{t_\diamond}(t_0)$. Note that $H_1(t_0)=\bigcup_{m\geq 0}H_{1,m}^{t_\diamond}(t_0)$. Therefore, it suffices to show that for every $m\geq 0$,
 $H_{1,m}^{t_\diamond}(t_0)$ is generated by $G_0(t_0)\cup H_{1,0}^{t_\diamond}(t_0)\cup H_{1,1}^{t_\diamond}(t_0)\cup H_{1,2}^{t_\diamond}(t_0)$. Assume this holds for $m\geq 2$, and let $f=f_{p(t_0,t)}\in H_{1,m+1}^{t_\diamond}(t_0)$, where $p(t_0,t)\in \mathcal{P}_{1,m+1}^{t_\diamond}(t_0)$. Recall that there exist $\sigma_t\in S_n$ and $\varepsilon_t\in \{\pm\}$ such that $f=(\Sigma_{t_0},\Sigma_{t},\sigma_t, \varepsilon_t)$. Let $t_1$ be the vertex lying on $p(t_0,t)$ such that $B_{t_1}\sim \pm B_{t_\diamond}$ and
\[
    \mathbf{w}_{t_0}(t_0,t_1)=1,  \mathbf{w}_{t_\diamond}(t_0.t_1)=2,
    \mathbf{w}_{t_0}(t_1,t)=1,  \mathbf{w}_{t_\diamond}(t_1,t)=m.
\]
By assumption ($\diamond$), there exists a vertex $s\in \mathbb{T}_n$ such that $B_s\sim \pm B_{t_0}$ and $\mathbf{w}_{t_0}(t_1,s)=1=\mathbf{w}_{t_\diamond}(t_1,s)$.  Specifically, suppose  that $\tau(B_{t_\diamond})=\pm B_{t_1}$, then $s=\mu_{\tau(p(t_\diamond,t_0))}(t_1)$.
Consequently, we have
\[
    \mathbf{w}_{t_0}(t_0,s)=2, 1\leq \mathbf{w}_{t_\diamond}(t_0,s)\leq 2, \mathbf{w}_{t_0}(s,t)=2, m-1\leq \mathbf{w}_{t_\diamond}(s,t)\leq m.
\]
Applying Lemma \ref{lem:factorization-cluster-automorphism}, we conclude that $f$ is the product of an element from $G_{1,1}^{t_\diamond}(t_0)\cup G_{1,2}^{t_\diamond}(t_0)$  and an elements from $G_{1,m-1}^{t_\diamond}\cup G_{1,m}^{t_\diamond}(t_0)$. By induction and Remark \ref{rem:factorization-trivial-path}, we conclude that $H_1(t_0)$ is generated by $G_0(t_0)\cup H_{1,0}^{t_\diamond}(t_0)\cup H_{1,1}^{t_\diamond}(t_0)\cup H_{1,2}^{t_\diamond}(t_0)$.

It remains to show that $H_{1,2}^{t_\diamond}(t_0)$ is generated by $G_0(t_0)\cup H_{1,0}^{t_\diamond}(t_0)\cup H_{1,1}^{t_\diamond}(t_0)\cup K_1^{t_\diamond}(t_0)$. Let $g=f_{p(t_0,t_2)}\in H_{1,2}^{t_\diamond}(t_0)$, where $p(t_0,t_2)\in \mathcal{P}_{1,2}^{t_\diamond}(t_0)$. Suppose that $\sigma_{t_2}\in S_n$ and $\varepsilon_{t_2}\in \{\pm\}$ such that $g=(\Sigma_{t_0},\Sigma_{t_2},\sigma_{t_2},\varepsilon_{t_2})$.

Let $t_3$ be the vertex lying on $p(t_0,t_2)$ such that $B_{t_3}\sim \pm B_{t_\diamond}$ and \[\mathbf{w}_{t_0}(t_0,t_3)=\mathbf{w}_{t_\diamond}(t_0,t_3)=\mathbf{w}_{t_0}(t_3,t_2)=1, \mathbf{w}_{t_\diamond}(t_3,t_2)=2.\]
Assume that $B_{t_3}=\varepsilon_3\sigma_3 B_{t_\diamond}$ for some $\sigma_3\in S_n$ and $\varepsilon_3\in \{\pm\}$. Denote by $s'=\mu_{\sigma_3(p(t_\diamond,t_0))}(t_3)$. It follows that $B_{s'}=\varepsilon_3 \sigma_3 B_{t_0}$ and $\mathbf{w}_{t_0}(t_3,s')=\mathbf{w}_{t_\diamond}(t_3,s')=1$. Let $g'=(\Sigma_{t_0},\Sigma_{s'},\sigma_3,\varepsilon_3)$, so that $g'\in  G_0(t_0)\cup G_{1,0}^{t_\diamond}(t_0)\cup G_{1,1}^{t_\diamond}(t_0)$. By applying Lemma \ref{lem:factorization-cluster-automorphism}, we obtain  $g=g'\circ g''$, where $g''=(\Sigma_{t_0},\mu_{\sigma_3^{-1}(p(s',t_2))}(\Sigma_{t_0}),\sigma_3^{-1}\sigma_{t_2},\varepsilon_3\varepsilon_{t_2})$. In particular, $g''\in G_{1,1}^{t_\diamond}(t_0)\cup G_{1,2}^{t_\diamond}(t_0)$ and
the path $p(t_0,\mu_{\sigma_3^{-1}(p(s',t_2))}(t_0))$ passes through the vertex $t_\diamond$ if $g''\in G_{1,2}^{t_\diamond}(t_0)$. By Remark \ref{rem:factorization-trivial-path}, we conclude that $g$ is generated by $G_0(t_0)\cup H_{1,0}^{t_\diamond}(t_0)\cup H_{1,1}^{t_\diamond}(t_0)\cup K_1^{t_\diamond}(t_0)$.
This completes the proof.
\end{proof}

\begin{remark}\label{rem:sharpen-K-2}
     Assume that there are two paths $p(t_0,t)$ and $p(t_0,s)\in \mathcal{P}_{1,2}^{t_\diamond}(t_0)$ passing through $t_\diamond$, such that
    \[
        \xymatrix{t_0\ar@{-}[r]&\cdot\ar@{.}[r]&t_\diamond\ar@{.}[r]&t_2\ar@{.}[r]\ar@{.}[d]&t\\&
        &&s},
    \]
    where $ B_{t_2}\sim \pm B_{t_\diamond}$. By Lemma \ref{lem:factorization-cluster-automorphism}, $f_{p(t_0,t)}=f_{p(t_0,s)}\circ h$, where $h\in G_{1,0}^{t_\diamond}(t_0)\cup G_{1,1}^{t_\diamond}(t_0)$. In this case, the set $K_1^{t_\diamond}(t_0)$ can be further refined.
\end{remark}

\section{Application to cluster algebras of rank $3$}\label{s:cluster-automorphism-rank-3}

In this section, we present the main results concerning the cluster automorphism groups of  cluster algebras of rank $3$ with indecomposable exchange matrices. Let $\mathcal{A}:=\mathcal{A}(\Gamma_0)$ be the cluster algebra associated with the diagram $\Gamma_0$ in Theorem \ref{thm:classification-rank-3-matrix} with skew-symmetrizer $D=\operatorname{diag}\{d_1,d_2,d_3\}$. We fix a cluster pattern $\mathbf{\Sigma}=\{\Sigma_t=(\mathbf{x}_t,B_t)\}_{t\in \mathbb{T}_3}$ of $\mathcal{A}$, such that $\Gamma(B_{t_0})=\Gamma_0$.

To state the main result, we first recall some relevant groups.


\begin{itemize}
\item $\langle \mathbf{id}\rangle$: the trivial group.
\item $D_{m}$: the dihedral group of order $2m$, which has a presentation as
\[
    \langle x,y |x^m=y^2=\mathbf{id},xy=yx^{-1}\rangle.
\]
    \item $D_\infty$: the infinite dihedral group, which has presentations as \[\langle r,s~|~r^2=s^2=\mathbf{id}\rangle\cong \langle x,y~|~y^2=\mathbf{id}, xy=yx^{-1}\rangle,\] and the isomorphism is given by $r\mapsto xy$ and $s\mapsto y$.
    \item $D_\infty\rtimes_{\rho}\mathbb{Z}_2$: the semidirect product of $D_\infty$ and $\mathbb{Z}_2=\langle z|z^2=\mathbf{id}\rangle$, where $\rho:\mathbb{Z}_2\to \operatorname{Aut}(D_\infty)$ is given by $\rho(z)(x)=x^{-1}$ and $\rho(z)(y)=y$. The group $D_\infty\rtimes_{\rho}\mathbb{Z}_2$ also has  a presentation as \[\langle x,y,z|y^2=z^2=\mathbf{id}, xy=yx^{-1},xz=zx^{-1},yz=zy\rangle.\]
    \item $N:=\langle f_1,f_2,f_3|f_1^2=f_2^2=f_3^2=\mathbf{id}\rangle$.
    \item $N\rtimes_{\tau} \mathbb{Z}_2$: the semidirect product of $N$ and $\mathbb{Z}_2=\langle z|z^2=\mathbf{id}\rangle$, where $\tau:\mathbb{Z}_2\to \operatorname{Aut}(N)$ is given by $\tau(z)(f_1)=f_1$, $\tau(z)(f_2)=f_3$ and $\tau(z)(f_3)=f_2$.
    \item $N\rtimes_{\varphi} S_3$: the semidirect product of $N$ and the permutation group $S_3$, where $\varphi:S_3\to \operatorname{Aut}(N)$ is given by $\varphi(\sigma)(f_i)=f_{\sigma(i)}$ for $\sigma\in S_3$ and $1\leq i\leq 3$.
\end{itemize}

\begin{theorem}\label{thm:cyclic-case}
    Assume that $\Gamma_0$ is mutation-cyclic, then the cluster automorphism group $\operatorname{Aut}(\mathcal{A})$ is given by Table \ref{tab:mutation-cyclic}.
\begin{table}
\centering
    \begin{tabular}{|c|c|c|c|c|}
\hline
\multicolumn{3}{|c|}{$2\leq a\leq b\leq c,ab\geq 2c$} & $d_1,d_2,d_3$ & $\textup{Aut} \,(\mathcal A)$ \\
\hline
\multirow{6}{*}{$a=2$}    & \multicolumn{2}{c|}{\multirow{5}{*}{$a=b=c$}}  &  $d_1=d_2=d_3$     &   $  N\rtimes_\varphi S_3 $ \\    
\cline{4-5}
    &  \multicolumn{2}{c|}{}     &    $d_1=d_2\ne d_3$     &     \multirow{3}{*}{$ N\rtimes_\tau \mathbb{Z}_2 $} \\
\cline{4-4}
    &  \multicolumn{2}{c|}{}     &    $d_1=d_3\ne d_2$     &    \\
\cline{4-4}
    &  \multicolumn{2}{c|}{}     &    $d_2=d_3\ne d_1$     &    \\
\cline{4-5}
    &  \multicolumn{2}{c|}{}     &    $d_1\ne d_2,d_2\ne d_3,d_3\ne d_1$     &    $N$ \\ 
\cline{2-5}
   &   \multicolumn{2}{c|}{$a<b=c$}    &          &   $D_\infty$ \\ 
\hline
\multirow{12}{*}{$a>2$}    &   \multirow{3}{*}{$ab=2c$}    &      \multirow{2}{*}{$a=b<c$}    &    $d_1= d_3$     &    $\mathbb{Z}_2\oplus \mathbb{Z}_2$ \\
\cline{4-5}
    &      &      &    $d_1\ne d_3$     &    $\mathbb{Z}_2$ \\
\cline{3-5}
    &        &  $a<b<c$     &         &    $\mathbb{Z}_2$ \\
\cline{2-5}
    &   \multirow{9}{*}{$ab>2c$}    &    \multirow{5}{*}{$a=b=c$}    &    $d_1=d_2= d_3$     &    $S_3$ \\
\cline{4-5}
   &       &      &    $d_1=d_2\ne d_3$     &     \multirow{3}{*}{$ \mathbb{Z}_2 $} \\
\cline{4-4}
   &       &    &    $d_1=d_3\ne d_2$     &     \\
\cline{4-4}
   &       &       &    $d_2=d_3\neq d_1$     &     \\
\cline{4-5}
&       &       &    $d_1\neq d_2, d_2\neq d_3, d_1\neq d_3$     &    $\langle \mathbf{id}\rangle$ \\
   
\cline{3-5}
   &       &     \multirow{2}{*}{$a<b=c$}   &    $d_2=d_3$     &    $\mathbb{Z}_2$ \\
\cline{4-5}
   &       &       &    $d_2\ne d_3$     &    $\langle \mathbf{id}\rangle$ \\
\cline{3-5}
   &       &     \multirow{2}{*}{$a=b<c$}   &    $d_1=d_3$     &    $\mathbb{Z}_2$ \\
\cline{4-5}
   &       &       &    $d_1\ne d_3$     &    $\langle \mathbf{id}\rangle$ \\
 \cline{3-5}
   &       &  $a<b<c$    &        &    $\langle \mathbf{id}\rangle$ \\
\hline
\end{tabular}
\caption{Cluster auotmorphism groups of mutation-cyclic cases}
\label{tab:mutation-cyclic}
\end{table}
\end{theorem}
The proof of Theorem \ref{thm:cyclic-case} will be presented in Section \ref{s:mutation-cyclic}.

\begin{theorem}\label{thm:acyclic-case}
Assume that $\Gamma_0$ is mutation-acyclic, then the cluster automorphism group $\operatorname{Aut}(\mathcal{A})$ is given by Table \ref{tab:mutation-acyclic}.

\begin{table}[ht]
\centering
\begin{tabular}{|c|c|c|c|c|}
\hline
\multicolumn{3}{|c|}{$0\le a\le b\le c$} & $d_1,d_2,d_3$ & $\textup{Aut} \,(\mathcal A)$ \\
\hline
\multirow{5}{*}{$a=0$}    & \multirow{2}{*}{$b=1$} & $c=1$ &    &     $D_6$  \\    
\cline{3-5}
    &  & $c=\sqrt 2$   &     &   $D_4$  \\
\cline{2-5}
    & \multirow{3}{*}{$b^2+c^2\ge 4$} & \multirow{2}{*}{$b=c$}   &    $d_2\neq d_3$     &  $D_\infty$   \\
\cline{4-5}
    & &    &    $d_2= d_3$     &  $D_\infty\rtimes_\rho \mathbb{Z}_2$  \\
\cline{3-5}
    &  &  $b<c$  &         &    $D_\infty$ \\ 
\hline
\multirow{4}{*}{$a=1$}    & \multirow{2}{*}{$b=1$} & $c=1$ &    &     $D_\infty$  \\    
\cline{3-5}
    &  & $c\ge 2$   &     &   $D_\infty$  \\
\cline{2-5}
    & \multicolumn{2}{c|}{$1<b=c$}   &        &  $ D_\infty$   \\
\cline{2-5}
    & \multicolumn{2}{c|}{$1<b<c$}   &       &  $\mathbb{Z}$  \\
\hline
\multirow{3}{*}{$a=\sqrt{2}$ or $\sqrt{3}$ }    &\multicolumn{2}{c|}{\multirow{2}{*}{$a=b$}} &  $d_1=d_3$  &     $D_\infty$  \\    
\cline{4-5}
    &  \multicolumn{2}{c|}{}   &   $d_1\ne d_3$  &   $\mathbb{Z}$   \\
\cline{2-5}
    & \multicolumn{2}{c|}{$a<b$}   &       &  $\mathbb{Z}$  \\
\hline
  \multirow{9}{*}{$a\ge 2$}    &  \multicolumn{2}{c|}{\multirow{5}{*}{$a=b=c$}}      &    $d_1=d_2=d_3$     &  $D_\infty$ \\
\cline{4-5}
   &     \multicolumn{2}{c|}{}     &    $d_1=d_2\ne d_3$     &    $D_\infty$ \\
   \cline{4-5}
   &     \multicolumn{2}{c|}{}     &    $d_1=d_3\ne d_2$     &    $D_\infty$ \\
\cline{4-5}
   &      \multicolumn{2}{c|}{}   &    $d_2=d_3\ne d_1$     &    $D_\infty$ \\
\cline{4-5}
   &      \multicolumn{2}{c|}{}      &    else     &    $\mathbb Z$ \\
\cline{2-5}
   &     \multicolumn{2}{c|}{\multirow{2}{*}{$a=b<c$}}  &    $d_1=d_3$     &    $D_\infty$ \\
\cline{4-5}
   &       \multicolumn{2}{c|}{}     &  $d_1\ne d_3$  & $\mathbb Z$ \\
\cline{2-5}
   &    \multicolumn{2}{c|}{\multirow{2}{*}{$a<b=c$}} &    $d_2=d_3$     &    $D_\infty$ \\
\cline{4-5}
   &       \multicolumn{2}{c|}{}      &    $d_2\ne d_3$     &    $\mathbb Z$ \\
\cline{2-5}
   &         \multicolumn{2}{c|}{$a<b<c$}    &        &    $\mathbb Z$ \\
\hline
\end{tabular}
\caption{Cluster automorphism groups of mutation-acyclic cases}
\label{tab:mutation-acyclic}
\end{table}

\end{theorem}
The proof of Theorem \ref{thm:acyclic-case} will be given in Section \ref{s:mutation-acyclic}.

\section{Proof of Theorem \ref{thm:cyclic-case}: Mutation-cyclic cases}
\label{s:mutation-cyclic}
In this section, we assume that $\Gamma_0$ is cyclic (cf. the diagram in Theorem \ref{thm:classification-rank-3-matrix} (1)). We have $2\leq a\leq b\leq c$ and $ab\geq 2c$.
Let \[B_{t_0}=\begin{pmatrix}
    0&-b_{12}&b_{13}\\ b_{21}&0&-b_{23}\\ -b_{31}&b_{32}&0
\end{pmatrix}\] be the exchange matrix at vetex $t_0$, where $b_{ij}\in \mathbb{N}$. By definition, we have $a^2=b_{23}b_{32}$, $b^2=b_{12}b_{21}$ and $c^2=b_{13}b_{31}$.
We will calculate $\operatorname{Aut}(\mathcal{A})$ case-by-case and finish the proof of Theorem \ref{thm:cyclic-case}.

\subsection{Case 1: $a>2$, $ab>2c$.}
By definition, $\Gamma_0$ has property $(M)$ at vertex $2$.
According to Lemma \ref{lem:weight-increasing}, for any vertex $t\in \mathbb{T}_3$ such that the first edge $\xymatrix{t_0\ar@{-}[r]^k&\mu_k(t_0)}$ of the path $p(t_0,t)$ is labeled by $k\neq 2$, we have $s(\Gamma_t)>s(\Gamma_0)$.

Now consider the mutation of $\Gamma_0$ in direction $2$. Since $ab>2c$, we know that $\Gamma_{\mu_2(t_0)}=\mu_2(\Gamma_0)$ is cyclic, i.e.,
\begin{align*}
\xymatrix@R=0.5cm@C=0.5cm
{
&1\ar[rd]^{c}&\\
2\ar[ru]^{b}&&3\ar[ll]^a\\
}
\xymatrix@R=0.5cm@C=0.6cm
{
&\\
\ar@{-}[r]^{\mu_2}&\\
}
\xymatrix@R=0.5cm@C=0.5cm
{
&1\ar[ld]_b&\\
2\ar[rr]_{a}&&3.\ar[lu]_{ab-c}\\
}
\end{align*}
Moreover, since $ab-2c>0$, it follows that $s(\mu_2(\Gamma_0))>s(\Gamma_0)$, and $\mu_2(\Gamma_0)$ also has the property $M$ at vertex $2$. By Lemma \ref{lem:weight-increasing} again, we conclude that for any vertex $t\in \mathbb{T}_3$ such that the first edge $\xymatrix{\mu_2(t_0)\ar@{-}[r]^k&\mu_k\mu_2(t_0)}$ of the path $p(\mu_2(t_0),t)$ is labeled by $k\neq 2$, we have $s(\Gamma_t)>s(\Gamma_{t_0})$. 
Consequently, for any $t_0\neq t\in \mathbb{T}_3$, $s(\Gamma_t)>s(\Gamma_0)$. It follows that $\operatorname{Aut}(\mathcal{A})=G_0(t_0)$. A direct computation shows that

\begin{itemize}
    \item[(1)]  Assume that $d_1,d_2,d_3$ are pairwise distinct. Then $G_0(t_0)=\{\mathbf{id}\}$;
    \item[(2)] Assume that exactly two of $d_1,d_2,d_3$ are equal. 
    \begin{itemize}
        \item[(i)] $d_1=d_2\neq d_3$: $G_0(t_0)\neq \{\mathbf{id}\}$ if and only if $a=c$, equivalently, $a=b=c$. 
        \item[(ii)]$d_1=d_3\neq d_2$: $G_0(t_0)\neq \{\mathbf{id}\}$ if and only if $a=b$. 
        \item[(iii)] $d_2=d_3\neq d_1$: $G_0(t_0)\neq \{\mathbf{id}\}$ if and only if $b=c$. 
    \end{itemize}
    In all of these cases, we have $G_0(t_0)\cong \mathbb{Z}_2$.
    
    \item[(3)] Assume that $d_1=d_2=d_3$. It follows that $b_{12}=b_{21}=b$, $b_{13}=b_{31}=c$ and $b_{23}=b_{32}=a$.
    \begin{itemize}
        \item[(i)] If $a,b,c$ are pairwise distinct, then $G_0(t_0)=\{\mathbf{id}\}$.
        \item[(ii)] If exactly two of $a,b,c$ are equal, then either $a=b\neq c$ or $a\neq b=c$. In this case,  $G_0(t_0)\cong \mathbb{Z}_2$.
        \item[(iii)] If $a=b=c$, then $G_0(t_0)\cong S_3$.
       
    \end{itemize}
\end{itemize}

\subsection{Case 2: $a>2,ab=2c$} In this case, $b<c$.
 Similar to the case 1, one can show that $\Gamma_t\sim \Gamma_0$ or $\Gamma_t\sim \Gamma_0^{\text{op}}$ if and only if $t=t_0$ or $t=\mu_2(t_0)$, by Lemma \ref{lem:weight-increasing}. Furthermore, $\Gamma_{\mu_2(t_0)}=\Gamma_0^{\text{op}}$. It follows that $|H_1(t_0)|=1$, and we may set \[H_1(t_0)=\{g_{2}^{\mathbf{id},-}\}.\]
It is easy to see that $(g_{2}^{\mathbf{id},-})^2=\mathbf{id}$.

 On the other hand, $G_0(t_0)\neq \{\mathbf{id}\}$ if and only if $a=b$ and $d_1=d_3$. In this case, $G_0(t_0)\cong \mathbb{Z}_2$ and $ \psi_{(13)}^{-}$ is a generator of $G_0(t_0)$. Moreover,  $g_{2}^{\mathbf{id},-}\circ \psi_{(13)}^-=\psi_{(13)}^-\circ g_{2}^{\mathbf{id},-}$, since $b_{12}=b_{32}$. According to  Theorem \ref{thm:generator-set}, we have
 \[
     \operatorname{Aut}(\mathcal{A})\cong\begin{cases}
         \mathbb{Z}_2\oplus \mathbb{Z}_2 & \text{$a=b$ and $d_1=d_3$;}\\ \mathbb{Z}_2& \text{else.}
     \end{cases}
 \]

\subsection{Case 3: $a=2$.} Since $ab\geq 2c$ and $a\leq b\leq c$, we obtain $b=c$. 

\subsubsection{\bf $a=b=c=2$.}\label{sss:a=b=c=2}
By mutating $\Gamma_0$ in three directions, we obtain 
{\tiny
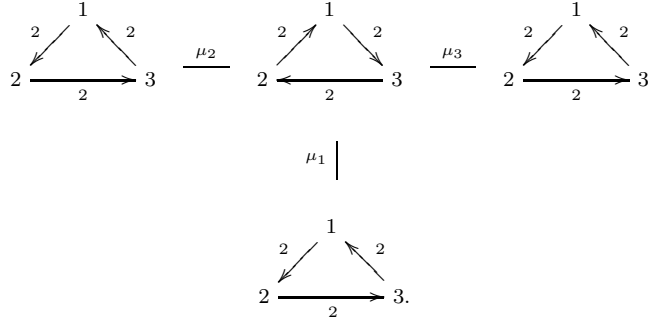
\begin{align*}
&
\xymatrix@R=0.5cm@C=0.5cm
{
&1\ar[ld]_2&\\
2\ar[rr]_{2}&&3\ar[lu]_{2}\\
}
\xymatrix@R=0.5cm@C=0.6cm
{
&\\
\ar@{-}[r]^{\mu_2}&\\
}
\xymatrix@R=0.5cm@C=0.5cm
{
&1\ar[rd]^{2}&\\
2\ar[ru]^{2}&&3\ar[ll]^2\\
}
\xymatrix@R=0.5cm@C=0.6cm
{
&\\
\ar@{-}[r]^{\mu_3}&\\
}
\xymatrix@R=0.5cm@C=0.5cm
{
&1\ar[ld]_2&\\
2\ar[rr]_{2}&&3\ar[lu]_{2}\\
}
\\
&\hspace{3.3cm}
\xymatrix@R=0.5cm@C=0.8cm
{
&\ar@{-}[d]_{\mu_1}& \\
&&\\
}
\\
&\hspace{3.3cm}
\xymatrix@R=0.5cm@C=0.5cm
{
&1\ar[ld]_2&\\
2\ar[rr]_{2}&&3.\ar[lu]_{2}\\
}
\end{align*}
}
It follows that $|H_1(t_0)|=3$ and we may take 
\[
    H_1(t_0)=\{g_i^{\bf{id},-}~|~i=1,2,3\}.
\]
    It is clear that $(g_i^{\bf{id},-})^2=\mathbf{id}$ for $i=1,2,3$. Now we turn to determine $G_0(t_0)$:
\begin{itemize}
	\item If $d_1=d_2\neq d_3$, then $G_0(t_0)\cong \mathbb{Z}_2$ and $\psi_{(12)}^-$ is a generator of $G_0(t_0)$.
	\item If $d_1=d_3\neq d_2$, then $G_0(t_0)\cong \mathbb{Z}_2$ and $\psi_{(13)}^-$ is a generator of $G_0(t_0)$.
	\item If $d_2=d_3\neq d_1$, then $G_0(t_0)\cong \mathbb{Z}_2$ and $\psi_{(23)}^-$ is a generator of $G_0(t_0)$.
	
	\item If $d_1=d_2=d_3$, then $G_0(t_0)\cong S_3$
\end{itemize}

For each $\sigma\in S_3$, let $\operatorname{sgn}(\sigma)$ be the sign of $\sigma$. 
It is straightforward to check that 
\[
\psi_\sigma^{\operatorname{sgn}(\sigma)}\circ g_i^{\bf{id},-}=g_{\sigma(i)}^{\bf{id},-} \circ \psi_\sigma^{\operatorname{sgn}(\sigma)}
\]
for any $\sigma\in S_3$ and $1\leq i\leq 3$. We claim that 
\[
\operatorname{Aut}(\mathcal{A})\cong \begin{cases}
    N\rtimes_{\tau} \mathbb{Z}_2& \text{exactly two of $d_1,d_2,d_3$ are equal;}\\
    N\rtimes_{\varphi}S_3 & d_1=d_2=d_3;\\
    N&\text{else}.
\end{cases}
\]
It suffices to show that for any finite reduced sequence $i_1,\dots, i_n$ of vertices of $\Gamma_0$ and $\sigma\in S_3$,  we have $g^{\bf{id},-}_{i_1}\circ\cdots \circ g^{\bf{id},-}_{i_n}=\psi_\sigma^{\operatorname{sgn}(\sigma)}$ if and only if  $n=0$ and $\sigma=\mathbf{id}$. 

Denote by $t_k=\mu_{i_k\dots i_1}(t_0)\in \mathbb{T}_3$ for $1\leq k\leq n$.
 We claim that the following inequalities hold:
\begin{equation}\label{ineq:d-vector}
	\begin{cases}
		{\bf d}_{i;t_{k}}\ge  {\bf d}_{i;t_{k-1}} & i\in \{1,2,3\},2\le k\le n, \\ 
		{\bf d}_{i_k;t_{k}}\ge \mathbf{0} & 1\le k\le n, \\
		b_{i_k i''_k} {\bf d}_{i_k;t_{k}}\ge b_{i'_ki''_k}{\bf d}_{i'_k;t_{k}} 
		& 1\le k\le n,\{i_k,i'_k,i''_k\}=\{1,2,3\}.   
	\end{cases}
\end{equation}
Since $\mathbf{d}_{i,t_0}=-e_i$ for $i=1,2,3$, where $e_1,e_2,e_3$ is the standard basis of $\mathbb{Z}^3$. It follows that (\ref{ineq:d-vector}) holds for $k=1$ by definition. Suppose that $(\ref{ineq:d-vector})$ holds for some $n>k\geq 1$. Let us prove it for $k+1$. Without loss of generality, we may assume that $i_k=1$
and $i_{k+1}=2$. Thus, by induction hypothesis, we have
\begin{align}
	b_{12}{\bf d}_{1;t_{k}} \ge b_{32} {\bf d}_{3;t_{k}},  \\
	b_{13}{\bf d}_{1;t_{k}} \ge b_{23} {\bf d}_{2;t_{k}}. 
\end{align}

It follows that 
\begin{align*}
	{\bf d}_{2;t_{k+1}} &=-{\bf d}_{2;t_{k}}+b_{12} {\bf d}_{1;t_{k}} \ge {\bf d}_{2;t_{k}} & \text{by $(5.3)\times b_{32}$,} \\
	b_{23} {\bf d}_{2;t_{k+1}} &=-b_{23}{\bf d}_{2;t_{k}}+2b_{13} {\bf d}_{1;t_{k}} \ge  b_{13} {\bf d}_{1;t_{k}}= b_{13}{\bf d}_{1;t_{k+1}} &  \\
	b_{21} {\bf d}_{2;t_{k+1}} &=-b_{21}{\bf d}_{2;t_{k}}+4 {\bf d}_{1;t_{k}} \ge b_{31}  {\bf d}_{3;t_{k}}=b_{31}  {\bf d}_{3;t_{k+1}}  
	&\text{by $ (5.2)\times b_{21}+(5.3)\times b_{31}$.}
\end{align*}
This finishes the proof of $(\ref{ineq:d-vector})$.
By Proposition \ref{prop:expression-mutation-composition}, 
\[f:=g^{\bf{id},-}_{i_1}\circ\cdots \circ g^{\bf{id},-}_{i_n}=(\Sigma_{t_0},\Sigma_{t_n},\mathbf{id},(-)^{n}).\] 
It follows that $f$ maps at least one of initial cluster variables to a non initial cluster variable. On the other hand, $\psi_\sigma^{\operatorname{sgn}(\sigma)}$ maps initial cluster variables to initial cluster variables. Therefore $f=\psi_\sigma^{\operatorname{sgn}(\sigma)}$ if and only if $n=0$ and $\sigma=\mathbf{id}$.

\subsubsection{ $2=a<b=c$.}
By mutating $\Gamma_0$ in three directions, we obtain
{\tiny \begin{align*}
&
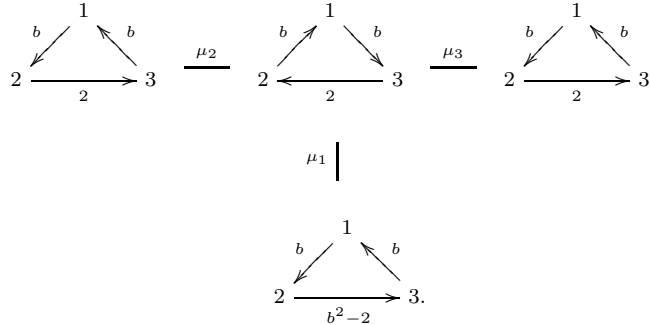
\xymatrix@R=0.5cm@C=0.5cm
{
&1\ar[ld]_b&\\
2\ar[rr]_{2}&&3\ar[lu]_{b}\\
}
\xymatrix@R=0.5cm@C=0.6cm
{
&\\
\ar@{-}[r]^{\mu_2}&\\
}
\xymatrix@R=0.5cm@C=0.5cm
{
&1\ar[rd]^{b}&\\
2\ar[ru]^{b}&&3\ar[ll]^2\\
}
\xymatrix@R=0.5cm@C=0.6cm
{
&\\
\ar@{-}[r]^{\mu_3}&\\
}
\xymatrix@R=0.5cm@C=0.5cm
{
&1\ar[ld]_b&\\
2\ar[rr]_{2}&&3\ar[lu]_{b}\\
}
\\
&\hspace{3.3cm}
\xymatrix@R=0.5cm@C=0.8cm
{
&\ar@{-}[d]_{\mu_1}& \\
&&\\
}
\\
&\hspace{3.5cm}
\xymatrix@R=0.5cm@C=0.5cm
{
&1\ar[ld]_b&\\
2\ar[rr]_{b^2-2}&&3.\ar[lu]_{b}\\
}
\end{align*}
}
Note that $\Gamma_{\mu_1(t_0)}$ has property $(M)$ at vertex $1$. It follows that for any $t\in \mathbb{T}_3$ such that the edge $\xymatrix{\mu_1(t_0)\ar@{-}[r]^k &\mu_k\mu_1(t_0)}$ of the path $p(\mu_1(t_0),t)$ is labeled by $k\neq 1$, $s(\Gamma_t)>s(\Gamma_{\mu_1(t_0)})>s(\Gamma_0)$ by Lemma \ref{lem:weight-increasing}. As a consequence, $|H_1(t_0)|=2$, and we may take
\[
    H_1(t_0)=\{g_2^{\bf{id},-},g_3^{\bf{id},-}\}.
\]
It is clear that $(g_2^{\bf{id},-})^2=(g_3^{\bf{id},-})^2=\mathbf{id}$.

On the other hand, $G_0(t_0)\neq \{\mathbf{id}\}$ if and only if $d_2=d_3$. In this case, we have $G_0(t_0)=\{\mathbf{id}, \psi_{(2,3)}^-\}$ and $(\psi_{(2,3)}^-)^2=\mathbf{id}$. Furthermore, $\psi_{(2,3)}^-\circ g_2^{\bf{id},-}=g_3^{\bf{id},-} \circ \psi_{(2,3)}^-$. Similar to Section \ref{sss:a=b=c=2}, one can show that $\operatorname{Aut}(\mathcal{A})\cong D_\infty $

\section{Proof of Theorem \ref{thm:acyclic-case}: Mutation-acyclic cases}\label{s:mutation-acyclic}
In this section, we prove Theorem \ref{thm:acyclic-case}. Let $\Gamma_0$ be the diagram in Theorem \ref{thm:classification-rank-3-matrix}(2).
Denote by $B_{t_0}=\begin{pmatrix}
	0&-b_{12}&b_{13}\\
	b_{21}& 0 & b_{23}\\
	-b_{31} & -b_{32} & 0
\end{pmatrix}$ the exchange matrix at vertex $t_0$, where $b_{ij}\in \mathbb{N}$. Let $\mathbf{x}_{t_0}=\mathbf{x}=(x_1,x_2,x_3)$. Recall that $b_{12}b_{21}=b^2$, $b_{23}b_{32}=a^2$, $b_{13}b_{31}=b^2$ and $0\leq a\leq b\leq c$.

\subsection{Preliminary results} In this subsection, we collect basic properties for certain cluster automorphisms that will be used frequently in the following sections.


\begin{lem}\label{lem:relation-rank-2-finite}
The following statements hold:
\begin{itemize}
    \item if $b=1$, then $\mu_{12121}(\Gamma(B))=\sigma(\Gamma(B))$, where $\sigma=(1 \,\, 2)$, and 
\[\mu_{12121}(x_1,x_2,x_3)=(x_2,x_1,x_3).\]
\item if $b=\sqrt 2$, then $(\mu_{21})^3(\Gamma(B))=\Gamma(B)$, and $(\mu_{21})^3(x_1,x_2,x_3)=(x_1,x_2,x_3)$.
\item if $b=\sqrt 3$, then $(\mu_{21})^4(\Gamma(B))=\Gamma(B)$, and $(\mu_{21})^4(x_1,x_2,x_3)=(x_1,x_2,x_3)$.
\end{itemize} 
\end{lem}
\begin{proof}
One can deduce  the statements from the periodicity of cluster algebras (with coefficients) of type $A_2,B_2$ and $G_2$, or alternatively, by direct calculation. For example,  to show $\mu_{12121}(x_1,x_2,x_3)=(x_2,x_1,x_3)$, one may use Lemma \ref{lem:d-vector-initial} and \ref{lem:mutation-d-vector} to compute the denominator vectors.

\end{proof}

Consider the mutation $\mu_{312}(\Gamma_0)$ of $\Gamma_0$, we obtain that 
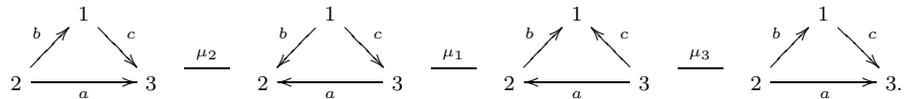
\[
{\tiny \xymatrix@R=0.5cm@C=0.5cm
{
&1\ar[rd]^c&\\
2\ar[ru]^b\ar[rr]_a&&3\\
}
\xymatrix@R=0.5cm@C=0.6cm
{
&\\
\ar@{-}[r]^{\mu_2}&\\
}
\xymatrix@R=0.5cm@C=0.5cm
{
&1\ar[rd]^c\ar[ld]_b&\\
2&&3\ar[ll]^a\\
}
\xymatrix@R=0.5cm@C=0.6cm
{
&\\
\ar@{-}[r]^{\mu_1}&\\
}
\xymatrix@R=0.5cm@C=0.5cm
{
&1&\\
2\ar[ru]^b&&3\ar[lu]_c\ar[ll]^a\\
}
\xymatrix@R=0.5cm@C=0.6cm
{
&\\
\ar@{-}[r]^{\mu_3}&\\
}
\xymatrix@R=0.5cm@C=0.5cm
{
&1\ar[rd]^c&\\
2\ar[ru]^b\ar[rr]_a&&3.\\
}}
\]
Namely, $\Gamma_0=\mu_{312}(\Gamma_0)$.
It follows that $g_{312}^{\mathbf{id},+}$ belongs to $\operatorname{Aut}(\mathcal{A})$. The following result can be proved using representation theory of valued quivers. Here we provide an elementary proof using denominator vectors.
\begin{lem}\label{lem:order-of-tau}
 We have
\[
|g_{312}^{\mathbf{id},+}|=\begin{cases}
6,& a=0,b=c=1;\\
4,&a=0,b=1,c=\sqrt 2;\\
\infty,&\text{else}.
\end{cases}
\]
Here $|g_{312}^{\mathbf{id},+}|$ is the order of $g_{312}^{\mathbf{id},+}$ in $\operatorname{Aut}(\mathcal{A})$.
\end{lem}
\begin{proof}
For any positive integer $n$, denote by $t_n:=(\mu_{312})^n(t_0)$.

We first assume that $a\geq 1$, hence $b_{ij}\geq 1$ for all $i,j$. We claim that
\begin{eqnarray}\label{ineq:d-vector-acyclic-1}
{\bf d}_{3;t_n}\ge {\bf d}_{1;t_n}\ge{\bf d}_{2;t_n} \ge\mathbf{0} \hspace{1cm} n\in \mathbb{N}_+.
\end{eqnarray}
It is evident that the claim is true for $n=1$ by direct calculation. We assume that the statement holds for $n=k$, and prove it for $n=k+1$. By Lemma \ref{lem:mutation-d-vector}, we have
\begin{align}
{\bf d}_{2;t_{k+1}}&=-{\bf d}_{2;t_k}+(b_{12}{\bf d}_{1;t_k}+b_{32}{\bf d}_{3;t_k})\ge {\bf d}_{3;t_k},\label{ineq:d-vector-acyc-2} \\
{\bf d}_{1;t_{k+1}}&=-{\bf d}_{1;t_k}+(b_{21}{\bf d}_{2;t_{k+1}}+b_{31}{\bf d}_{3;t_k})\ge {\bf d}_{2;t_{k+1}},\label{ineq:d-vector-acyc-3} \\
{\bf d}_{3;t_{k+1}}&=-{\bf d}_{3;t_k}+(b_{23}{\bf d}_{2;t_{k+1}}+b_{13}{\bf d}_{1;t_{k+1}})\ge {\bf d}_{1;t_{k+1}}, \label{ineq:d-vector-acyc-4}
\end{align}
which completes the proof of (\ref{ineq:d-vector-acyclic-1}). Consequently, $(g_{312}^{\mathbf{id},+})^n\neq \mathbf{id}$ for any $n\in \mathbb{N}_+$.

Now assume that $a=0$ and $m:=b^2+c^2\geq 4$. We claim that
\begin{equation}\label{ineq:d-vector-acyc-5}
  \begin{cases}
  {\bf d}_{1;t_{n+1}}\ge {\bf d}_{1;t_n} \ge\mathbf{0} & n\in \mathbb{N}_+;\\
  b_{12}{\bf d}_{1;t_{n}}\ge {\bf d}_{2;t_n} \ge\mathbf{0} &  n\in \mathbb{N}_+; \\
  b_{13}{\bf d}_{1;t_{n+1}}\ge {\bf d}_{3;t_n} \ge\mathbf{0} &  n\in \mathbb{N}_+.
\end{cases}  
\end{equation}

A direct calculation shows that (\ref{ineq:d-vector-acyc-5}) holds for $n=1$.
 Assume the statement is true for $n=k$. According to (\ref{ineq:d-vector-acyc-2})-(\ref{ineq:d-vector-acyc-4}), we obtain
\begin{align*}
{\bf d}_{1;t_{k+1}}&=(b^2-1){\bf d}_{1;t_k}-b_{21}{\bf d}_{2;t_{k}}+b_{31}{\bf d}_{3;t_k},  \\
{\bf d}_{2;t_{k+1}}&=b_{12}{\bf d}_{1;t_k} -{\bf d}_{2;t_k}, \\
{\bf d}_{3;t_{k+1}}&=b_{13}(b^2-1){\bf d}_{1;t_{k}}-b_{13}b_{21}{\bf d}_{2;t_{k}}+(c^2-1){\bf d}_{3;t_k}.
\end{align*}
By induction hypothesis, we further have
\begin{align*}
{\bf d}_{1;t_{k+1}}-{\bf d}_{1;t_{k}}&=(b^2-2){\bf d}_{1;t_k}-b_{21}{\bf d}_{2;t_{k}}+b_{31}{\bf d}_{3;t_k}\ge \mathbf{0}, \\
{\bf d}_{2;t_{k+1}}&=b_{12}{\bf d}_{1;t_{k}}-{\bf d}_{2;t_{k}} \ge \mathbf{0}, \\
{\bf d}_{3;t_{k+1}} &= b_{13}{\bf d}_{1;t_{k+1}}-{\bf d}_{3;t_{k}} \ge \mathbf{0}.
\end{align*}
Putting all of these together, we compute
\begin{align*}
{\bf d}_{1;t_{k+2}}-{\bf d}_{1;t_{k+1}} &=(b^2-2)[(b^2-1){\bf d}_{1;t_k}-b_{21}{\bf d}_{2;t_{k}}+b_{31}{\bf d}_{3;t_k}]-b_{21}(b_{12}{\bf d}_{1;t_k} -{\bf d}_{2;t_k}) \\
&+b_{31}[b_{13}(b^2-1){\bf d}_{1;t_{k}}-b_{13}b_{21}{\bf d}_{2;t_{k}}+(c^2-1){\bf d}_{3;t_k}]  \\
&=(mb^2-3b^2-m+2) {\bf d}_{1;t_{k}}-b_{21}(m-3){\bf d}_{2;t_{k}}+b_{31}(m-3){\bf d}_{3;t_{k}} \\
&\ge (mb^2-3b^2-2m+6){\bf d}_{1;t_{k}}-b_{21}(m-3){\bf d}_{2;t_{k}}+b_{31}(m-3){\bf d}_{3;t_{k}} \\
&=(m-3)[(b^2-2){\bf d}_{1;t_k}-b_{21}{\bf d}_{2;t_{k}}+b_{31}{\bf d}_{3;t_k}] \ge \mathbf{0},\\
b_{12}{\bf d}_{1;t_{k+1}}-{\bf d}_{2;t_{k+1}}&=b_{12}{\bf d}_{1;t_{k+1}}-b_{12}{\bf d}_{1;t_{k}}+{\bf d}_{2;t_{k}}=b_{12}({\bf d}_{1;t_{k+1}}-{\bf d}_{1;t_{k}})+{\bf d}_{2;t_{k}} \ge \mathbf{0},\\
b_{13}{\bf d}_{1;t_{k+2}}-{\bf d}_{3;t_{k+1}}&=b_{13}(b^2-1)[(b^2-1){\bf d}_{1;t_k}-b_{21}{\bf d}_{2;t_{k}}+b_{31}{\bf d}_{3;t_k}]-b_{13}b_{21}(b_{12}{\bf d}_{1;t_k} -{\bf d}_{2;t_k}) \\
&+(c^2-1)[b_{13}(b^2-1){\bf d}_{1;t_{k}}-b_{13}b_{21}{\bf d}_{2;t_{k}}+(c^2-1){\bf d}_{3;t_k}] \\
&=b_{13}(mb^2-3b^2-m+2){\bf d}_{1;t_{k}}-b_{13}b_{21}(m-3){\bf d}_{2;t_{k}}+(mc^2-3c^2+1){\bf d}_{3;t_{k}} \\
&\ge b_{13}(mb^2-3b^2-2m+6){\bf d}_{1;t_{k}}-b_{13}b_{21}(m-3){\bf d}_{2;t_{k}}+(mc^2-3c^2){\bf d}_{3;t_{k}} \\
&=b_{13}(m-3)[(b^2-2){\bf d}_{1;t_k}-b_{21}{\bf d}_{2;t_{k}}+b_{31}{\bf d}_{3;t_k}] \ge \mathbf{0}.
\end{align*}
This completes the proof of (\ref{ineq:d-vector-acyc-5}). As a consequence, $(g_{312}^{\mathbf{id},+})^n\neq \mathbf{id}$ for any $n\in \mathbb{N}_+$.

It remains to consider the cases $a=0,b=1,c=1$ and $a=0,b=1,c=\sqrt{2}$, which can be calculated directly via denominator vector using Lemma \ref{lem:d-vector-initial} and \ref{lem:mutation-d-vector}.

\end{proof}

\subsection{Case 1: $a=0$.}\label{ss:a=0}
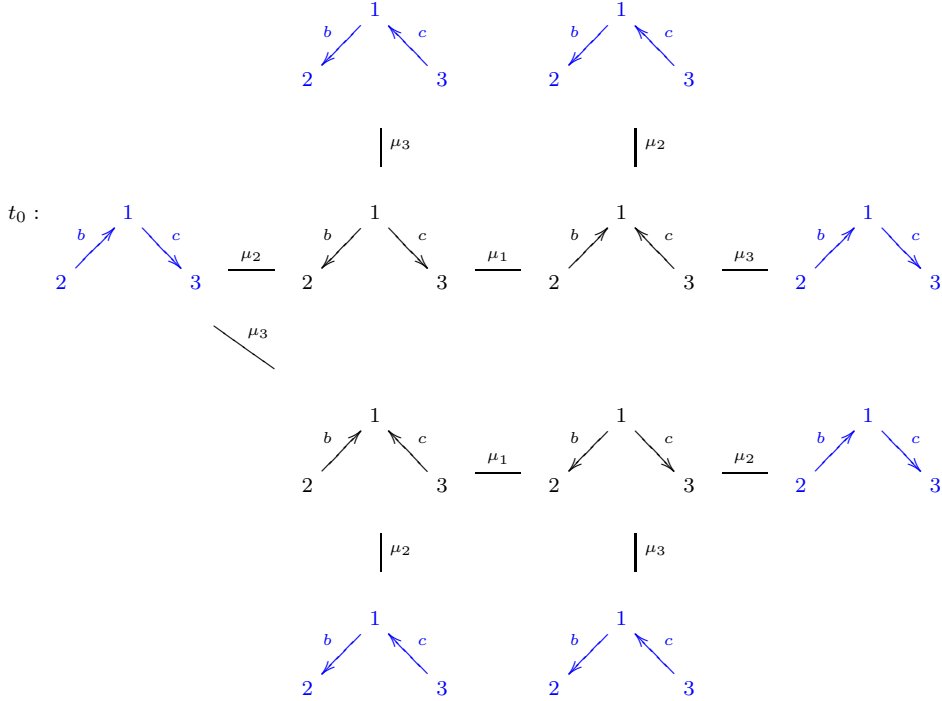
\begin{figure}
\centering
{\tiny
\begin{align*}
&
{\color{blue}\xymatrix@R=0.5cm@C=0.5cm
{
&1\ar[ld]_b&\\
2&&3\ar[lu]_c\\
}}
\xymatrix@R=0.5cm@C=0.6cm
{
&\\
&\\
}
{\color{blue}\xymatrix@R=0.5cm@C=0.5cm
{
&1\ar[ld]_b&\\
2&&3\ar[lu]_c\\
}}
\\
&
\xymatrix@R=0.5cm@C=0.8cm
{
&\ar@{-}[d]^{\mu_3}&\\
&&\\
}
\hspace{1.05cm}
\xymatrix@R=0.5cm@C=0.8cm
{
&\ar@{-}[d]^{\mu_2}&\\
&&\\
}
\\
t_0:{\color{blue}\xymatrix@R=0.5cm@C=0.5cm
{
&1\ar[rd]^c&\\
2\ar[ru]^b&&3\\
}}
\xymatrix@R=0.5cm@C=0.6cm
{
&\\
\ar@{-}[r]^{\mu_2}&\\
}&
\xymatrix@R=0.5cm@C=0.5cm
{
&1\ar[rd]^c\ar[ld]_b&\\
2&&3\\
}
\xymatrix@R=0.5cm@C=0.6cm
{
&\\
\ar@{-}[r]^{\mu_1}&\\
}
\xymatrix@R=0.5cm@C=0.5cm
{
&1&\\
2\ar[ru]^b&&3\ar[lu]_c\\
}
\xymatrix@R=0.5cm@C=0.6cm
{
&\\
\ar@{-}[r]^{\mu_3}&\\
}
{\color{blue}\xymatrix@R=0.5cm@C=0.5cm
{
&1\ar[rd]^c&\\
2\ar[ru]^b&&3\\
}}
\\
\xymatrix@R=0.5cm@C=0.8cm{
&\ar@{-}[rd]^{\mu_3}&\\
&&\\
}
\\
\xymatrix@R=0.5cm@C=0.5cm
{
&&\\
&&\\
}
\xymatrix@R=0.5cm@C=0.6cm
{
&\\
&\\
}&
\xymatrix@R=0.5cm@C=0.5cm
{
&1&\\
2\ar[ru]^b&&3\ar[lu]_c\\
}
\xymatrix@R=0.5cm@C=0.6cm
{
&\\
\ar@{-}[r]^{\mu_1}&\\
}
\xymatrix@R=0.5cm@C=0.5cm
{
&1\ar[ld]_b\ar[rd]^c&\\
2&&3\\
}
\xymatrix@R=0.5cm@C=0.6cm
{
&\\
\ar@{-}[r]^{\mu_2}&\\
}
{\color{blue}\xymatrix@R=0.5cm@C=0.5cm
{
&1\ar[rd]^c&\\
2\ar[ru]^b&&3\\
}}
\\
&
\xymatrix@R=0.5cm@C=0.8cm
{
&\ar@{-}[d]^{\mu_2}&\\
&&\\
}
\hspace{1.05cm}
\xymatrix@R=0.5cm@C=0.8cm
{
&\ar@{-}[d]^{\mu_3}&\\
&&\\
}
\\
&
{\color{blue}\xymatrix@R=0.5cm@C=0.5cm
{
&1\ar[ld]_b&\\
2&&3\ar[lu]_c\\
}}
\xymatrix@R=0.5cm@C=0.6cm
{
&\\
&\\
}
{\color{blue}\xymatrix@R=0.5cm@C=0.5cm
{
&1\ar[ld]_b&\\
2&&3\ar[lu]_c\\
}}
\end{align*}
}
 \caption{Mutation subgraph of $\Gamma_0$  of $(0,b,c)$ in directions $2$ and $3$. All the diagrams at vertex $t$ such that $B_t\sim \pm B_{t_0}$ are colored blue.}
    \label{fig:0bc-23}
\end{figure}

Since we have assumed that the exchange matrix $B_{t_0}$ is indecomposable, we have $1\leq b$.
By mutating $\Gamma_0$ in directions $2$ and $3$, we obtain the mutation subgraph of $\Gamma_0$ as shown in Figure \ref{fig:0bc-23}.


Recall that each path in $\mathcal{P}_{1}(t_0)$ induces an element in $H_1(t_0)$.
Hence, we obtain the following $6$ elements in $H_1(t_0)$:
\[
    g_{312}^{\mathbf{id},+}, g_{213}^{\mathbf{id},+}, g_{32}^{\mathbf{id},-}, g_{23}^{\mathbf{id},-}, g_{212}^{\mathbf{id},-}, g_{313}^{\mathbf{id},-}.
\]
By Lemma \ref{prop:expression-mutation-composition} and \ref{lem:mutation-d-vector}, a direct computation shows that
\begin{align*}
&g_{312}^{\mathbf{id},+}\circ g_{213}^{\mathbf{id},+}=g_{32}^{\mathbf{id},-}\circ g_{23}^{\mathbf{id},-}=(g_{32}^{\mathbf{id},-})^2=\mathbf{id}, \\
&g_{212}^{\mathbf{id},-}=g_{312}^{\mathbf{id},+}\circ g_{23}^{\mathbf{id},-},\quad 
g_{313}^{\mathbf{id},-}=g_{213}^{\mathbf{id},+}\circ g_{32}^{\mathbf{id},-},\quad  g_{312}^{\mathbf{id},+}\circ g_{32}^{\mathbf{id},-}=g_{32}^{\mathbf{id},-}\circ (g_{312}^{\mathbf{id},+})^{-1}.
\end{align*}
In particular, these $6$ elements are generated by $g_{312}^{\mathbf{id},+}$ and $g_{32}^{\mathbf{id},-}$. More elements in $H_1(t_0)$ and $G_0(t_0)$ depend on the values of $b$ and $c$. Hence, we split the discussion by considering the different choices of $b$ and $c$.

\subsubsection{$a=0,b=1$.}
 Applying $\mu_2\mu_1$ to $\Gamma_0$, we obtain
{\tiny
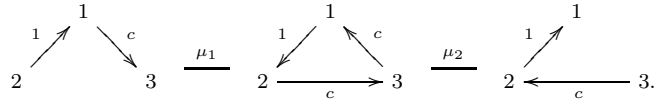
\begin{align*}
\xymatrix@R=0.5cm@C=0.5cm
{
&1\ar[rd]^{c}&\\
2\ar[ru]^1&&3\\
}
\xymatrix@R=0.5cm@C=0.6cm
{
&\\
\ar@{-}[r]^{\mu_1}&\\
}
\xymatrix@R=0.5cm@C=0.5cm
{
&1\ar[dl]_1&\\
2\ar[rr]_{c}&&3\ar[lu]_{c}\\
}
\xymatrix@R=0.5cm@C=0.6cm
{
&\\
\ar@{-}[r]^{\mu_2}&\\
}
\xymatrix@R=0.5cm@C=0.5cm
{
&1&\\
2\ar[ru]^1&&3.\ar[ll]^{c}\\
}
\end{align*}}
It follows that $g_{21}^{(12),-}\in G_1(t_0)$. Moreover, \[g_{21}^{(12),-}\circ g_{212}^{\mathbf{id},-}=(\Sigma_{t_0},\Sigma_{\mu_1\mu_2\mu_1\mu_2\mu_1(t_0)},(1 2),+)=\mathbf{id}\] by Proposition \ref{prop:expression-mutation-composition} and Lemma \ref{lem:relation-rank-2-finite}. As a consequence, $g_{21}^{(12),-}=g_{32}^{\mathbf{id},-}\circ (g_{312}^{\mathbf{id},+})^{-1}$.

\subsubsection{$a=0,b=c=1$.} It turns out that  $\mathcal{A}$ is of type $A_3$. Its cluster automorphism group $\operatorname{Aut}(\mathcal{A})$ is isomorphic to $D_6$, cf. \cite{ASS2012}.

\subsubsection{$a=0,b=1,c=\sqrt{2}$.}\label{ss:a=0-b=1-c=sqrt2}
It turns out that $\mathcal{A}$ is of type $B_3$ or $C_3$, and its cluster automorphism group $\operatorname{Aut}(\mathcal{A})$ is isomorphic to $D_4$, cf. \cite{CZ16a}.

\subsubsection{$a=0,b=1,c=\sqrt{3}$.}\label{ss:a=0-b=1-c=sqrt3}
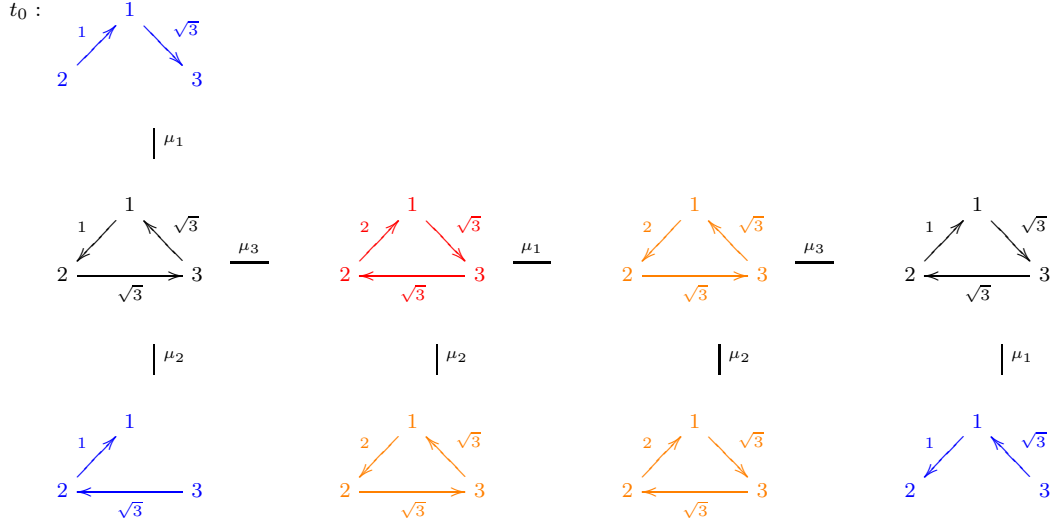
\begin{figure}
\centering
{\tiny
\begin{align*}
t_0:
{\color{blue} \xymatrix@R=0.5cm@C=0.5cm
{
&1\ar[rd]^{\sqrt 3}&\\
2\ar[ru]^1&&3\\
}}&&&&&&\\
\xymatrix@R=0.4cm@C=0.4cm
{
&\ar@{-}[d]^{\mu_1}&\\
&&\\
}
&&&&&&\\
\xymatrix@R=0.5cm@C=0.5cm
{
&1\ar[ld]_1&\\
2\ar[rr]_{\sqrt 3}&&3\ar[lu]_{\sqrt 3}\\
}
&
\xymatrix@R=0.5cm@C=0.5cm
{
&\\
\ar@{-}[r]^{\mu_3}&\\
}
&
{\color{red}\xymatrix@R=0.5cm@C=0.5cm
{
&1\ar[rd]^{\sqrt3}&\\
2\ar[ru]^2 &&3\ar[ll]^{\sqrt3}\\
} }
&
\xymatrix@R=0.5cm@C=0.5cm
{
&\\
\ar@{-}[r]^{\mu_1}&\\
}
&
{\color{orange}\xymatrix@R=0.5cm@C=0.5cm
{
&1\ar[ld]_{2}&\\
2\ar[rr]_{\sqrt 3}&&3\ar[lu]_{\sqrt3}\\
}}
&
\xymatrix@R=0.5cm@C=0.5cm
{
&\\
\ar@{-}[r]^{\mu_3}&\\
}
&
\xymatrix@R=0.5cm@C=0.5cm
{
&1\ar[rd]^{\sqrt 3}&\\
2\ar[ru]^1&&3\ar[ll]^{\sqrt3}
}\\
\xymatrix@R=0.4cm@C=0.4cm
{
&\ar@{-}[d]^{\mu_2}&\\
&&\\
}
&&
\xymatrix@R=0.4cm@C=0.4cm
{
&\ar@{-}[d]^{\mu_2}&\\
&&\\
}
&&
\xymatrix@R=0.4cm@C=0.4cm
{
&\ar@{-}[d]^{\mu_2}&\\
&&\\
}
&&
\xymatrix@R=0.4cm@C=0.4cm
{
&\ar@{-}[d]^{\mu_1}&\\
&&\\
} \\
{\color{blue}\xymatrix@R=0.5cm@C=0.5cm
{
&1&\\
2\ar[ru]^1&&3\ar[ll]^{\sqrt 3}\\
}}
&&
{\color{orange}\xymatrix@R=0.5cm@C=0.5cm
{
&1\ar[ld]_{2}&\\
2\ar[rr]_{\sqrt 3}&&3\ar[lu]_{\sqrt3}\\
}}
&&
{\color{orange}\xymatrix@R=0.5cm@C=0.5cm
{
&1\ar[rd]^{\sqrt3}&\\
2\ar[ru]^2 &&3\ar[ll]^{\sqrt3}\\
}}
&&
{\color{blue}\xymatrix@R=0.5cm@C=0.5cm
{
&1\ar[ld]_1&\\
2&&3\ar[lu]_{\sqrt3}
}}\\
\end{align*}
}
\caption{Mutation subgraph of $\Gamma_0$  of $(0,1,\sqrt{3})$ in direction $1$. The diagran of $t_\diamond$ is colored red and the diagrams at vertex $t$ such that $B_t\sim \pm B_{t_\diamond}$ are colored in orange.}
    \label{fig:01sqrt3-1}
\end{figure}

 It is clear that $d_1=d_2$ and $G_0(t_0)=\langle \mathbf{id}\rangle$. By mutating in direction $1$, we obtain the mutation subgraph of $\Gamma_0$, cf. Figure \ref{fig:01sqrt3-1}.

It is clear that $t_\diamond:=\mu_3\mu_1(t_0)$ satisfies the condition $(\diamond)$. Combining with Figure \ref{fig:0bc-23} in Section \ref{ss:a=0}, we conclude that $|\mathcal{P}_{1,0}^{t_\diamond}(t_0)|=7$, $|\mathcal{P}_{1,1}^{t_\diamond}(t_0)|=0$ and $|\mathcal{P}_{1,2}^{t_\diamond}(t_0)|=4$.
By Proposition \ref{prop:generator-set-assumption}, Remark \ref{rem:factorization-trivial-path} and Remark \ref{rem:sharpen-K-2}, $\operatorname{Aut}(\mathcal{A})$ is generated by the following elements:
\begin{align*}
    &g_{312}^{\mathbf{id},+}, g_{213}^{\mathbf{id},+}, g_{32}^{\mathbf{id},-}, g_{23}^{\mathbf{id},-}, g_{212}^{\mathbf{id},-}, g_{313}^{\mathbf{id},-}, g_{21}^{(12),-}&\in H_{1,0}^{t_\diamond}(t_0),\\
   & g_{13131}^{\mathbf{id},-}, g_{13231}^{\mathbf{id},-}&\in H_{1,2}^{t_\diamond}(t_0).
\end{align*}
By Proposition \ref{prop:expression-mutation-composition} and Lemma \ref{lem:relation-rank-2-finite}, we obtain
\[
    g_{21}^{(12),-}\circ g_{212}^{\mathbf{id},-}=\mathbf{id}=g_{13131}^{\mathbf{id},-}\circ g_{313}^{\mathbf{id},-}.
\]
On the other hand, \[g_{21}^{(12),-}\circ g_{13231}^{\mathbf{id},-}\circ g_{21}^{(12),-}=g_{13131}^{\mathbf{id},-}\] by Proposition \ref{prop:expression-mutation-composition}. As a consequence, $\operatorname{Aut}(\mathcal{A})$ is generated by $g_{312}^{\mathbf{id},+}$ and $g_{32}^{\mathbf{id},-}$.

We claim that for any $i\in \mathbb{Z}$ and $0\leq j\leq 1$ such that $(g_{312}^{\mathbf{id},+})^i\circ (g_{32}^{\mathbf{id},-})^j=\mathbf{id}$ if and only if $i=j=0$. Indeed, $(g_{312}^{\mathbf{id},+})^i\circ (g_{32}^{\mathbf{id},-})^j(x_1)=(g_{312}^{\mathbf{id},+})^i(x_1)=x_1$, which implies that $i=0$ by Lemma \ref{lem:order-of-tau} and its proof. Consequently, $j=0$. It follows that $\operatorname{Aut}(\mathcal{A})\cong D_\infty$.


\subsubsection{$a=0,b=1,c\geq 2$.}

By mutating in direction $1$, we obtain the following mutation subgraph of $\Gamma_0$:
{\tiny 
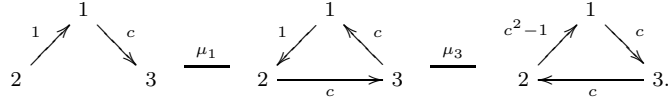
\begin{align*}
\xymatrix@R=0.5cm@C=0.5cm
{
&1\ar[rd]^c&\\
2\ar[ru]^1&&3\\
}
\xymatrix@R=0.5cm@C=0.6cm
{
&\\
\ar@{-}[r]^{\mu_1}&\\
}
\xymatrix@R=0.5cm@C=0.5cm
{
&1\ar[dl]_1&\\
2\ar[rr]_c&&3\ar[lu]_c\\
}
\xymatrix@R=0.5cm@C=0.6cm
{
&\\
\ar@{-}[r]^{\mu_3}&\\
}
\xymatrix@R=0.5cm@C=0.5cm
{
&1\ar[dr]^c&\\
2\ar[ur]^{c^2-1}&&3.\ar[ll]^c\\
}
\end{align*}}
It is clear that $\Gamma_{\mu_{31}(t_0)}$ has property $(M)$ at vertex $3$. It follows that  for any vertex $t\in \mathbb{T}_3$ such that the first edge $\xymatrix{\mu_{31}(t_0)\ar@{-}[r]^k&\mu_{k31}(t_0)}$ of the path $p(\mu_{31}(t_0),t)$ is labeled by $k\neq 3$, we have $s(\Gamma_t)>s(\Gamma_{t_0})$. Consequently, $\operatorname{Aut}(\mathcal{A})$ is generated by $g_{312}^{\mathbf{id},+}$ and $g_{32}^{\mathbf{id},-}$. The same proof in Section \ref{ss:a=0-b=1-c=sqrt3} yields that $\operatorname{Aut}(\mathcal{A})\cong D_\infty$.

\subsubsection{$a=0,b=\sqrt{2},c=\sqrt{2}$.} The mutation subgraph of $\Gamma_0$ in direction $1$ is given in Figure \ref{fig:0sqrt2sqrt2-1}.
\begin{figure}
\centering
{\tiny
\begin{align*}
t_0:{\color{blue}\xymatrix@R=0.5cm@C=0.5cm
{
&1\ar[rd]^{\sqrt 2}&\\
2\ar[ru]^{\sqrt 2}&&3\\
}}
\xymatrix@R=0.5cm@C=0.6cm
{
&\\
\ar@{-}[r]^{\mu_1}&\\
}&
{\color{red}\xymatrix@R=0.5cm@C=0.5cm
{
&1\ar[ld]_{\sqrt 2}&\\
2\ar[rr]_{2}&&3\ar[lu]_{\sqrt 2}\\
}}&
\xymatrix@R=0.5cm@C=0.6cm
{
&\\
\ar@{-}[r]^{\mu_2}&\\
}&
{\color{orange}\xymatrix@R=0.5cm@C=0.5cm
{
&1\ar[rd]^{\sqrt 2}&\\
2\ar[ru]^{\sqrt 2}&&3\ar[ll]^{2} \\
}}
&
\xymatrix@R=0.5cm@C=0.6cm
{
&\\
\ar@{-}[r]^{\mu_1}&\\
}
{\color{blue}\xymatrix@R=0.5cm@C=0.5cm
{
&1\ar[ld]_{\sqrt 2}&\\
2&&3\ar[lu]_{\sqrt2}\\
}}
\hspace{1cm}
\xymatrix@R=0.5cm@C=0.5cm
{
&&\\
&&\\
}\\
&
\xymatrix@R=0.5cm@C=0.8cm
{
&\ar@{-}[d]^{\mu_3}&\\
&&\\
}
&&
\xymatrix@R=0.5cm@C=0.8cm
{
&\ar@{-}[d]^{\mu_3}&\\
&&\\
}\\
&
{\color{orange}\xymatrix@R=0.5cm@C=0.5cm
{
&1\ar[rd]^{\sqrt 2}&\\
2\ar[ru]^{\sqrt 2}&&3\ar[ll]^{2} \\
}}
&&
{\color{orange}\xymatrix@R=0.5cm@C=0.5cm
{
&1\ar[ld]_{\sqrt 2}&\\
2\ar[rr]_{ 2 }&&3.\ar[lu]_{\sqrt 2}\\
}}
&
\end{align*}
}
\caption{Mutation graph of $\Gamma_0$ of $(0,\sqrt{2},\sqrt{2})$ in direction $1$.}
    \label{fig:0sqrt2sqrt2-1}
\end{figure}
Clearly, $t_\diamond=\mu_1(t_0)$ satisfies the condition $(\diamond)$. Similar to Section \ref{ss:a=0-b=1-c=sqrt3}, we obtain the following elements in $H_1(t_0)$:
\begin{align*}
&g_{312}^{\mathbf{id},+}, g_{213}^{\mathbf{id},+}, g_{32}^{\mathbf{id},-}, g_{23}^{\mathbf{id},-}, g_{212}^{\mathbf{id},-}, g_{313}^{\mathbf{id},-}&\in H_{1,0}^{t_\diamond}(t_0),\\
&g_{121}^{\mathbf{id},-}, g_{131}^{\mathbf{id},-}&\in H_{1,2}^{t_\diamond}(t_0).
\end{align*}
By Lemma \ref{lem:relation-rank-2-finite}, we have 
\[
    g_{121}^{\mathbf{id},-}\circ g_{212}^{\mathbf{id},-}=\mathbf{id}=g_{131}^{\mathbf{id},-}\circ g_{313}^{\mathbf{id},-}.
\]
Consequently, $\operatorname{Aut}(\mathcal{A})$ is generated by $g_{312}^{\mathbf{id},+}$, $g_{32}^{\mathbf{id},-}$ and $G_0(t_0)$.

If $d_2\neq d_3$, then $G_0(t_0)=\langle\mathbf{id}\rangle$. It follows that $\operatorname{Aut}(\mathcal{A})\cong D_\infty$.

If $d_2=d_3$, then $G_0(t_0)\cong \mathbb{Z}_2$ and $\psi_{(23)}^-$ is a generator of $G_0(t_0)$. A direct computation shows that 
\[g_{312}^{\mathbf{id},+}\circ \psi_{(23)}^-=\psi_{(23)}^-\circ (g_{312}^{\mathbf{id},+})^{-1}, \psi_{(23)}^-\circ g_{32}^{\mathbf{id},-}=g_{32}^{\mathbf{id},-}\circ \psi_{(23)}^-.
\]
Similar to Section \ref{ss:a=0-b=1-c=sqrt3}, one can show that $\operatorname{Aut}(\mathcal{A})\cong D_\infty\rtimes_\rho \mathbb{Z}_2$.


\subsubsection{$a=0,b=\sqrt{2},c=\sqrt{3}$.}

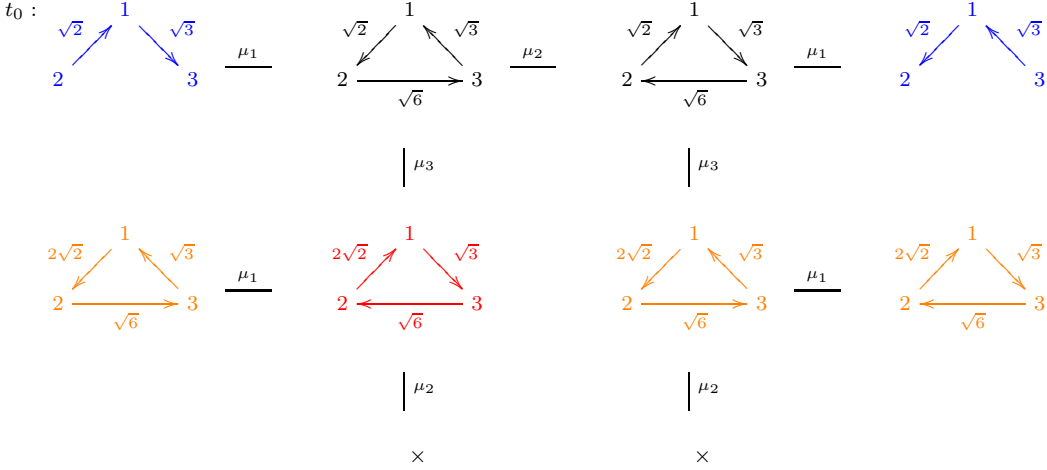
\begin{figure}
\centering
{\tiny
\begin{align*}
t_0:{\color{blue}\xymatrix@R=0.5cm@C=0.5cm
{
&1\ar[rd]^{\sqrt 3}&\\
2\ar[ru]^{\sqrt 2}&&3\\
}}
&
\xymatrix@R=0.5cm@C=0.6cm
{
&\\
\ar@{-}[r]^{\mu_1}&\\
}&
\xymatrix@R=0.5cm@C=0.5cm
{
&1\ar[ld]_{\sqrt 2}&\\
2\ar[rr]_{\sqrt 6}&&3\ar[lu]_{\sqrt 3}\\
}&
\xymatrix@R=0.5cm@C=0.6cm
{
&\\
\ar@{-}[r]^{\mu_2}&\\
}&
\xymatrix@R=0.5cm@C=0.5cm
{
&1\ar[rd]^{\sqrt 3}&\\
2\ar[ru]^{\sqrt 2}&&3\ar[ll]^{\sqrt 6} \\
}&
\xymatrix@R=0.5cm@C=0.6cm
{
&\\
\ar@{-}[r]^{\mu_1}&\\
}&
{\color{blue}\xymatrix@R=0.5cm@C=0.5cm
{
&1\ar[ld]_{\sqrt 2}&\\
2&&3\ar[lu]_{\sqrt 3}\\
}}
&
\xymatrix@R=0.5cm@C=0.5cm
{
&&\\
&&\\
}\\
&&
\xymatrix@R=0.5cm@C=0.8cm
{
&\ar@{-}[d]^{\mu_3}&\\
&&\\
}
&
&
\xymatrix@R=0.5cm@C=0.8cm
{
&\ar@{-}[d]^{\mu_3}&\\
&&\\
}
&
&&&
\\
{\color{orange}\xymatrix@R=0.5cm@C=0.5cm
{
&1\ar[ld]_{2\sqrt 2}&\\
2\ar[rr]_{\sqrt 6}&&3\ar[lu]_{\sqrt 3}\\
}}
&
\xymatrix@R=0.5cm@C=0.6cm
{
&\\
\ar@{-}[r]^{\mu_1}&\\
}
&
{\color{red}\xymatrix@R=0.5cm@C=0.5cm
{
&1\ar[rd]^{\sqrt 3}&\\
2\ar[ru]^{2\sqrt 2}&&3\ar[ll]^{\sqrt 6} \\
}}
&
&
{\color{orange}\xymatrix@R=0.5cm@C=0.5cm
{
&1\ar[ld]_{2\sqrt 2}&\\
2\ar[rr]_{\sqrt 6}&&3\ar[lu]_{\sqrt 3}\\
}}&
\xymatrix@R=0.5cm@C=0.6cm
{
&\\
\ar@{-}[r]^{\mu_1}&\\
}&
{\color{orange}\xymatrix@R=0.5cm@C=0.5cm
{
&1\ar[rd]^{\sqrt 3}&\\
2\ar[ru]^{2\sqrt 2}&&3\ar[ll]^{\sqrt 6} \\
}}
\\
&&
\xymatrix@R=0.5cm@C=0.8cm
{
&\ar@{-}[d]^{\mu_2}&\\
&&\\
}
&&
\xymatrix@R=0.5cm@C=0.8cm
{
&\ar@{-}[d]^{\mu_2}&\\
&&\\
}
&&
\\
&&
\xymatrix@R=0.5cm@C=0.5cm
{
&\times&\\
&&\\
}
&&
\xymatrix@R=0.5cm@C=0.5cm
{
&\times&\\
&&\\
}&&
\end{align*}
}
\caption{Mutation subgraph of $\Gamma_0$ of $(0,\sqrt{2},\sqrt{3})$ in direction $1$. Here and in the following, we use $\times$ to indicate that for any vertex $t$ lies in the branch $\times$, one has $B_{t}\not\sim \pm B_{t_0}$.}
    \label{fig:0sqrt2sqrt3-1}
\end{figure}

By mutating in direction $1$ first, we obtain the mutation subgraph of mutations of $\Gamma_0$, cf. Figure \ref{fig:0sqrt2sqrt3-1}.
It clear that $t_\diamond:=\mu_{31}(t_0)$ satisfies the condition $(\diamond)$.
By Proposition \ref{prop:generator-set-assumption}, Remark \ref{rem:factorization-trivial-path} and Remark \ref{rem:sharpen-K-2}, we deduce that $\operatorname{Aut}(\mathcal{A})$ is generated by
\begin{align*}
&g_{312}^{\mathbf{id},+}, g_{213}^{\mathbf{id},+}, g_{32}^{\mathbf{id},-}, g_{23}^{\mathbf{id},-}, g_{212}^{\mathbf{id},-}, g_{313}^{\mathbf{id},-}, g_{121}^{\mathbf{id},-}&\in H_{1,0}^{t_\diamond}(t_0),\\
& g_{13131}^{\mathbf{id},-}&\in H_{1,2}^{t_\diamond}(t_0).
\end{align*}
According to Lemma \ref{lem:relation-rank-2-finite}, we have
\[
    g_{121}^{\mathbf{id},-}\circ g_{212}^{\mathbf{id},-}=\mathbf{id}=g_{13131}^{\mathbf{id},-}\circ g_{313}^{\mathbf{id},-}.
\]
Thus, $\operatorname{Aut}(\mathcal{A})$ is generated by $g_{312}^{\mathbf{id},+}$ and $g_{32}^{\mathbf{id},-}$. Hence $\operatorname{Aut}(\mathcal{A})\cong D_\infty$.


\subsubsection{$a=0,b=\sqrt{2},c\geq 2$.}

The mutation subgraph of $\Gamma_0$ in direction $1$ is given in Figure \ref{fig:0sqrt2c-1}.
\begin{figure}
\centering
{\tiny
\begin{align*}
\hspace{1cm}
t_0:{\color{blue}\xymatrix@R=0.5cm@C=0.5cm
{
&1\ar[rd]^c&\\
2\ar[ru]^{\sqrt 2}&&3\\
}}
\xymatrix@R=0.5cm@C=0.6cm
{
&\\
\ar@{-}[r]^{\mu_1}&\\
}&
\xymatrix@R=0.5cm@C=0.5cm
{
&1\ar[ld]_{\sqrt 2}&\\
2\ar[rr]_{\sqrt 2 c}&&3\ar[lu]_{c}\\
}&
\xymatrix@R=0.5cm@C=0.6cm
{
&\\
\ar@{-}[r]^{\mu_2}&\\
}&
\xymatrix@R=0.5cm@C=0.5cm
{
&1\ar[rd]^c&\\
2\ar[ru]^{\sqrt 2}&&3\ar[ll]^{\sqrt2c} \\
}&
\xymatrix@R=0.5cm@C=0.6cm
{
&\\
\ar@{-}[r]^{\mu_1}&\\
}
{\color{blue}\xymatrix@R=0.5cm@C=0.5cm
{
&1\ar[ld]_{\sqrt 2}&\\
2&&3\ar[lu]_c\\
}}
\hspace{1cm}
\xymatrix@R=0.5cm@C=0.5cm
{
&&\\
&&\\
}\\
&
\xymatrix@R=0.5cm@C=0.8cm
{
&\ar@{-}[d]^{\mu_3}&\\
&&\\
}
&&
\xymatrix@R=0.5cm@C=0.8cm
{
&\ar@{-}[d]^{\mu_3}&\\
&&\\
}\\
&
\xymatrix@R=0.5cm@C=0.5cm
{
&\times&\\
&&\\
}
&&
\xymatrix@R=0.5cm@C=0.5cm
{
&\times&\\
&&\\
}
&
\end{align*}
}
 \caption{Mutation subgraph of $\Gamma_0$ of $(0,\sqrt{2},c\ge 2)$ in direction $1$.}
    \label{fig:0sqrt2c-1}
\end{figure}
It follows that $\operatorname{Aut}(\mathcal{A})$ is generated by 
\begin{align*}
g_{312}^{\bf{id},+},g_{32}^{\bf{id},-},g_{121}^{\mathbf{id},-} &\in H_{1}(t_0).
\end{align*}
Consequently, $\operatorname{Aut}(\mathcal{A})$ is generated by
$g_{312}^{\bf{id},+},g_{32}^{\bf{id},-}$ and hence $\operatorname{Aut}(\mathcal{A})\cong D_\infty$.

\subsubsection{$a=0,b=\sqrt{3},c=\sqrt{3}$.}

By mutating in direction $1$, we get the  mutation subgraph of $\Gamma_0$, cf. Figure \ref{fig:0sqrt3sqrt3-1}. Clearly, $t_\diamond:=\mu_{21}(t_0)$ satisfies the condition $(\diamond)$. 

\begin{figure}
\centering
{\tiny
 \begin{align*}
 \hspace{1cm}
 t_0:
 {\color{blue}\xymatrix@R=0.5cm@C=0.5cm
 {
 &1\ar[rd]^{\sqrt 3}&\\
 2\ar[ru]^{\sqrt 3}&&3\\
 }}
 &&&&&&
 {\color{blue}\xymatrix@R=0.5cm@C=0.5cm
 {
 &1\ar[ld]_{\sqrt 3}&\\
 2&&3\ar[lu]_{\sqrt 3}\\
 }}
 \\
 \xymatrix@R=0.5cm@C=0.5cm
 {
 &\ar@{-}[d]^{\mu_1}&\\
 && \\
 }
 &&&&&&
 \xymatrix@R=0.5cm@C=0.5cm
 {
 &\ar@{-}[d]^{\mu_1}&\\
 && \\
 }\\
 \xymatrix@R=0.5cm@C=0.5cm
 {
 &1\ar[ld]_{\sqrt 3}&\\
 2\ar[rr]_{3}&&3\ar[lu]_{\sqrt 3}\\
 }&
 \xymatrix@R=0.5cm@C=0.6cm
 {
 &\\
 \ar@{-}[r]^{\mu_2}&\\
 }&
 {\color{red}\xymatrix@R=0.5cm@C=0.5cm
 {
 &1\ar[rd]^{2\sqrt 3}&\\
 2\ar[ru]^{\sqrt 3}&&3\ar[ll]^{3} \\
 }}
 &
 \xymatrix@R=0.5cm@C=0.6cm
 {
 &\\
 \ar@{-}[r]^{\mu_1}&\\
 }&
 {\color{orange}\xymatrix@R=0.5cm@C=0.5cm
 {
 &1\ar[ld]_{\sqrt 3}&\\
 2\ar[rr]_{3}&&3\ar[lu]_{2 \sqrt 3}\\
 }}
 &
 \xymatrix@R=0.5cm@C=0.6cm
 {
 &\\
 \ar@{-}[r]^{\mu_2}&\\
 }
 &
 \xymatrix@R=0.5cm@C=0.5cm
 {
 &1\ar[rd]^{\sqrt 3}&\\
 2\ar[ru]^{\sqrt 3}&&3\ar[ll]^{3} \\
 }
 \\
 \xymatrix@R=0.5cm@C=0.6cm
 {
 &\ar@{-}[d]^{\mu_3}&\\
 &&\\
 }
 &&
 \xymatrix@R=0.5cm@C=0.6cm
 {
 &\ar@{-}[d]^{\mu_3}&\\
 &&\\
 }
 &&
 \xymatrix@R=0.5cm@C=0.6cm
 {
 &\ar@{-}[d]^{\mu_3}&\\
 &&\\
 }
 &&
 \xymatrix@R=0.5cm@C=0.6cm
 {
 &\ar@{-}[d]^{\mu_3}&\\
 &&\\
 }
 \\
 {\color{orange}\xymatrix@R=0.5cm@C=0.5cm
 {
 &1\ar[rd]^{\sqrt 3}&\\
 2\ar[ru]^{2\sqrt 3}&&3\ar[ll]^{3} \\
 }}
 &&
 \xymatrix@R=0.5cm@C=0.5cm
 {
 &\times& \\
 && \\
 }
 &&
 \xymatrix@R=0.5cm@C=0.5cm
 {
 &\times& \\
 && \\
 }
 &&
 {\color{orange}\xymatrix@R=0.5cm@C=0.5cm
 {
 &1\ar[ld]_{2\sqrt 3}&\\
 2\ar[rr]_{3}&&3\ar[ul]_{\sqrt{3}} \\
 }}
 \end{align*}
 }
 \caption{Mutation subgraph of $\Gamma_0$ of $(0,\sqrt{3},\sqrt{3})$ in direction $1$.}
    \label{fig:0sqrt3sqrt3-1}
 \end{figure}


Assume that $d_2=d_3$, then $|G_0(t_0)|=2$ and 
$\operatorname{Aut}(\mathcal{A})$ is generated by 
\begin{align*}
    &\psi_{(23)}^{-}&\in G_0(t_0),\\
    &g_{312}^{\mathbf{id},+}, g_{213}^{\mathbf{id},+}, g_{32}^{\mathbf{id},-}, g_{23}^{\mathbf{id},-}, g_{212}^{\mathbf{id},-}, g_{313}^{\mathbf{id},-} &\in H_{1,0}^{t_\diamond}(t_0),\\
& g_{12121}^{\mathbf{id},-}&\in K_{1}^{t_\diamond}(t_0).
\end{align*}

By Lemma \ref{lem:relation-rank-2-finite}, we obtain
\[
   g_{12121}^{\mathbf{id},-}\circ g_{212}^{\bf{id,-}}=\mathbf{id}.
\]
It follows that $\operatorname{Aut}(\mathcal{A})$ is generated by $g_{312}^{\bf{id},+},g_{32}^{\bf{id},-}$ and $\psi_{(23)}^{-}$. A direct calculation shows that 
\[
g_{312}^{\mathbf{id},+}\circ \psi_{(23)}^-=\psi_{(23)}^-\circ (g_{312}^{\mathbf{id},+})^{-1}, \psi_{(23)}^-\circ g_{32}^{\mathbf{id},-}=g_{32}^{\mathbf{id},-}\circ \psi_{(23)}^-.
\]
Similar to Section \ref{ss:a=0-b=1-c=sqrt3}, one can show that $\operatorname{Aut}(\mathcal{A})\cong D_\infty\rtimes_{\rho}\mathbb{Z}_2$.

Now assume that $d_2\neq d_3$, then $G_0(t_0)=\langle\mathbf{id}\rangle$. Note that every cluster automorphism in this case is also a cluster automorphism in the former case.
It is routine to check that $\operatorname{Aut}(\mathcal{A})$ is the subgroup of the cluster automorphsim group of the former case generated by $g_{312}^{\bf{id},+}$ and 
 $g_{32}^{\bf{id},-}$. Consequently, $\operatorname{Aut}(\mathcal{A})\cong D_\infty$.


\subsubsection{$a=0,b=\sqrt{3},c\geq 2$.}

By mutating in direction $1$, we obtain the  mutation subgraph of $\Gamma_0$ cf. Figure \ref{fig:0sqrt3c-1}. It follows that $\operatorname{Aut}(\mathcal{A})$ is generated by
\begin{figure}
\centering
{\tiny
\begin{align*}
\hspace{1cm}
t_0:
{\color{blue}\xymatrix@R=0.5cm@C=0.5cm
{
&1\ar[rd]^{c}&\\
2\ar[ru]^{\sqrt 3}&&3\\
}}&&&&&&
{\color{blue}\xymatrix@R=0.5cm@C=0.5cm
{
&1\ar[ld]_{\sqrt 3}&\\
2&&3\ar[lu]_{c}\\
}}
\\
\xymatrix@R=0.5cm@C=0.5cm
{
&\ar@{-}[d]^{\mu_1}&\\
&& \\
}
&&&&&&
\xymatrix@R=0.5cm@C=0.5cm
{
&\ar@{-}[d]^{\mu_1}&\\
&& \\
}
\\
\xymatrix@R=0.5cm@C=0.5cm
{
&1\ar[ld]_{\sqrt 3}&\\
2\ar[rr]_{\sqrt 3 c}&&3\ar[lu]_{c}\\
}&
\xymatrix@R=0.5cm@C=0.6cm
{
&\\
\ar@{-}[r]^{\mu_2}&\\
}&
\xymatrix@R=0.5cm@C=0.5cm
{
&1\ar[rd]^{2c}&\\
2\ar[ru]^{\sqrt 3}&&3\ar[ll]^{\sqrt 3 c} \\
}&
\xymatrix@R=0.5cm@C=0.6cm
{
&\\
\ar@{-}[r]^{\mu_1}&\\
}&
\xymatrix@R=0.5cm@C=0.5cm
{
&1\ar[ld]_{\sqrt 3}&\\
2\ar[rr]_{\sqrt 3c}&&3\ar[lu]_{2 c}\\
}
&
\xymatrix@R=0.5cm@C=0.6cm
{
&\\
\ar@{-}[r]^{\mu_2}&\\
}
&
\xymatrix@R=0.5cm@C=0.5cm
{
&1\ar[rd]^{c}&\\
2\ar[ru]^{\sqrt 3}&&3\ar[ll]^{\sqrt 3 c} \\
}
\\
\xymatrix@R=0.5cm@C=0.5cm
{
&\ar@{-}[d]^{\mu_3}&\\
&& \\
}
&
\xymatrix@R=0.5cm@C=0.6cm
{
&\\
&\\
}
&
\xymatrix@R=0.5cm@C=0.5cm
{
&\ar@{-}[d]^{\mu_3}&\\
&& \\
}
&
\xymatrix@R=0.5cm@C=0.6cm
{
&\\
&\\
}
&
\xymatrix@R=0.5cm@C=0.5cm
{
&\ar@{-}[d]^{\mu_3}&\\
&& \\
}
&
\xymatrix@R=0.5cm@C=0.6cm
{
&\\
&\\
}
&
\\
\xymatrix@R=0.5cm@C=0.5cm
{
&\times&\\
&& \\
}
&
\xymatrix@R=0.5cm@C=0.6cm
{
&\\
&\\
}&
\xymatrix@R=0.5cm@C=0.5cm
{
&\times& \\
&& \\
}&
\xymatrix@R=0.5cm@C=0.6cm
{
&\\
&\\
}&
\xymatrix@R=0.5cm@C=0.5cm
{
&\times&\\
&&\\
}&
\xymatrix@R=0.5cm@C=0.6cm
{
&\\
&\\
}&
\end{align*}
}
\caption{Mutation subgraph of $\Gamma_0$ of $(0,\sqrt{3},c\ge 2)$ in direction $1$.}
    \label{fig:0sqrt3c-1}
\end{figure}


\[
 g_{312}^{\mathbf{id},+}, g_{213}^{\mathbf{id},+}, g_{32}^{\mathbf{id},-}, g_{23}^{\mathbf{id},-}, g_{212}^{\mathbf{id},-}, g_{313}^{\mathbf{id},-},
 g_{12121}^{\mathbf{id},-}
 \in H_{1}(t_0)
\]
Moreover, we have $g_{12121}^{\mathbf{id},-}\circ g_{212}^{\mathbf{id},-}=\mathbf{id}$ by Lemma \ref{lem:relation-rank-2-finite}. Consequently, $\operatorname{Aut}(\mathcal{A})\cong D_\infty$.


\subsubsection{$a=0,b\geq 2,c\geq 2$.}
By mutating in direction $1$, we obtain
\begin{align*}
\xymatrix@R=0.5cm@C=0.5cm
{
&1\ar[rd]^{c}&\\
2\ar[ru]^{b}&&3\\
}
\xymatrix@R=0.5cm@C=0.6cm
{
&\\
\ar@{-}[r]^{\mu_1}&\\
}
\xymatrix@R=0.5cm@C=0.5cm
{
&&\\
\times.&& \\
} 
\end{align*}
If $b=c$ and $d_2=d_3$, then $|G_0(t_0)|=2$ and $\psi_{(23)}^-\in G_0(t_0)$. Moreover, 
\[
g_{312}^{\mathbf{id},+}\circ \psi_{(23)}^-=\psi_{(23)}^-\circ (g_{312}^{\mathbf{id},+})^{-1}, \psi_{(23)}^-\circ g_{32}^{\mathbf{id},-}=g_{32}^{\mathbf{id},-}\circ \psi_{(23)}^-.
\]
It follows that $\operatorname{Aut}(\mathcal{A})\cong D_\infty\rtimes_{\rho}\mathbb{Z}_2$ in this case. 
Otherwise, $\operatorname{Aut}(\mathcal{A})\cong D_\infty$.


\subsection{Case 2: $a=1$.}

Since $b_{12}b_{23}b_{31}=b_{21}b_{32}b_{13}$, $abc=b_{12}b_{23}b_{32}$ is a positive integer.

\subsubsection{$a=b=c=1$.} In this case, $d_1=d_2=d_3$. It turns out that $\mathcal{A}$ is of Euclidean type $\tilde{A}_{2,1}$. Hence, $\operatorname{Aut}(\mathcal{A})\cong D_\infty$, cf. \cite{ASS2012}.

\subsubsection{$a=1,b=1,c\geq 2$.}
Clearly, $d_1=d_2=d_3$ and $G_0(t_0)=\langle\mathbf{id}\rangle$.
Mutation subgraphs of $\Gamma_0$ are given in Figure \ref{fig:11c-1} and \ref{fig:11c-23}. It is obvious that $t_\diamond=\mu_1(t_0)$ satisfies the condition $(\diamond)$.
\begin{figure}
\centering
{\tiny\bf
\begin{align*}
&
&
\xymatrix@R=0.5cm@C=0.5cm
{
&&\\
&\times& \\
}
&
&
&
&
{\color{blue}\xymatrix@R=0.5cm@C=0.5cm
{
&1&\\
2\ar[ru]^1\ar[rr]_c&&3\ar[lu]_1\\
}}
\\
&&
\xymatrix@R=0.5cm@C=0.8cm
{
&\ar@{-}[d]_{\mu_3}&\\
&&\\
}
&
&
&
&
\xymatrix@R=0.5cm@C=0.8cm
{
&\ar@{-}[d]_{\mu_2}&\\
&&\\
}\\
t_0:
{\color{blue}\xymatrix@R=0.5cm@C=0.5cm
{
&1\ar[rd]^c&\\
2\ar[ru]^1\ar[rr]_1&&3\\
}}
&
\xymatrix@R=0.5cm@C=0.6cm
{
&\\
\ar@{-}[r]^{\mu_1}&\\
}
&
{\color{red}\xymatrix@R=0.5cm@C=0.5cm
{
&1\ar[ld]_1&\\
2\ar[rr]_{1+c}&&3\ar[lu]_c\\
}}
&
\xymatrix@R=0.5cm@C=0.6cm
{
&\\
\ar@{-}[r]^{\mu_2}&\\
}
&
\xymatrix@R=0.5cm@C=0.5cm
{
&1\ar[rd]^1&\\
2\ar[ru]^1&&3\ar[ll]^{1+c}\\
}
&
\xymatrix@R=0.5cm@C=0.6cm
{
&\\
\ar@{-}[r]^{\mu_1}&\\
}
&
\xymatrix@R=0.5cm@C=0.5cm
{
&1\ar[ld]_1&\\
2&&3\ar[ll]^c\ar[lu]_1\\
}
\\
\xymatrix@R=0.5cm@C=0.8cm
{
&&\\
&&\\
}
&
\xymatrix@R=0.5cm@C=0.6cm
{
&\\
&\\
}
&
\xymatrix@R=0.5cm@C=0.6cm
{
&\\
&\\
}
&
\xymatrix@R=0.5cm@C=0.6cm
{
&\\
&\\
}
&
\xymatrix@R=0.5cm@C=0.8cm
{
&\ar@{-}[d]_{\mu_3}&\\
&&\\
}
&
\xymatrix@R=0.5cm@C=0.6cm
{
&\\
&\\
}
&
\xymatrix@R=0.5cm@C=0.8cm
{
&\ar@{-}[d]_{\mu_3}&\\
&&\\
}
\\
&
&
{\color{blue}\xymatrix@R=0.5cm@C=0.5cm
{
&1\ar[rd]^1&\\
2\ar[ru]^{c}\ar[rr]_1&&3\\
}}
&
\xymatrix@R=0.5cm@C=0.6cm
{
&\\
\ar@{-}[r]^{\mu_1}&\\
}
&
{\color{orange}\xymatrix@R=0.5cm@C=0.5cm
{
&1\ar[ld]_c&\\
2\ar[rr]_{1+c}&&3\ar[lu]_1\\
}}
&
\xymatrix@R=0.5cm@C=0.6cm
{
&\\
&\\
}
&
{\color{blue}\xymatrix@R=0.5cm@C=0.5cm
{
&1\ar[rd]^1\ar[ld]_1&\\
2\ar[rr]_c&&3\\
}}
\end{align*}
}
\caption{Mutation subgraph of $\Gamma_0$ of $(1,1,c\ge 2)$ in direction $1$.}
    \label{fig:11c-1}
\end{figure}

\begin{figure}
\centering
{\tiny\bf 
\begin{align*}
t_0:
{\color{blue}\xymatrix@R=0.5cm@C=0.5cm
{
&1\ar[rd]^c&\\
2\ar[ru]^1\ar[rr]_1&&3\\
}}&
\xymatrix@R=0.5cm@C=0.6cm
{
&\\
\ar@{-}[r]^{\mu_3}&\\
}&
\xymatrix@R=0.5cm@C=0.5cm
{
&1&\\
2\ar[ru]^{1}&&3\ar[lu]_c\ar[ll]^1\\
}&
\xymatrix@R=0.5cm@C=0.6cm
{
&\\
\ar@{-}[r]^{\mu_2}&\\
}&
\xymatrix@R=0.5cm@C=0.5cm
{
&1\ar[ld]_1&\\
2\ar[rr]_1&&3\ar[lu]_{1+c}\\
}&
\xymatrix@R=0.5cm@C=0.6cm
{
&\\
\ar@{-}[r]^{\mu_1}&\\
}&
{\color{orange}\xymatrix@R=0.5cm@C=0.5cm
{
&1\ar[rd]^{1+c}&\\
2\ar[ru]^{1}&&3\ar[ll]^c\\
}}
&
\xymatrix@R=0.5cm@C=0.6cm
{
&\\
&\\
}
&
\xymatrix@R=0.5cm@C=0.5cm
{
&&\\
&& \\
}
\\
\xymatrix@R=0.5cm@C=0.8cm
{
&\ar@{-}[d]_{\mu_2}&\\
&&\\
}
&
\xymatrix@R=0.5cm@C=0.6cm
{
&\\
&\\
}
&
\xymatrix@R=0.5cm@C=0.8cm
{
&\ar@{-}[d]_{\mu_1}&\\
&&\\
}
&
\xymatrix@R=0.5cm@C=0.6cm
{
&\\
&\\
}
&
\xymatrix@R=0.5cm@C=0.8cm
{
&\ar@{-}[d]_{\mu_3}&\\
&&\\
}
&
\xymatrix@R=0.5cm@C=0.6cm
{
&\\
&\\
}
&
\xymatrix@R=0.5cm@C=0.8cm
{
&& \\
&&\\
}
&
\xymatrix@R=0.5cm@C=0.6cm
{
&\\
& \\
}
&
\xymatrix@R=0.5cm@C=0.8cm
{
&&\\
&& \\
}
\\
{\color{blue}\xymatrix@R=0.5cm@C=0.5cm
{
&1\ar[rd]^c\ar[ld]_1&\\
2&&3\ar[ll]^{1}\\
}}
&
\xymatrix@R=0.5cm@C=0.6cm
{
&\\
&\\
}
&
{\color{blue}\xymatrix@R=0.5cm@C=0.5cm
{
&1\ar[rd]^c\ar[ld]_1&\\
2&&3\ar[ll]^{1}\\
}}
&
\xymatrix@R=0.5cm@C=0.6cm
{
&\\
&\\
}
&
{\color{orange}\xymatrix@R=0.5cm@C=0.5cm
{
&1\ar[rd]^{1+c}&\\
2\ar[ru]^{c}&&3\ar[ll]^1\\
}}
&
\xymatrix@R=0.5cm@C=0.6cm
{
&\\
&\\
}
&
\xymatrix@R=0.5cm@C=0.5cm
{
&&\\
&&\\
}
&
\xymatrix@R=0.5cm@C=0.6cm
{
&\\
&\\
}
&
\xymatrix@R=0.5cm@C=0.5cm
{
&&\\
&&\\
}
\end{align*}
}
\caption{Mutation subgraph of $\Gamma_0$ of $(1,1,c\ge 2)$ in direction $2$ and $3$.}
    \label{fig:11c-23}
\end{figure}


According to Proposition \ref{prop:generator-set-assumption},
$\operatorname{Aut}(\mathcal{A})$ is generated by
\begin{align*}
& g_{2}^{(1 3),-},g_{13}^{(1 3),-}
&\in H_{1,0}^{t_\diamond}(t_0), \\
& g_{2121}^{(1 3 2),-},
g_{3121}^{(1 2),+},
g_{2123}^{(1 2),+},
g_{2323}^{(1 2 3),-}
&\in H_{1,1}^{t_\diamond}(t_0),\\
& g_{1321}^{(2 3),-} 
&\in K_{1}^{t_\diamond}(t_0).
\end{align*}
By Lemma \ref{lem:relation-rank-2-finite}, we have
\[
    g_{2121}^{(1 3 2),-}\circ g_2^{(1 3),-}=\mathbf{id}=g_2^{(1 3),-}\circ g_{2323}^{(1 2 3),-}.
\]
On the other hand, by direct computation, we obtain
\begin{align*}
     &g_{3121}^{(1 2),+}\circ g_{13}^{(1 3),-}=g_{2121}^{(1 3 2),-},& g_{13}^{(1 3),-}\circ g_{2123}^{(1 2),+}=g_{2323}^{(1 2 3),-}, \\
     &g_{1321}^{(2 3),-}\circ g_{3121}^{(1 2),+}=g_{2121}^{(1 3 2),-},& (g_{2}^{(1 3),-})^2=(g_{13}^{(1 3),-})^2=\mathbf{id}.  
\end{align*}
It follows that $\operatorname{Aut}(\mathcal{A})$ is generated by $g_2^{(13),-}$ and $g_{13}^{(13),-}$. We claim that $\operatorname{Aut}(\mathcal{A})\cong D_\infty$. It suffices to show that \[(g_{2}^{(1 3),-}\circ g_{13}^{(1 3),-})^m\neq \mathbf{id}\  \text{and}\  g_{13}^{(1 3),-} \circ (g_{2}^{(1 3),-}\circ g_{13}^{(1 3),-})^m\neq \mathbf{id}\] for any positive integer $m$. Note that $g_{2}^{(1 3),-}\circ g_{13}^{(1 3),-}=g_{312}^{\bf{id},+}$. By Lemma \ref{lem:order-of-tau}, we know that the order of $g_{312}^{\bf{id},+}$ is infinite, which implies the result.

\subsubsection{$a=1,b=\sqrt{2},c=\sqrt{2}$.}
 It is clear that $d_2=d_3$ and $G_0(t_0)=\{\mathbf{id},\psi_{(23)}^-\}$. Mutation sugraphs of $\Gamma_0$ in directions $1,2,3$ are given in Figure \ref{fig:1sqrt2sqrt2-1} and \ref{fig:1sqrt2sqrt2-23}. It is clear that $t_\diamond:=\mu_1\mu_2\mu_1(t_0)$ satisfies the condition $(\diamond)$. 
\begin{figure}[t]
      \centering
      {\tiny\bf 
\begin{align*}
&&&&
{\color{orange}\xymatrix@R=0.5cm@C=0.5cm
{
&1\ar[ld]_{2\sqrt 2}&\\
2\ar[rr]_{1}&&3\ar[lu]_{\sqrt 2}\\
}}
&&
{\color{blue}\xymatrix@R=0.5cm@C=0.5cm
{
&1\ar[rd]^{\sqrt 2}&\\
2\ar[ru]^{\sqrt 2}\ar[rr]_1&&3 \\
}}
\\
&&&&
\xymatrix@R=0.5cm@C=0.8cm
{
&\ar@{-}[d]^{\mu_3}&\\
&&\\
}
&&
\xymatrix@R=0.5cm@C=0.8cm
{
&\ar@{-}[d]^{\mu_3}&\\
&&\\
}
\\
&&
{\color{blue}\xymatrix@R=0.5cm@C=0.5cm
{
&1\ar[rd]^{\sqrt 2}&\\
2\ar[ru]^{\sqrt 2}\ar[rr]_1&&3\\
}}
&
\xymatrix@R=0.5cm@C=0.6cm
{
&\\
\ar@{-}[r]^{\mu_2}&\\
}
&
\xymatrix@R=0.5cm@C=0.5cm
{
&1\ar[ld]_{\sqrt 2}\ar[rd]^{\sqrt 2}&\\
2&&3\ar[ll]^{1}\\
}
&
\xymatrix@R=0.5cm@C=0.6cm
{
&\\
\ar@{-}[r]^{\mu_1}&\\
}
&
\xymatrix@R=0.5cm@C=0.5cm
{
&1&\\
2\ar[ru]^{\sqrt{2}}&&3\ar[ll]^{1}\ar[lu]_{\sqrt{2}}\\
}\\
&&&&&&
\xymatrix@R=0.5cm@C=0.8cm
{
&\ar@{-}[d]^{\mu_2}&\\
&&\\
}\\
t_0:{\color{blue}\xymatrix@R=0.5cm@C=0.5cm
{
&1\ar[rd]^{\sqrt 2}&\\
2\ar[ru]^{\sqrt 2}\ar[rr]_1&&3\\
}}
&
\xymatrix@R=0.5cm@C=0.6cm
{
&\\
\ar@{-}[r]^{\mu_1}&\\
}
&
\xymatrix@R=0.5cm@C=0.5cm
{
&1\ar[ld]_{\sqrt 2}&\\
2\ar[rr]_{3}&&3\ar[lu]_{\sqrt 2}\\
}
&
\xymatrix@R=0.5cm@C=0.6cm
{
&\\
\ar@{-}[r]^{\mu_2}&\\
}
&
\xymatrix@R=0.5cm@C=0.5cm
{
&1\ar[rd]^{2\sqrt 2}&\\
2\ar[ru]^{\sqrt 2}&&3\ar[ll]^3\\
}
&
\xymatrix@R=0.5cm@C=0.6cm
{
&\\
\ar@{-}[r]^{\mu_1}&\\
}
&
{\color{red}\xymatrix@R=0.5cm@C=0.5cm
{
&1\ar[ld]_{\sqrt 2}&\\
2\ar[rr]_{1}&&3\ar[lu]_{2\sqrt 2}\\
}}
&
\xymatrix@R=0.5cm@C=0.6cm
{
&\\
&\\
}
&
\xymatrix@R=0.5cm@C=0.5cm
{
&&\\
&& \\
}
\\
\xymatrix@R=0.5cm@C=0.8cm
{
&&\\
&&\\
}
&
\xymatrix@R=0.5cm@C=0.5cm
{
&\\
& \\
}
&
\xymatrix@R=0.5cm@C=0.8cm
{
&\ar@{-}[d]^{\mu_3}&\\
&&\\
}
&
\xymatrix@R=0.5cm@C=0.5cm
{
&\\
& \\
}
&
\xymatrix@R=0.5cm@C=0.8cm
{
&\ar@{-}[d]^{\mu_3}&\\
&&\\
}
&
\xymatrix@R=0.5cm@C=0.5cm
{
&\\
& \\
}
&
\xymatrix@R=0.5cm@C=0.8cm
{
&\ar@{-}[d]^{\mu_3}&\\
&&\\
}
&
\xymatrix@R=0.5cm@C=0.5cm
{
&\\
& \\
}
&
\xymatrix@R=0.5cm@C=0.8cm
{
&&\\
&&\\
}
\\
{\color{orange}\xymatrix@R=0.5cm@C=0.5cm
{
&1\ar[ld]_{2\sqrt 2}&\\
2\ar[rr]_{1}&&3\ar[lu]_{\sqrt 2}\\
}}
&
\xymatrix@R=0.5cm@C=0.6cm
{
&\\
\ar@{-}[r]^{\mu_1}&\\
}
&
\xymatrix@R=0.5cm@C=0.5cm
{
&1\ar[rd]^{\sqrt 2}&\\
2\ar[ru]^{2\sqrt 2}&&3\ar[ll]^3\\
}
&
\xymatrix@R=0.5cm@C=0.6cm
{
&\\
\ar@{-}[r]^{\mu_2}&\\
}
&
\xymatrix@R=0.5cm@C=0.5cm
{
&\times&\\
\times&&\\
}
&
\xymatrix@R=0.5cm@C=0.6cm
{
&\\
&\\
}
&
{\color{orange}\xymatrix@R=0.5cm@C=0.5cm
{
&1\ar[rd]^{2\sqrt 2}&\\
2\ar[ru]^{\sqrt 2}&&3\ar[ll]^1\\
}}
&
\xymatrix@R=0.5cm@C=0.6cm
{
&\\
&\\
}
&
\xymatrix@R=0.5cm@C=0.5cm
{
&&\\
&&\\
}
\end{align*}
}
      \caption{Mutation subgraph of $\Gamma_0$ of $(1,\sqrt{2},\sqrt{2})$ in direction $1$.}
      \label{fig:1sqrt2sqrt2-1}
  \end{figure}

\begin{figure}[ht]
    \centering 
   {\tiny\bf 
\begin{align*}
\xymatrix@R=0.5cm@C=0.5cm
{
&&\\
&&\\
}
&
\xymatrix@R=0.5cm@C=0.6cm
{
&\\
&\\
}
&
{\color{orange}\xymatrix@R=0.5cm@C=0.5cm
{
&1\ar[ld]_{2\sqrt 2}&\\
2\ar[rr]_{1}&&3\ar[lu]_{\sqrt 2}\\
}}
&
\xymatrix@R=0.5cm@C=0.6cm
{
&\\
&\\
}
&
{\color{orange}\xymatrix@R=0.5cm@C=0.5cm
{
&1\ar[ld]_{\sqrt 2}&\\
2\ar[rr]_{1}&&3\ar[lu]_{2\sqrt 2}\\
}}
&
\xymatrix@R=0.5cm@C=0.6cm
{
&\\
&\\
}
&
\xymatrix@R=0.5cm@C=0.5cm
{
&&\\
&&\\
}
&
\xymatrix@R=0.5cm@C=0.6cm
{
&\\
&\\
}
&
\xymatrix@R=0.5cm@C=0.5cm
{
&&\\
&&\\
}\\
\xymatrix@R=0.5cm@C=0.8cm
{
&&\\
&&\\
}
&
\xymatrix@R=0.5cm@C=0.6cm
{
&\\
&\\
}
&
\xymatrix@R=0.5cm@C=0.8cm
{
&\ar@{-}[d]^{\mu_3}&\\
&&\\
}
&
\xymatrix@R=0.5cm@C=0.6cm
{
&\\
&\\
}
&
\xymatrix@R=0.5cm@C=0.8cm
{
&\ar@{-}[d]^{\mu_2}&\\
&&\\
}
&
\xymatrix@R=0.5cm@C=0.6cm
{
&\\
&\\
}
&
\xymatrix@R=0.5cm@C=0.8cm
{
&&\\
&&\\
}
&
\xymatrix@R=0.5cm@C=0.6cm
{
&\\
&\\
}
&\xymatrix@R=0.5cm@C=0.8cm
{
&&\\
&&\\
}\\
t_0:{\color{blue}\xymatrix@R=0.5cm@C=0.5cm
{
&1\ar[rd]^{\sqrt 2}&\\
2\ar[ru]^{\sqrt 2}\ar[rr]_1&&3\\
}}
&
\xymatrix@R=0.5cm@C=0.6cm
{
&\\
\ar@{-}[r]^{\mu_2}&\\
}
&
\xymatrix@R=0.5cm@C=0.5cm
{
&1\ar[ld]_{\sqrt 2}\ar[rd]^{\sqrt 2}&\\
2&&3\ar[ll]^{1}\\
}
&
\xymatrix@R=0.5cm@C=0.6cm
{
&\\
\ar@{-}[r]^{\mu_1}&\\
}
&
\xymatrix@R=0.5cm@C=0.5cm
{
&1&\\
2\ar[ru]^{\sqrt 2}&&3\ar[ll]^1\ar[lu]_{\sqrt 2}\\
}
&
\xymatrix@R=0.5cm@C=0.6cm
{
&\\
\ar@{-}[r]^{\mu_3}&\\
}
&
{\color{blue}\xymatrix@R=0.5cm@C=0.5cm
{
&1\ar[rd]^{\sqrt 2}&\\
2\ar[ru]^{\sqrt 2}\ar[rr]_1&&3\\
}}
&
\xymatrix@R=0.5cm@C=0.6cm
{
&\\
&\\
}
&
\xymatrix@R=0.5cm@C=0.8cm
{
&&\\
&&\\
}\\
\xymatrix@R=0.5cm@C=0.8cm
{
&\ar@{-}[d]^{\mu_3}&\\
&&\\
}
&&&&&&\\
\xymatrix@R=0.5cm@C=0.5cm
{
&1&\\
2\ar[ru]^{\sqrt 2}&&3\ar[ll]^{1}\ar[lu]_{\sqrt 2}\\
}
&
\xymatrix@R=0.5cm@C=0.6cm
{
&\\
\ar@{-}[r]^{\mu_1}&\\
}
&
\xymatrix@R=0.5cm@C=0.5cm
{
&1\ar[ld]_{\sqrt 2}\ar[rd]^{\sqrt 2}&\\
2&&3\ar[ll]^{1}\\
}
&
\xymatrix@R=0.5cm@C=0.6cm
{
&\\
\ar@{-}[r]^{\mu_2}&\\
}
&
{\color{blue}\xymatrix@R=0.5cm@C=0.5cm
{
&1\ar[rd]^{\sqrt 2}&\\
2\ar[ru]^{\sqrt 2}\ar[rr]_1&&3\\
}}
&
\xymatrix@R=0.5cm@C=0.6cm
{
&\\
&\\
}
&
\xymatrix@R=0.5cm@C=0.5cm
{
&&\\
&& \\
}
&&
\\
\xymatrix@R=0.5cm@C=0.8cm
{
&\ar@{-}[d]^{\mu_2}&\\
&&\\
}
&
\xymatrix@R=0.5cm@C=0.6cm
{
&\\
&\\
}
&
\xymatrix@R=0.5cm@C=0.8cm
{
&\ar@{-}[d]^{\mu_3}&\\
&&\\
}
&&&&
\\
{\color{orange}\xymatrix@R=0.5cm@C=0.5cm
{
&1\ar[ld]_{\sqrt 2}&\\
2\ar[rr]_{1}&&3\ar[lu]_{2\sqrt 2}\\
}}
&
\xymatrix@R=0.5cm@C=0.6cm
{
&\\
&\\
}
&
{\color{orange}\xymatrix@R=0.5cm@C=0.5cm
{
&1\ar[ld]_{2\sqrt 2}&\\
2\ar[rr]_{1}&&3\ar[lu]_{\sqrt 2}\\
}}&&&&
\end{align*}
}
 \caption{Mutation subgraph of $\Gamma_0$  of  $(1,\sqrt{2},\sqrt{2})$ in directions $2$ and $3$.}
    \label{fig:1sqrt2sqrt2-23}
\end{figure}



According to Proposition \ref{prop:generator-set-assumption} and Remark \ref{rem:sharpen-K-2}, we conclude that $\operatorname{Aut}(\mathcal{A})$ is generated by
\begin{align*}
&   \psi_{(23)}^-&\in G_0(t_0),\\
&    g_{312}^{\bf{id},+}, 
     g_{213}^{\bf{id},+}
    &\in H_{1,0}^{t_\diamond}(t_0),\\
 &   g_{212121}^{\bf{id},+},
    g_{32121}^{\bf{id},+},
    g_{313131}^{\bf{id},+},
    g_{23131}^{\bf{id},+},
    g_{13132}^{\bf{id},+},
    g_{121212}^{\bf{id},+},
    g_{12123}^{\bf{id},+},
    g_{131313}^{\bf{id},+},
    &\in H_{1,1}^{t_\diamond}(t_0),\\
&    g_{323121}^{(2 3),+},
    g_{131312121}^{\bf{id},+}
    &\in K_{1}^{t_\diamond}(t_0),\\
\end{align*}
By Lemma \ref{lem:relation-rank-2-finite},
\[
    g_{121212}^{\bf{id},+}=g_{212121}^{\bf{id},+}=\mathbf{id}, g_{131313}^{\bf{id},+}=g_{313131}^{\bf{id},+}=\mathbf{id}, g_{12123}^{\bf{id},+}\circ g_{323121}^{(2 3),+}=\mathbf{id}.
\]
Moreover, we have
\begin{align*}
 &g_{212121}^{\bf{id},+}=g_{32121}^{\bf{id},+}\circ g_{213}^{\bf{id},+}, g_{313131}^{\bf{id},+}=g_{23131}^{\bf{id},+}\circ g_{312}^{\bf{id},+}, \\
 & g_{131313}^{\bf{id},+}=g_{213}^{\bf{id},+}\circ g_{13132}^{\bf{id},+}, g_{121212}^{\bf{id},+}=g_{312}^{\bf{id},+}\circ g_{12123}^{\bf{id},+},\\
  &g_{131312121}^{\bf{id},+}\circ g_{23131}^{\bf{id},+}=g_{212121}^{\bf{id},+},
g_{312}^{\bf{id},+}\circ g_{213}^{\bf{id},+}=\mathbf{id}, \\
&g_{312}^{\bf{id},+}\circ \psi_{(23)}^{-}=\psi_{(23)}^-\circ g_{213}^{\bf{id},+}.
\end{align*}
As a consequence, $\operatorname{Aut}(\mathcal{A})$ is generated by $g_{312}^{\mathbf{id},+}$ and $\psi_{(23)}^-$. By Lemma \ref{lem:order-of-tau}, the order of $g_{312}^{\mathbf{id},+}$ is infinity which implies that for any $i\in \mathbb{Z}$ and $0\leq j\leq 1$, $(g_{312}^{\bf{id},+})^i\circ(\psi_{(23)}^-)^j=\mathbf{id}$ if and only if $i=j=0$. Hence $\operatorname{Aut}(\mathcal{A})\cong D_\infty$.

\subsubsection{$a=1,b=\sqrt{2},c\geq 2$.}
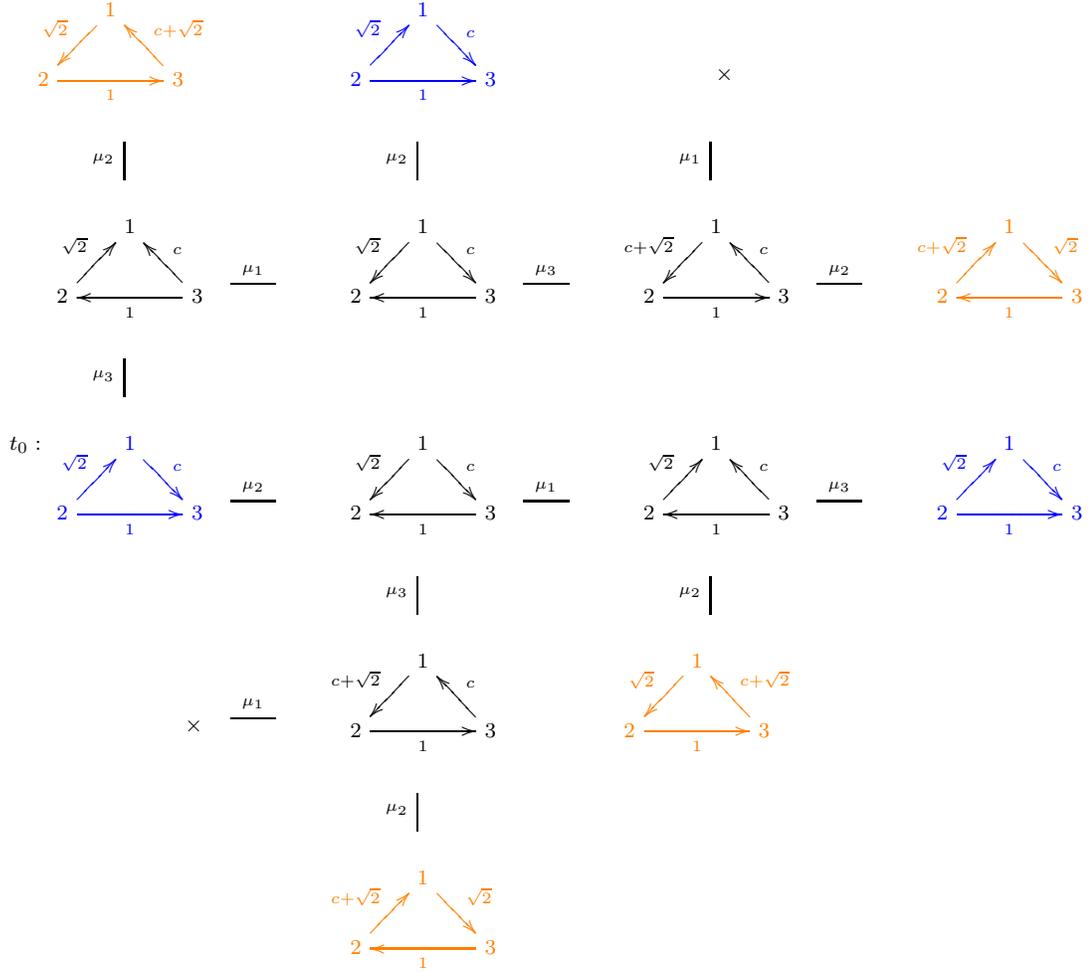
\begin{figure}
\centering
{\tiny\bf
\begin{align*}
{\color{orange}\xymatrix@R=0.5cm@C=0.5cm
{
&1\ar[ld]_{\sqrt{2}}&\\
2\ar[rr]_{1}&&3\ar[ul]_{c+\sqrt{2}}\\
}}
&
&
{\color{blue}\xymatrix@R=0.5cm@C=0.5cm
{
&1\ar[rd]^{c}&\\
2\ar[ru]^{\sqrt 2}\ar[rr]_1&&3\\
}}
&
&
\xymatrix@R=0.5cm@C=0.5cm
{
&&\\
&\times& \\
}
&&\\
\xymatrix@R=0.5cm@C=0.8cm
{
&\ar@{-}[d]_{\mu_2}& \\
&&\\
}
&
&
\xymatrix@R=0.5cm@C=0.8cm
{
&\ar@{-}[d]_{\mu_2}& \\
&&\\
}
&
&
\xymatrix@R=0.5cm@C=0.8cm
{
&\ar@{-}[d]_{\mu_1}& \\
&&\\
}
&
&
\\
\xymatrix@R=0.5cm@C=0.5cm
{
&1&\\
2\ar[ru]^{\sqrt 2}&&3\ar[ll]^1\ar[ul]_c\\
}
&
\xymatrix@R=0.5cm@C=0.6cm
{
&\\
\ar@{-}[r]^{\mu_1}&\\
}
&
\xymatrix@R=0.5cm@C=0.5cm
{
&1\ar[dl]_{\sqrt{2}}\ar[dr]^c&\\
2&&3\ar[ll]^1\\
}
&
\xymatrix@R=0.5cm@C=0.6cm
{
&\\
\ar@{-}[r]^{\mu_3}&\\
}
&
\xymatrix@R=0.5cm@C=0.5cm
{
&1\ar[ld]_{c+\sqrt{2}}&\\
2\ar[rr]_{1}&&3\ar[ul]_{c}\\
}
&
\xymatrix@R=0.5cm@C=0.6cm
{
&\\
\ar@{-}[r]^{\mu_2}&\\
}
&
{\color{orange}\xymatrix@R=0.5cm@C=0.5cm
{
&1\ar[rd]^{\sqrt{2}}&\\
2\ar[ru]^{c+\sqrt 2}&&3\ar[ll]^1\\
}}
\\
\xymatrix@R=0.5cm@C=0.8cm
{
&\ar@{-}[d]_{\mu_3}& \\
&&\\
}&
&&&&&
\\
t_0:{\color{blue}\xymatrix@R=0.5cm@C=0.5cm
{
&1\ar[rd]^{c}&\\
2\ar[ru]^{\sqrt 2}\ar[rr]_1&&3\\
}}
&
\xymatrix@R=0.5cm@C=0.6cm
{
&\\
\ar@{-}[r]^{\mu_2}&\\
}
&
\xymatrix@R=0.5cm@C=0.5cm
{
&1\ar[dl]_{\sqrt{2}}\ar[dr]^c&\\
2&&3\ar[ll]^1\\
}
&
\xymatrix@R=0.5cm@C=0.6cm
{
&\\
\ar@{-}[r]^{\mu_1}&\\
}
&
\xymatrix@R=0.5cm@C=0.5cm
{
&1&\\
2\ar[ur]^{\sqrt{2}}&&3\ar[ll]^1\ar[ul]_c\\
}
&
\xymatrix@R=0.5cm@C=0.6cm
{
&\\
\ar@{-}[r]^{\mu_3}&\\
}
&
{\color{blue}\xymatrix@R=0.5cm@C=0.5cm
{
&1\ar[rd]^{c}&\\
2\ar[ru]^{\sqrt 2}\ar[rr]_1&&3\\
}}
\\
&&
\xymatrix@R=0.5cm@C=0.8cm
{
&\ar@{-}[d]_{\mu_3}& \\
&&\\
}
&&
\xymatrix@R=0.5cm@C=0.8cm
{
&\ar@{-}[d]_{\mu_2}& \\
&&\\
}
&&
\\
\xymatrix@R=0.5cm@C=0.5cm
{
&& \\
&&\times \\
}&
\xymatrix@R=0.5cm@C=0.6cm
{
&\\
\ar@{-}[r]^{\mu_1}&\\
}
&
\xymatrix@R=0.5cm@C=0.5cm
{
&1\ar[ld]_{c+\sqrt{2}}&\\
2\ar[rr]_{1}&&3\ar[ul]_c\\
}
&
&
{\color{orange}\xymatrix@R=0.5cm@C=0.5cm
{
&1\ar[ld]_{\sqrt{2}}&\\
2\ar[rr]_1&&3\ar[lu]_{c+\sqrt{2}}\\
}}
&&
\\
&&
\xymatrix@R=0.5cm@C=0.8cm
{
&\ar@{-}[d]_{\mu_2}& \\
&&\\
}
&&
&&
\\
&&
{\color{orange}\xymatrix@R=0.5cm@C=0.5cm
{
&1\ar[rd]^{\sqrt{2}}&\\
2\ar[ru]^{c+\sqrt 2}&&3\ar[ll]^1\\
}}
&&
&
&
\end{align*}
}
\caption{Mutation subgraph of $\Gamma_0$ of $(1,\sqrt{2},c\geq 2)$ in directions $2$ and $3$.}
\label{fig:1sqrt2c-23}
\end{figure}

Clearly, $d_2=d_3$ and $G_0(t_0)=\langle\mathbf{id}\rangle$.  Mutation subgraphs of $\Gamma_0$ in directions $1,2,3$ are given in Figure \ref{fig:1sqrt2c-1} and \ref{fig:1sqrt2c-23}. It is clear that $t_\diamond:=\mu_{121}(t_0)$ satisfies the condition $(\diamond)$.



\begin{figure}
    \centering
    {\tiny\bf
\begin{align*}
&&
{\color{orange}\xymatrix@R=0.5cm@C=0.5cm
{
&1\ar[rd]^{\sqrt 2}&\\
2\ar[ur]^{\sqrt2 +c}&&3\ar[ll]^{1}\\
}}
&&
&
&
\\
&&\xymatrix@R=0.5cm@C=0.8cm
{
&\ar@{-}[d]^{\mu_3}&\\
&&
\\
}&&&&\\
{\color{blue}\xymatrix@R=0.5cm@C=0.5cm
{
&1\ar[ld]_{c}&\\
2&&3\ar[ll]^{1}\ar[lu]_{\sqrt{2}}\\
}}
&
\xymatrix@R=0.5cm@C=0.6cm
{
&\\
\ar@{-}[r]^{\mu_2}&\\
}
&
\xymatrix@R=0.5cm@C=0.5cm
{
&1&\\
2\ar[ur]^c\ar[rr]_1&&3\ar[lu]_{\sqrt{2}}\\
}
&
\xymatrix@R=0.5cm@C=0.6cm
{
&\\
\ar@{-}[r]^{\mu_1}&\\
}
&
\xymatrix@R=0.5cm@C=0.5cm
{
&1\ar[ld]_c\ar[dr]^{\sqrt{2}}&\\
2\ar[rr]_1&&3\\
}&
\xymatrix@R=0.5cm@C=0.6cm
{
&\\
\ar@{-}[r]^{\mu_3}&\\
}&
{\color{blue}\xymatrix@R=0.5cm@C=0.5cm
{
&1\ar[ld]_c&\\
2&&3\ar[ul]_{\sqrt{2}}\ar[ll]^1\\
}}
\\
&&&&
\xymatrix@R=0.5cm@C=0.8cm
{
&\ar@{-}[d]^{\mu_2}&\\
&&\\
}
&&
\\
t_0:{\color{blue}\xymatrix@R=0.5cm@C=0.5cm
{
&1\ar[rd]^{c}&\\
2\ar[ur]^{\sqrt{2}}\ar[rr]_{1}&&3\\
}}
&
&
\xymatrix@R=0.5cm@C=0.5cm
{
& &\\
&\times&\times\\
}
&
\xymatrix@R=0.5cm@C=0.6cm
{
&\\
\ar@{-}[r]^{\mu_1}&\\
}
&
\xymatrix@R=0.5cm@C=0.5cm
{
&1\ar[dr]^{c+\sqrt{2}}&\\
2\ar[ur]^c&&3\ar[ll]^1\\
}
&&
{\color{blue}\xymatrix@R=0.5cm@C=0.5cm
{
&1\ar[dr]^{c}&\\
2\ar[ur]^{\sqrt{2}}\ar[rr]_1&&3\\
}}
\\
\xymatrix@R=0.5cm@C=0.8cm
{
&\ar@{-}[d]^{\mu_1}&\\
&&\\
}
&&
\xymatrix@R=0.5cm@C=0.8cm
{
&\ar@{-}[d]^{\mu_3}&\\
&&\\
}
&&
\xymatrix@R=0.5cm@C=0.8cm
{
&\ar@{-}[d]^{\mu_3}&\\
&&\\
}
&&
\xymatrix@R=0.5cm@C=0.8cm
{
&\ar@{-}[d]^{\mu_3}&\\
&&\\
}
\\
\xymatrix@R=0.5cm@C=0.5cm
{
&1\ar[dl]_{\sqrt{2}}&\\
2\ar[rr]_{\sqrt{2}c+1}&&3\ar[ul]_c\\
}
&
\xymatrix@R=0.5cm@C=0.6cm
{
&\\
\ar@{-}[r]^{\mu_2}&\\
}
&
\xymatrix@R=0.5cm@C=0.5cm
{
&1\ar[dr]^{\sqrt{2}+c}&\\
2\ar[ur]^{\sqrt{2}}&&3\ar[ll]^{\sqrt{2}c+1}\\
}
&
\xymatrix@R=0.5cm@C=0.6cm
{
&\\
\ar@{-}[r]^{\mu_1}&\\
}
&
{\color{red}\xymatrix@R=0.5cm@C=0.5cm
{
&1\ar[dl]_{\sqrt{2}}&\\
2\ar[rr]_{1}&&3\ar[ul]_{c+\sqrt{2}}\\
}}
&
\xymatrix@R=0.5cm@C=0.6cm
{
&\\
\ar@{-}[r]^{\mu_2}&\\
}
&
\xymatrix@R=0.5cm@C=0.5cm
{
&1&\\
2\ar[ru]^{\sqrt{2}}&&3\ar[ll]^{1}\ar[ul]_{c}\\
}
\\
\xymatrix@R=0.5cm@C=0.8cm
{
&\ar@{-}[d]^{\mu_3}&\\
&&\\
}
&&&&&&
\xymatrix@R=0.5cm@C=0.8cm
{
&\ar@{-}[d]^{\mu_1}&\\
&&\\
}
\\
\xymatrix@R=0.5cm@C=0.5cm
{
&\times &\\
&&\\
}
&&
{\color{orange}\xymatrix@R=0.5cm@C=0.5cm
{
&1\ar[dr]^{\sqrt{2}}&\\
2\ar[ur]^{c+\sqrt{2}}&&3\ar[ll]^{1}\\
}}
&
\xymatrix@R=0.5cm@C=0.6cm
{
&\\
\ar@{-}[r]^{\mu_2}&\\
}
&
\xymatrix@R=0.5cm@C=0.5cm
{
&1\ar[dl]_{c+\sqrt{2}}&\\
2\ar[rr]_{1}&&3\ar[ul]_{c}\\
}
&
\xymatrix@R=0.5cm@C=0.6cm
{
&\\
\ar@{-}[r]^{\mu_3}&\\
}
&
\xymatrix@R=0.5cm@C=0.5cm
{
&1\ar[dl]_{\sqrt{2}}\ar[dr]^c&\\
2&&3\ar[ll]^{1}\\
}\\
&&
\xymatrix@R=0.5cm@C=0.8cm
{
&\ar@{-}[d]^{\mu_3}&\\
&&\\
}
&&
\xymatrix@R=0.5cm@C=0.8cm
{
&\ar@{-}[d]^{\mu_1}&\\
&&\\
}
&&
\xymatrix@R=0.5cm@C=0.8cm
{
&\ar@{-}[d]^{\mu_2}&\\
&&\\
}\\
{\color{blue}\xymatrix@R=0.5cm@C=0.5cm
{
&1\ar[dl]_c&\\
2&&3\ar[ul]_{\sqrt{2}}\ar[ll]^1\\
}}
&
\xymatrix@R=0.5cm@C=0.6cm
{
&\\
\ar@{-}[r]^{\mu_2}&\\
}
&
\xymatrix@R=0.5cm@C=0.5cm
{
&1&\\
2\ar[ur]^{c}\ar[rr]_1&&3\ar[ul]_{\sqrt{2}}\\
}
&&
\xymatrix@R=0.5cm@C=0.5cm
{
&\times &\\
&&\\
}&&
{\color{blue}\xymatrix@R=0.5cm@C=0.5cm
{
&1\ar[dr]^{c}&\\
2\ar[ur]^{\sqrt{2}}\ar[rr]_{1}&&3\\
}}
\end{align*}
}
\caption{Mutation subgraph of $\Gamma_0$  of $(1,\sqrt{2},c\geq 2)$ in direction $1$.}
   \label{fig:1sqrt2c-1}
\end{figure}

According to Proposition \ref{prop:generator-set-assumption} and Remark \ref{rem:sharpen-K-2}, $\operatorname{Aut}(\mathcal{A})$ is generated by
\begin{align*}
   & g_{312}^{\mathbf{id},+},
    g_{213}^{\mathbf{id},+}
   & \in H_{1,0}^{t_\diamond}(t_0),\\
   & g_{212121}^{\mathbf{id},+}, 
    g_{32121}^{\mathbf{id},+},
    g_{2123121}^{(2 3),+},
    g_{323121}^{(2 3),+},
    g_{131232}^{(2 3),+},
    g_{313232}^{(2 3),+},
    g_{23232}^{(2 3),+},
    g_{121212}^{\mathbf{id},+}
    &\in H_{1,1}^{t_\diamond}(t_0),\\
   & g_{2123212}^{(2 3),+},
    g_{323212}^{(2 3),+},
    g_{1312313}^{(2 3),+},
    g_{3132313}^{(2 3),+},
    g_{232313}^{(2 3),+},
    g_{12123}^{\mathbf{id},+},
    g_{212323}^{(2 3),+},
    g_{32323}^{(2 3),+}
   & \in H_{1,1}^{t_\diamond}(t_0),\\
   & g_{232312121}^{(2 3),+},
    g_{2323123121}^{\mathbf{id},+}
    &\in K_{1}^{t_\diamond}(t_0).
\end{align*}
By Lemma \ref{lem:relation-rank-2-finite}, we have 
\[
    g_{212121}^{\mathbf{id},+}=g_{121212}^{\mathbf{id},+}=g_{23232}^{(2 3),+}=g_{32323}^{(2 3),+}=\mathbf{id}.
\]
On the other hand, 
\begin{equation*}
    \begin{aligned}
    &g_{212121}^{\mathbf{id},+}\circ g_{312}^{\mathbf{id},+}=g_{32121}^{\mathbf{id},+}, \\
     &g_{323121}^{(2 3),+}\circ g_{313232}^{(2 3),+}=g_{212121}^{\mathbf{id},+},\\
     &g_{213}^{\mathbf{id},+}\circ g_{131232}^{(2 3),+}=g_{1312313}^{(2 3),+},\\ &g_{213}^{\mathbf{id},+}\circ g_{313232}^{(2 3),+}=g_{3132313}^{(2 3),+},\\
&g_{213}^{\mathbf{id},+}\circ g_{23232}^{(2 3),+}=g_{232313}^{(2 3),+},
    \end{aligned}
    \qquad
    \begin{aligned}
        &g_{323121}^{(2 3),+}\circ g_{312}^{\mathbf{id},+}=g_{2123121}^{(2 3),+},\\
         &g_{323121}^{(2 3),+}\circ g_{23232}^{(2 3),+}=g_{32121}^{\mathbf{id},+},\\
         &g_{213}^{\mathbf{id},+}\circ g_{121212}^{\mathbf{id},+}=g_{12123}^{\mathbf{id},+},\\
&g_{213}^{\mathbf{id},+}\circ g_{2123212}^{(2 3),+}=g_{212323}^{(2 3),+},\\
&g_{213}^{\mathbf{id},+}\circ g_{323212}^{(2 3),+}=g_{32323}^{(2 3),+},\\
    \end{aligned}
    \qquad
    \begin{aligned}
        &g_{323121}^{(2 3),+}\circ g_{131232}^{(2 3),+}=\mathbf{id},\\
        &g_{312}^{\mathbf{id},+}\circ g_{213}^{\mathbf{id},+}=\mathbf{id},\\
       & g_{2123212}^{(2 3),+}\circ g_{1312313}^{(2 3),+}=g_{121212}^{\mathbf{id},+},\\
        &g_{323121}^{(2 3),+}\circ g_{323212}^{(2 3),+}=g_{2323123121}^{\mathbf{id},+},\\ &g_{232312121}^{(2 3),+}\circ g_{32323}^{(2 3),+}=g_{212121}^{\mathbf{id},+}.
    \end{aligned}
\end{equation*}
It follows  that $\operatorname{Aut}(\mathcal{A}) $ is generated by $g_{312}^{\mathbf{id},+}$. Hence, we conclude that $\operatorname{Aut}(\mathcal{A})\cong \mathbb{Z}$.

\subsubsection{$a=1,b=\sqrt{3},c=\sqrt{3}$.}
In this case, we have $d_2=d_3$ and $t_\diamond:=\mu_{12121}(t_0)$ satisfies the condition $(\diamond)$. The mutation subgraphs of $\Gamma_0$  are given in Figure \ref{fig:1sqrt3sqrt3-1} and \ref{fig:1sqrt3sqrt3-23}. Here,  for a path $\xymatrix{t\ar@{-}[r]^{i_1}&\cdots\ar@{-}[r]&\cdots\cdot\ar@{-}[r]^{i_k}&t'}$ in $\mathbb{T}_3$ such that for any non-endpoint vertex $s$, $B_s\not \sim \pm B_{t_0}$, $B_s\not\sim \pm B_{t_\diamond}$, and by mutating the diagram $\Gamma_s$ at $s$ in the third direction $k$, the resulting diagram has property $(M)$ at vertex $k$, we will abbreviate it by a wavy arrow labeled by $\mu_{i_k\cdots i_2i_1}$ in the mutation graph of $\Gamma_0$.



\begin{figure}[ht]
    \centering
    {\tiny\bf
\begin{align*}
t_0:{\color{blue}\xymatrix@R=0.5cm@C=0.5cm
{
&1\ar[rd]^{\sqrt 3}&\\
2\ar[ru]^{\sqrt 3}\ar[rr]_1&&3\\
}}
&&&&&&
{\color{blue}\xymatrix@R=0.5cm@C=0.5cm
{
&1\ar[rd]^{\sqrt 3}&\\
2\ar[ru]^{\sqrt 3}\ar[rr]_1&&3\\
}}
\\
\xymatrix@R=0.5cm@C=0.6cm
{
&\ar@{-}[d]^{\mu_1}&\\
&&\\
}
&&&&&&
\xymatrix@R=0.5cm@C=0.6cm
{
&\ar@{-}[d]^{ \mu_{2}}&\\
&&\\
}
\\
\xymatrix@R=0.5cm@C=0.5cm
{
&1\ar[ld]_{\sqrt 3}&\\
2\ar[rr]_{4}&&3\ar[lu]_{\sqrt 3}\\
}
&
\xymatrix@R=0.5cm@C=0.6cm
{
&\\
\ar@{~>}[r]^{\mu_{1212}}&\\
}
&
{\color{red}\xymatrix@R=0.5cm@C=0.5cm
{
&1\ar[ld]_{\sqrt 3}&\\
2\ar[rr]_{1}&&3\ar[lu]_{2\sqrt 3}\\
}}
&
\xymatrix@R=0.5cm@C=0.6cm
{
&\\
\ar@{-}[r]^{\mu_2}&\\
}
&
\xymatrix@R=0.5cm@C=0.5cm
{
&1&\\
2\ar[ru]^{\sqrt 3}&&3\ar[ll]^1\ar[lu]_{\sqrt 3}\\
}
&
\xymatrix@R=0.5cm@C=0.6cm
{
&\\
\ar@{-}[r]^{\mu_1}&\\
}
&
\xymatrix@R=0.5cm@C=0.5cm
{
&1\ar[ld]_{\sqrt 3}\ar[rd]^{\sqrt 3}&\\
2&&3\ar[ll]^{1}\\
}
\\
\xymatrix@R=0.5cm@C=0.6cm
{
&\ar@{~>}[d]_{ \mu_{1313}}&\\
&&\\
}
&
&
\xymatrix@R=0.5cm@C=0.6cm
{
&\ar@{-}[d]^{ \mu_{3}}&\\
&&\\
}
&
&
\xymatrix@R=0.5cm@C=0.6cm
{
&\ar@{-}[d]^{ \mu_{3}}&\\
&&\\
}
&
&
\xymatrix@R=0.5cm@C=0.6cm
{
&\ar@{-}[d]^{ \mu_{3}}&\\
&&\\
}
\\
{\color{orange}\xymatrix@R=0.5cm@C=0.5cm
{
&1\ar[ld]_{2\sqrt 3}&\\
2\ar[rr]_{1}&&3\ar[lu]_{\sqrt 3}\\
}}
&
\xymatrix@R=0.5cm@C=0.6cm
{
&\\
&\\
}
&
{\color{orange}\xymatrix@R=0.5cm@C=0.5cm
{
&1\ar[rd]^{2\sqrt 3}&\\
2\ar[ru]^{\sqrt 3}&&3\ar[ll]^1\\
}}
&
\xymatrix@R=0.5cm@C=0.6cm
{
&\\
&\\
}
&
{\color{blue}\xymatrix@R=0.5cm@C=0.5cm
{
&1\ar[rd]^{\sqrt 3}&\\
2\ar[ru]^{\sqrt 3}\ar[rr]_1&&3 \\
}}
&
\xymatrix@R=0.5cm@C=0.6cm
{
&\\
&\\
}
&
{\color{orange}\xymatrix@R=0.5cm@C=0.5cm
{
&1\ar[ld]_{2\sqrt 3}&\\
2\ar[rr]_{1}&&3\ar[lu]_{\sqrt 3}\\
}}
&
\xymatrix@R=0.5cm@C=0.6cm
{
&\\
&\\
}
\end{align*}
}
\caption{Mutation subgraph of $\Gamma_0$  of $(1,\sqrt{3},\sqrt{3})$ in direction $1$.}
    \label{fig:1sqrt3sqrt3-1}
\end{figure}
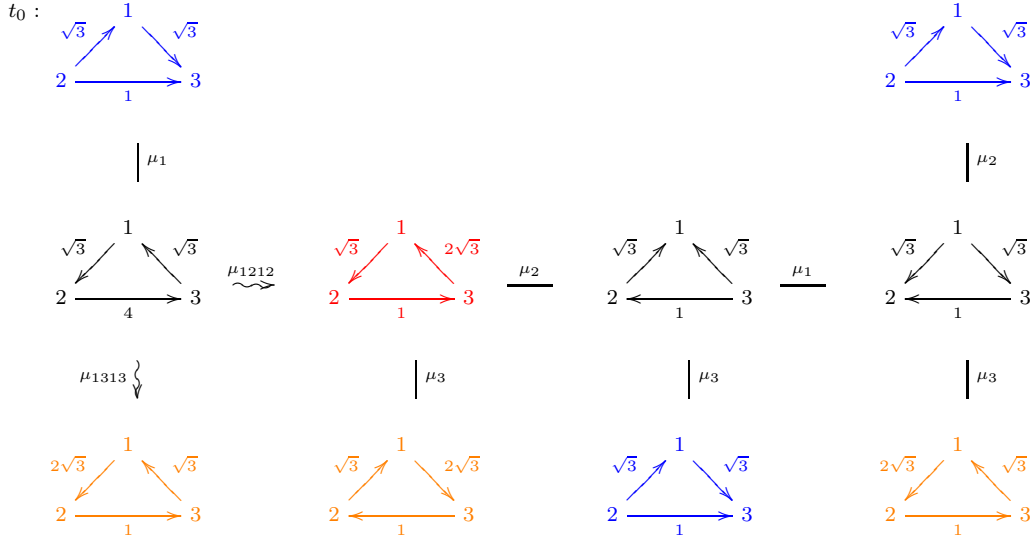

\begin{figure}[t]
  \centering
  {\tiny\bf
\begin{align*}
\xymatrix@R=0.5cm@C=0.5cm
{
&&\\
&&\\
}
&
\xymatrix@R=0.5cm@C=0.6cm
{
&\\
&\\
}
&
{\color{orange}\xymatrix@R=0.5cm@C=0.5cm
{
&1\ar[ld]_{2\sqrt 3}&\\
2\ar[rr]_{1}&&3\ar[lu]_{\sqrt 3}\\
}}
&
\xymatrix@R=0.5cm@C=0.6cm
{
&\\
&\\
}
&
{\color{orange}\xymatrix@R=0.5cm@C=0.5cm
{
&1\ar[ld]_{\sqrt 3}&\\
2\ar[rr]_{1}&&3\ar[lu]_{2\sqrt 3}\\
}}
&
\xymatrix@R=0.5cm@C=0.6cm
{
&\\
&\\
}
&
\xymatrix@R=0.5cm@C=0.5cm
{
&&\\
&&\\
}
&
\xymatrix@R=0.5cm@C=0.6cm
{
&\\
&\\
}
&
\xymatrix@R=0.5cm@C=0.5cm
{
&&\\
&&\\
}
\\
&&
\xymatrix@R=0.5cm@C=0.8cm
{
&\ar@{-}[d]^{\mu_3}&\\
&&\\
}
&&
\xymatrix@R=0.5cm@C=0.8cm
{
&\ar@{-}[d]^{\mu_2}&\\
&&\\
}
&&
\\
t_0:{\color{blue}\xymatrix@R=0.5cm@C=0.5cm
{
&1\ar[rd]^{\sqrt 3}&\\
2\ar[ru]^{\sqrt 3}\ar[rr]_1&&3\\
}}
&
\xymatrix@R=0.5cm@C=0.6cm
{
&\\
\ar@{-}[r]^{\mu_2}&\\
}
&
\xymatrix@R=0.5cm@C=0.5cm
{
&1\ar[ld]_{\sqrt 3}\ar[rd]^{\sqrt 3}&\\
2&&3\ar[ll]^{1}\\
}
&
\xymatrix@R=0.5cm@C=0.6cm
{
&\\
\ar@{-}[r]^{\mu_1}&\\
}
&
\xymatrix@R=0.5cm@C=0.5cm
{
&1&\\
2\ar[ru]^{\sqrt 3}&&3\ar[ll]^1\ar[lu]_{\sqrt 3}\\
}
&
\xymatrix@R=0.5cm@C=0.6cm
{
&\\
\ar@{-}[r]^{\mu_3}&\\
}
&
{\color{blue}\xymatrix@R=0.5cm@C=0.5cm
{
&1\ar[rd]^{\sqrt 3}&\\
2\ar[ru]^{\sqrt 3}\ar[rr]_1&&3\\
}}
\\
\xymatrix@R=0.5cm@C=0.8cm
{
&\ar@{-}[d]^{\mu_3}&\\
&&\\
}
&&&&&&
\\
\xymatrix@R=0.5cm@C=0.5cm
{
&1&\\
2\ar[ru]^{\sqrt 3}&&3\ar[ll]^{1}\ar[lu]_{\sqrt 3}\\
}
&
\xymatrix@R=0.5cm@C=0.6cm
{
&\\
\ar@{-}[r]^{\mu_1}&\\
}
&
\xymatrix@R=0.5cm@C=0.5cm
{
&1\ar[ld]_{\sqrt 3}\ar[rd]^{\sqrt 3}&\\
2&&3\ar[ll]^{1}\\
}
&
\xymatrix@R=0.5cm@C=0.6cm
{
&\\
\ar@{-}[r]^{\mu_2}&\\
}
&
{\color{blue}\xymatrix@R=0.5cm@C=0.5cm
{
&1\ar[rd]^{\sqrt 3}&\\
2\ar[ru]^{\sqrt 3}\ar[rr]_1&&3\\
}}
&
\xymatrix@R=0.5cm@C=0.6cm
{
&\\
&\\
}
&
\xymatrix@R=0.5cm@C=0.5cm
{
&&\\
&&\\
}
&&
\\
\xymatrix@R=0.5cm@C=0.8cm
{
&\ar@{-}[d]^{\mu_2}&\\
&&\\
}
&
\xymatrix@R=0.5cm@C=0.6cm
{
&\\
&\\
}
&
\xymatrix@R=0.5cm@C=0.8cm
{
&\ar@{-}[d]^{\mu_3}&\\
&&\\
}
&&&&
\\
{\color{orange}\xymatrix@R=0.5cm@C=0.5cm
{
&1\ar[ld]_{\sqrt 3}&\\
2\ar[rr]_{1}&&1\ar[lu]_{2\sqrt 3}\\
}}
&
\xymatrix@R=0.5cm@C=0.6cm
{
&\\
&\\
}
&
{\color{orange}\xymatrix@R=0.5cm@C=0.5cm
{
&1\ar[ld]_{2\sqrt 3}&\\
2\ar[rr]_{1}&&1\ar[lu]_{\sqrt 3}\\
}}
&&&&\\
\end{align*}
}
    \caption{Mutation subgraph of $\Gamma_0$ of $(1,\sqrt{3},\sqrt{3})$ in directions $2$ and $3$.}
    \label{fig:1sqrt3sqrt3-23}
\end{figure}

By Proposition \ref{prop:generator-set-assumption} and Remark \ref{rem:sharpen-K-2}, $\operatorname{Aut}(\mathcal{A})$ is generated by the following generators:
\begin{align*}
   & \psi_{(23)}^-&\in G_0(t_0),\\
   & g_{312}^{\mathbf{id},+}, 
    g_{213}^{\mathbf{id},+}
    &\in H_{1,0}^{t_\diamond}(t_0),\\
   & g_{21212121}^{\mathbf{id},+},
    g_{3212121}^{\mathbf{id},+},
    g_{31313131}^{\mathbf{id},+},
    g_{2313131}^{\mathbf{id},+},
    g_{1313132}^{\mathbf{id},+},
    g_{12121212}^{\mathbf{id},+},
    g_{1212123}^{\mathbf{id},+},
    g_{13131313}^{\mathbf{id},+}
    &\in H_{1,1}^{t_\diamond}(t_0),\\
   & g_{32312121}^{(2 3),+},
    g_{1313131212121}^{\mathbf{id},+}
    &\in K_{1}^{t_\diamond}(t_0).
\end{align*}
By Lemma \ref{lem:relation-rank-2-finite}, we deduce that 
\[
 g_{12121212}^{\mathbf{id},+}=g_{21212121}^{\mathbf{id},+}=\mathbf{id}, g_{13131313}^{\mathbf{id},+}=g_{31313131}^{\mathbf{id},+}=\mathbf{id},g_{1212123}^{\mathbf{id},+}\circ g_{32312121}^{(2 3),+}=\mathbf{id}.
\]
On the other hand, one can check that 
\[
    g_{21212121}^{\mathbf{id},+}=g_{3212121}^{\mathbf{id},+}\circ g_{213}^{\mathbf{id},+}, g_{31313131}^{\mathbf{id},+}=g_{2313131}^{\mathbf{id},+}\circ g_{312}^{\mathbf{id},+}, g_{13131313}^{\mathbf{id},+}=g_{213}^{\mathbf{id},+}\circ g_{1313132}^{\mathbf{id},+},
\]
\[
    g_{12121212}^{\mathbf{id},+}=g_{312}^{\mathbf{id},+}\circ g_{1212123}^{\mathbf{id},+}, g_{1313131212121}^{\mathbf{id},+}=g_{21212121}^{\mathbf{id},+}\circ g_{1313132}^{\mathbf{id},+},
\]
\[
    g_{312}^{\mathbf{id},+}\circ g_{213}^{\mathbf{id},+}=\mathbf{id}, g_{312}^{\mathbf{id},+}\circ \psi_{(23)}^-=\psi_{(23)}^-\circ (g_{312}^{\mathbf{id},+})^{-1}.
\]
It follows that $\operatorname{Aut}(\mathcal{A})$ is generated by $\psi_{(23)}^-$ and $g_{312}^{\mathbf{id},+}$. As before one can show that $\operatorname{Aut}(\mathcal{A})\cong D_\infty.$

\subsubsection{$a=1,b=\sqrt{3},c\geq 2$.}
It is clear that $d_2=d_3$ and $t_\diamond:=\mu_{12121}(t_0)$ satisfies the condition $(\diamond)$. The mutation subgraphs of $\Gamma_0$ are given in Figure \ref{fig:1sqrt3c-1} and \ref{fig:1sqrt3c-23}. 




\begin{figure}[ht]
{\tiny\bf
\begin{align*}
t_0:{\color{blue}\xymatrix@R=0.5cm@C=0.5cm
{
&1\ar[rd]^{c}&\\
2\ar[ru]^{\sqrt 3}\ar[rr]_1&&3\\
}}
&&
{\color{blue}\xymatrix@R=0.5cm@C=0.5cm
{
&1\ar[rd]^{c}&\\
2\ar[ru]^{\sqrt 3}\ar[rr]_1&&3\\
}}
&&&&
\\
\xymatrix@R=0.5cm@C=0.6cm
{
&\ar@{~>}[d]^{\mu_{12121}}&\\
&&\\
}
&&
\xymatrix@R=0.5cm@C=0.6cm
{
&\ar@{-}[d]^{\mu_3}&\\
&&\\
}
&&&&
\\
{\color{red}\xymatrix@R=0.5cm@C=0.5cm
{
&1\ar[ld]_{\sqrt 3}&\\
2\ar[rr]_{1}&&3\ar[lu]_{c+\sqrt 3}\\
}}
&
\xymatrix@R=0.5cm@C=0.6cm
{
&\\
\ar@{-}[r]^{\mu_2}&\\
}
&
\xymatrix@R=0.5cm@C=0.5cm
{
&1&\\
2\ar[ru]^{\sqrt 3}&&3\ar[lu]_{c}\ar[ll]^1\\
}
&
\xymatrix@R=0.5cm@C=0.6cm
{
&\\
\ar@{-}[r]^{\mu_1}&\\
}
&
\xymatrix@R=0.5cm@C=0.5cm
{
&1\ar[rd]^{c}\ar[ld]_{\sqrt 3}&\\
2&&3\ar[ll]^{1}\\
}
&
\xymatrix@R=0.5cm@C=0.6cm
{
&\\
\ar@{-}[r]^{\mu_2}&\\
}
&
{\color{blue}\xymatrix@R=0.5cm@C=0.5cm
{
&1\ar[rd]^{c}&\\
2\ar[ru]^{\sqrt 3}\ar[rr]_1&&3\\
}}
\\
\xymatrix@R=0.5cm@C=0.8cm
{
&\ar@{-}[d]^{\mu_3}&\\
&&\\
}
&&&&
\xymatrix@R=0.5cm@C=0.8cm
{
&\ar@{-}[d]^{\mu_3}&\\
&&\\
}
&
&
\\
\xymatrix@R=0.5cm@C=0.5cm
{
&1\ar[rd]^{c+\sqrt 3}&\\
2\ar[ru]^{c}&&3\ar[ll]^1\\
}
&
\xymatrix@R=0.5cm@C=0.6cm
{
&\\
\ar@{-}[r]^{\mu_1}&\\
}
&
\xymatrix@R=0.5cm@C=0.5cm
{
&&\\
\times&&\times \\
}
&
\xymatrix@R=0.5cm@C=0.6cm
{
&\\
\ar@{-}[r]^{\mu_1}&\\
}
&
\xymatrix@R=0.5cm@C=0.5cm
{
&1\ar[ld]_{c+\sqrt 3}&\\
2\ar[rr]_{1}&&3\ar[lu]_{c}\\
}
&
\xymatrix@R=0.5cm@C=0.6cm
{
&\\
\ar@{-}[r]^{\mu_2}&\\
}
&
{\color{orange}\xymatrix@R=0.5cm@C=0.5cm
{
&1\ar[rd]^{\sqrt 3}&\\
2\ar[ru]^{c+\sqrt 3}&&3\ar[ll]^1\\
}}
&\\
\xymatrix@R=0.5cm@C=0.8cm
{
&\ar@{-}[d]^{\mu_2}&\\
&&\\
}
&&&&&&
\\
\xymatrix@R=0.5cm@C=0.5cm
{
&1\ar[ld]_{c}\ar[dr]^{\sqrt{3}}&\\
2\ar[rr]_{1}&&3\\
}
&
\xymatrix@R=0.5cm@C=0.6cm
{
&\\
\ar@{-}[r]^{\mu_1}&\\
}
&
\xymatrix@R=0.5cm@C=0.5cm
{
&1&\\
2\ar[ru]^{c}\ar[rr]_1&&3\ar[ul]_{\sqrt 3}\\
}
&
\xymatrix@R=0.5cm@C=0.6cm
{
&\\
\ar@{-}[r]^{\mu_3}&\\
}
&
{\color{orange}\xymatrix@R=0.5cm@C=0.5cm
{
&1\ar[rd]^{\sqrt 3}&\\
2\ar[ru]^{c+\sqrt{3}}&&3\ar[ll]^1\\
}}
&&
\\
\xymatrix@R=0.5cm@C=0.8cm
{
&\ar@{-}[d]^{\mu_3}&\\
&&\\
}
&&
\xymatrix@R=0.5cm@C=0.8cm
{
&\ar@{-}[d]^{\mu_2}&\\
&&\\
}
&&&&
\\
{\color{blue}\xymatrix@R=0.5cm@C=0.5cm
{
&1\ar[ld]_{c}&\\
2&&3\ar[ll]^1\ar[ul]_{\sqrt{3}}\\
}}
&&
{\color{blue}\xymatrix@R=0.5cm@C=0.5cm
{
&1\ar[ld]_{c}&\\
2&&3\ar[ll]^1\ar[ul]_{\sqrt{3}}\\
}}
&&&&
\end{align*}
}
\caption{Mutation subgraph of $\Gamma_0$ of  $(1,\sqrt{3},c\geq 2)$ in direction $1$.}
\label{fig:1sqrt3c-1}
\end{figure}

\begin{figure}[t]
    \centering
    {\tiny\bf
    \begin{align*}
{\color{orange}\xymatrix@R=0.5cm@C=0.5cm
{
&1\ar[rd]^{\sqrt 3}&\\
2\ar[ru]^{c+\sqrt 3}&&3\ar[ll]^{1}\\
}}
&
\xymatrix@R=0.5cm@C=0.6cm
{
&\\
\ar@{-}[r]^{\mu_2}&\\
}
&
\xymatrix@R=0.5cm@C=0.5cm
{
&1\ar[ld]_{c+\sqrt 3}&\\
2\ar[rr]_{1}&&3\ar[lu]_{c}\\
}
&
\xymatrix@R=0.5cm@C=0.6cm
{
&\\
\ar@{-}[r]^{\mu_1}&\\
}
&
\xymatrix@R=0.5cm@C=0.5cm
{
&& \\
\times&& \\
}
&&
{\color{orange}\xymatrix@R=0.5cm@C=0.5cm
{
&1\ar[ld]_{\sqrt 3}&\\
2\ar[rr]_{1}&&3\ar[lu]_{c+\sqrt 3}\\
}}
&
\\
&&
\xymatrix@R=0.5cm@C=0.8cm
{
&\ar@{-}[d]_{\mu_3}& \\
&&\\
}
&&&
\xymatrix@R=0.5cm@C=0.6cm
{
&\ar@{-}[ld]_{\mu_2} \\
&\\
}&&
    \\
t_0:{\color{blue}\xymatrix@R=0.5cm@C=0.5cm
{
&1\ar[rd]^{c}&\\
2\ar[ru]^{\sqrt 3}\ar[rr]_1&&3\\
}}
&
\xymatrix@R=0.5cm@C=0.6cm
{
&\\
\ar@{-}[r]^{\mu_2}&\\
}
&
\xymatrix@R=0.5cm@C=0.5cm
{
&1\ar[rd]^{c}\ar[ld]_{\sqrt 3}&\\
2&&3\ar[ll]^1 \\
}
&
\xymatrix@R=0.5cm@C=0.6cm
{
&\\
\ar@{-}[r]^{\mu_1}&\\
}
&
\xymatrix@R=0.5cm@C=0.5cm
{
&1&\\
2\ar[ru]^{\sqrt 3}&&3\ar[lu]_{c}\ar[ll]^1\\
}
&
\xymatrix@R=0.5cm@C=0.6cm
{
&\\
\ar@{-}[r]^{\mu_3}&\\
}
&
{\color{blue}\xymatrix@R=0.5cm@C=0.5cm
{
&1\ar[rd]^{c}&\\
2\ar[ru]^{\sqrt 3}\ar[rr]_1&&3\\
}}
\\
\xymatrix@R=0.5cm@C=0.8cm
{
&\ar@{-}[d]_{\mu_3}& \\
&&\\
}
&&&&&&
\\
\xymatrix@R=0.5cm@C=0.5cm
{
&1&\\
2\ar[ru]^{\sqrt 3}&&3\ar[ll]^1\ar[ul]_c\\
}
&
\xymatrix@R=0.5cm@C=0.6cm
{
&\\
\ar@{-}[r]^{\mu_1}&\\
}
&
\xymatrix@R=0.5cm@C=0.5cm
{
&1\ar[ld]_{\sqrt{3}}\ar[dr]^c&\\
2&&3\ar[ll]^1\\
}
&
\xymatrix@R=0.5cm@C=0.6cm
{
&\\
\ar@{-}[r]^{\mu_3}&\\
}
&
\xymatrix@R=0.5cm@C=0.5cm
{
&1\ar[ld]_{c+\sqrt{3}}&\\
2\ar[rr]_1&&3\ar[lu]_c\\
}
&
\xymatrix@R=0.5cm@C=0.6cm
{
&\\
\ar@{-}[r]^{\mu_2}&\\
}
&
{\color{orange}\xymatrix@R=0.5cm@C=0.5cm
{
&1\ar[rd]^{\sqrt{3}}&\\
2\ar[ur]^{c+\sqrt{3}}&&3\ar[ll]^1\\
}}\\
\xymatrix@R=0.5cm@C=0.8cm
{
&\ar@{-}[d]_{\mu_2}& \\
&&\\
}
&&
\xymatrix@R=0.5cm@C=0.8cm
{
&\ar@{-}[d]_{\mu_2}& \\
&&\\
}
&&
\xymatrix@R=0.5cm@C=0.8cm
{
&\ar@{-}[d]_{\mu_1}& \\
&&\\
}
&&
\\
{\color{orange}\xymatrix@R=0.5cm@C=0.5cm
{
&1\ar[ld]_{\sqrt{3}}&\\
2\ar[rr]_{1}&&3\ar[lu]_{c+\sqrt{3}}\\
}}&&
{\color{blue}\xymatrix@R=0.5cm@C=0.5cm
{
&1\ar[rd]^{c}&\\
2\ar[ur]^{\sqrt{3}}\ar[rr]_1&&3\\
}}
&&
\xymatrix@R=0.5cm@C=0.5cm
{
&\times& \\
&& \\
}
&&
\\
\end{align*}
}
    \caption{Mutation subgraph of $\Gamma_0$  of $(1,\sqrt{3},c\ge 2)$ in directions $2$ and $3$.}
    \label{fig:1sqrt3c-23}
\end{figure}

By Proposition \ref{prop:generator-set-assumption} and Remark \ref{rem:sharpen-K-2}, we obtain the following generators of $\operatorname{Aut}(\mathcal{A})$:
\begin{align*}
    & g_{312}^{\mathbf{id},+}, 
    g_{213}^{\mathbf{id},+}
    &\in H_{1,0}^{t_\diamond}(t_0),\\
    & g_{21212121}^{\mathbf{id},+},
     g_{3212121}^{\mathbf{id},+},
     g_{212312121}^{(2 3),+}, 
     g_{32312121}^{(2 3),+},
     g_{13131232}^{(2 3),+}, 
     g_{313232}^{(2 3),+},
     g_{23232}^{(2 3),+},
     g_{12121212}^{\mathbf{id},+}
     &\in H_{1,1}^{t_\diamond}(t_0),\\
    & g_{2123212}^{(2 3),+},
     g_{323212}^{(2 3),+},
     g_{131312313}^{(2 3),+},
     g_{3132313}^{(2 3),+},
     g_{232313}^{(2 3),+},
     g_{1212123}^{\mathbf{id},+},
     g_{212323}^{(2 3),+},
     g_{32323}^{(2 3),+}
     &\in H_{1,1}^{t_\diamond}(t_0),\\
    & g_{23231212121}^{(2 3),+},
     g_{232312312121}^{\mathbf{id},+}
     &\in K_{1}^{t_\diamond}(t_0).
\end{align*}
As a consequence of Lemma \ref{lem:relation-rank-2-finite}, we have
\[
    g_{21212121}^{\mathbf{id},+}=g_{12121212}^{\mathbf{id},+}=\mathbf{id}, g_{23232}^{(2 3),+}=g_{32323}^{(2 3),+}=\mathbf{id}.
\]
Moreover, by direction calculation, one can show that the above generators can be generated by $g_{312}^{\mathbf{id},+}$. Therefore, $\operatorname{Aut}(A)$ is generated by $g_{312}^{\mathbf{id},+}$ and the order of $g_{312}^{\mathbf{id},+}$ is infinite by Lemma \ref{lem:order-of-tau}, which implies that $\operatorname{Aut}(\mathcal{A})\cong \mathbb{Z}$.

\subsubsection{$a=1,b\geq 2,c\geq 2$}\label{ss:1bgeq2cgeq2}
Clearly, we have $d_2=d_3$ and $t_\diamond=\mu_2(t_0)$ satisfies the condition $(\diamond)$. Mutation subgraphs of $\Gamma_0$ are given in Figure \ref{fig:1bc-12} and \ref{fig:1bc-3}.



\begin{figure}
    \centering
    {\tiny\bf
\begin{align*}
\xymatrix@R=0.5cm@C=0.5cm
{
&&\\
&\times&\\
}&&&&
{\color{orange}\xymatrix@R=0.5cm@C=0.5cm
{
&1\ar[rd]^{b}\ar[ld]_{c}&\\
2\ar[rr]_{1}&&3\\
}}
&&\\
\xymatrix@R=0.5cm@C=0.8cm
{
&\ar@{-}[d]_{\mu_1}& \\
&&\\
}
&&&&
\xymatrix@R=0.5cm@C=0.8cm
{
&&\\
&\ar@{~>}[u]_{\mu_{232}}&\\
}
&&\\
t_0:{\color{blue}\xymatrix@R=0.5cm@C=0.5cm
{
&1\ar[rd]^{c}&\\
2\ar[ru]^{b}\ar[rr]_1&&3\\
}}
&
\xymatrix@R=0.5cm@C=0.6cm
{
&\\
\ar@{-}[r]^{\mu_2}&\\
}
&
{\color{red}\xymatrix@R=0.5cm@C=0.5cm
{
&1\ar[ld]_{b}\ar[rd]^c&\\
2&&3\ar[ll]^1 \\
}}
&
\xymatrix@R=0.5cm@C=0.6cm
{
&\\
\ar@{-}[r]^{\mu_1}&\\
}
&
\xymatrix@R=0.5cm@C=0.5cm
{
&1&\\
2\ar[ru]^{b}&&3\ar[ll]^{1}\ar[lu]_c\\
}
&
\xymatrix@R=0.5cm@C=0.6cm
{
&\\
\ar@{-}[r]^{\mu_3}&\\
}
&
{\color{blue}\xymatrix@R=0.5cm@C=0.5cm
{
&1\ar[rd]^{c}&\\
2\ar[rr]_1\ar[ru]^{b}&&3\\
}}
&
\xymatrix@R=0.5cm@C=0.6cm
{
&\\
&\\
}
&
\xymatrix@R=0.5cm@C=0.5cm
{
&&\\
&&\\
}
\\
\xymatrix@R=0.5cm@C=0.8cm
{
&&\\
&&\\
}
&
\xymatrix@R=0.5cm@C=0.6cm
{
&\\
&\\
}
&
\xymatrix@R=0.5cm@C=0.8cm
{
&\ar@{~>}[d]^{\mu_{323}}&\\
&&\\
}
&
\xymatrix@R=0.5cm@C=0.6cm
{
&\\
&\\
}
&
\xymatrix@R=0.5cm@C=0.8cm
{
&&\\
&&\\
}
&
&
\xymatrix@R=0.5cm@C=0.8cm
{
&&\\
&&\\
}
&
\xymatrix@R=0.5cm@C=0.6cm
{
&\\
&\\
}
&
\xymatrix@R=0.5cm@C=0.8cm
{
&&\\
&&\\
}
\\
{\color{orange}\xymatrix@R=0.5cm@C=0.5cm
{
&1\ar[rd]^{b}\ar[ld]_{c}&\\
2\ar[rr]_{1}&&3\\
}}
&
\xymatrix@R=0.5cm@C=0.6cm
{
&\\
\ar@{-}[r]^{\mu_1}&\\
}
&
\xymatrix@R=0.5cm@C=0.5cm
{
&1&\\
2\ar[rr]_{1}\ar[ru]^c&&3\ar[lu]_{b}\\
}
&
\xymatrix@R=0.5cm@C=0.6cm
{
&\\
\ar@{-}[r]^{\mu_2}&\\
}
&
{\color{blue}\xymatrix@R=0.5cm@C=0.5cm
{
&1\ar[ld]_{c}&\\
2&&3\ar[ll]^1\ar[lu]_b \\
}}
&
\xymatrix@R=0.5cm@C=0.6cm
{
&\\
&\\
}
&
&
\xymatrix@R=0.5cm@C=0.6cm
{
&\\
&\\
}
&
\xymatrix@R=0.5cm@C=0.5cm
{
&&\\
&&\\
}
\end{align*}
}
\caption{Mutation subgraph of $\Gamma_0$ of $(1,b\ge 2,c\ge 2)$ in directions $1$ and $2$.}
    \label{fig:1bc-12}
\end{figure}
\begin{figure}
    \centering
{\tiny\bf
\begin{align*}
t_0:{\color{blue}\xymatrix@R=0.5cm@C=0.5cm
{
&1\ar[rd]^{c}&\\
2\ar[ru]^{b}\ar[rr]_1&&3\\
}}
&
\xymatrix@R=0.5cm@C=0.6cm
{
&\\
\ar@{-}[r]^{\mu_3}&\\
}
&
\xymatrix@R=0.5cm@C=0.5cm
{
&1&\\
2\ar[ru]^{b}&&3\ar[lu]_{c}\ar[ll]^1\\
}
&
\xymatrix@R=0.5cm@C=0.6cm
{
&\\
\ar@{-}[r]^{\mu_1}&\\
}
&
{\color{orange}\xymatrix@R=0.5cm@C=0.5cm
{
&1\ar[rd]^{c}\ar[ld]_{b}&\\
2&&3\ar[ll]^1 \\
}}
&
\xymatrix@R=0.5cm@C=0.6cm
{
&\\
\ar@{-}[r]^{\mu_2}&\\
}
&
{\color{blue}\xymatrix@R=0.5cm@C=0.5cm
{
&1\ar[rd]^{c}&\\
2\ar[ru]^{b}\ar[rr]_1&&3\\
}}
&
\xymatrix@R=0.5cm@C=0.6cm
{
&\\
&\\
}
&
\xymatrix@R=0.5cm@C=0.5cm
{
&&\\
&&\\
}
\\
\xymatrix@R=0.5cm@C=0.8cm
{
&& \\
&&\\
}
&
\xymatrix@R=0.5cm@C=0.6cm
{
& \ar@{~>}[ld]_{\mu_{232}}\\
&\\
}
&
\xymatrix@R=0.5cm@C=0.8cm
{
&& \\
&&\\
}
&
\xymatrix@R=0.5cm@C=0.6cm
{
&\\
&\\
}
&
\xymatrix@R=0.5cm@C=0.8cm
{
&\ar@{~>}[d]^{\mu_{323}}&\\
&&\\
}
&
\xymatrix@R=0.5cm@C=0.6cm
{
&\\
&\\
}
&
\xymatrix@R=0.5cm@C=0.8cm
{
&& \\
&&\\
}
&
\xymatrix@R=0.5cm@C=0.6cm
{
&\\
&\\
}
&
\xymatrix@R=0.5cm@C=0.8cm
{
&& \\
&&\\
}
\\
{\color{orange}\xymatrix@R=0.5cm@C=0.5cm
{
&1\ar[rd]^{b}\ar[ld]_{c}&\\
2\ar[rr]_{1}&&3\\
}}
&
\xymatrix@R=0.5cm@C=0.6cm
{
&\\
&\\
}
&
{\color{orange}\xymatrix@R=0.5cm@C=0.5cm
{
&1\ar[rd]^{b}\ar[ld]_{c}&\\
2\ar[rr]_{1}&&3\\
}}
&
\xymatrix@R=0.5cm@C=0.6cm
{
&\\
\ar@{-}[r]^{\mu_1}&\\
}
&
\xymatrix@R=0.5cm@C=0.5cm
{
&1&\\
2\ar[rr]_{1}\ar[ru]^c&&3\ar[lu]_{b}\\
}
&
\xymatrix@R=0.5cm@C=0.6cm
{
&\\
\ar@{-}[r]^{\mu_2}&\\
}
&
{\color{blue}\xymatrix@R=0.5cm@C=0.5cm
{
&1\ar[ld]_{c}&\\
2&&3\ar[ll]^1\ar[lu]_b \\
}}
&
\xymatrix@R=0.5cm@C=0.6cm
{
&\\
&\\
}
&
\xymatrix@R=0.5cm@C=0.5cm
{
&&\\
&& \\
}
&
\end{align*}
}
\caption{Mutation subgraph of $\Gamma_0$ of $(1,b\ge 2,c\ge 2)$ in direction $3$.}
    \label{fig:1bc-3}
\end{figure}
Assume that $b\neq c$, then $G_0(t_0)=\langle\mathbf{id}\rangle$. It follows from Proposition \ref{prop:generator-set-assumption} and Remark \ref{rem:sharpen-K-2} that $\operatorname{Aut}(\mathcal{A})$ is generated by
\begin{align*}
    &g_{312}^{\mathbf{id},+}, 
    g_{213}^{\mathbf{id},+},
    g_{23232}^{(2 3), +}, 
    g_{232313}^{(2 3), +},
    g_{32323}^{(2 3), +}, g_{212323}^{(23),+} 
    &\in H_{1,1}^{t_\diamond}(t_0),\\
    &g_{323212}^{(2 3), +},
    g_{313232}^{(2 3), +}
    &\in K_{1}^{t_\diamond}(t_0).
\end{align*}    

Assume that $b=c$, then $|G_0(t_0)|=2$ and $\psi_{(23)}^-\in G_0(t_0)$. Moreover, $\mu_{12}(B_{t_0})\sim \pm B_{t_\diamond}$, $\mu_{32}^2(B_{t_0})\sim \pm B_{t_\diamond}$, $\mu_3(B_{t_0})\sim \pm B_{t_\diamond}$, $\mu_{13}(B_{t_0})\sim \pm B_{t_\diamond}$, and $\mu_{23}^2(B_{t_0})\sim \pm B_{t_\diamond}$. Consequently, $\operatorname{Aut}(\mathcal{A})$ is generated by 
\begin{align*}
       \psi_{(23)}^-\in G_0(t_0) \quad \text{and}\quad 
        g_{312}^{\mathbf{id},+},
        g_{23232}^{(2 3), +}
        \in K_{1}^{t_\diamond}(t_0).
\end{align*}
By Lemma \ref{lem:relation-rank-2-finite} and direct calculation, one can show that $\operatorname{Aut}(\mathcal{A})$ is generated by $g_{312}^{\mathbf{id},+}$ for $b\neq c$, and by  $g_{312}^{\mathbf{id},+}$ and $\psi_{(23)}^-$ for $b=c$. Similar to the previous cases, we have 
\[
    \operatorname{Aut}(\mathcal{A})\cong \begin{cases}
    \mathbb{Z} & \text{if $b\neq c$;}\\
    D_\infty & \text{if $b=c$.}
    \end{cases}
\]

\subsection{Case 3: $a=\sqrt{2}$.}
Note that  $abc$ is an integer and $a\leq b\leq c$. If $b=\sqrt{2}$, then $c\geq 2$. If $b=\sqrt{3}$, then $c\geq \sqrt{6}>2$.
\subsubsection{$a=\sqrt{2},b=\sqrt{2},c\geq 2$.}\label{ss:a=sqrt2-b=sqrt2-cge2}


 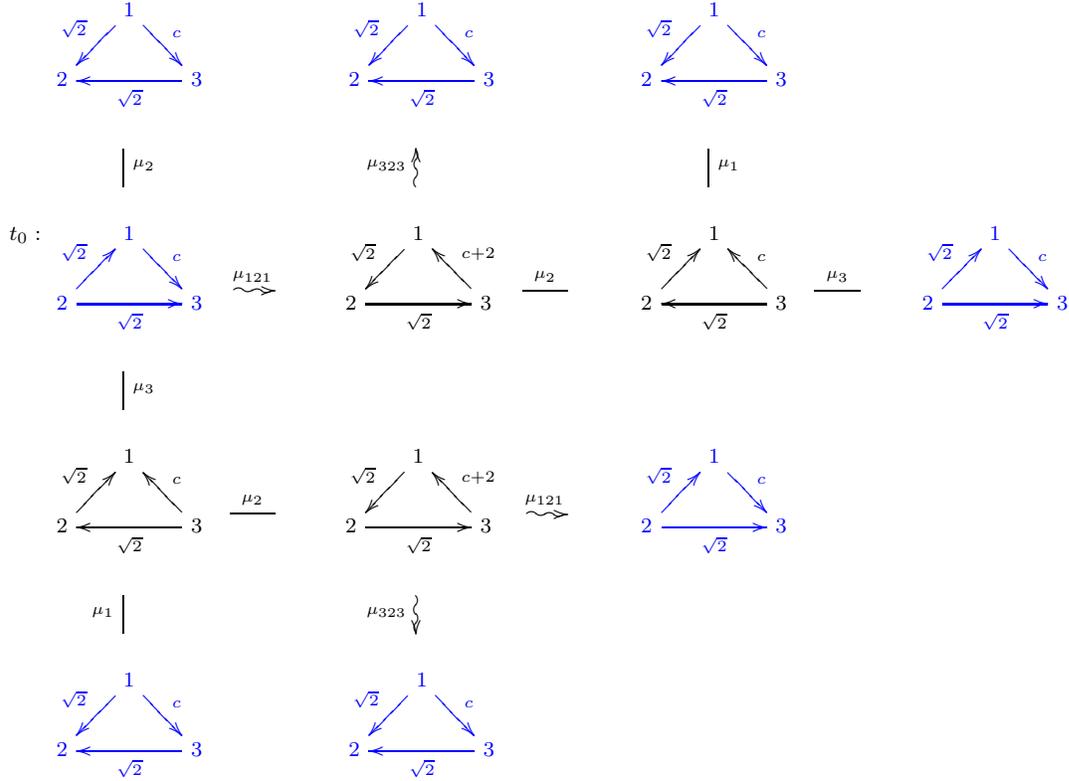
\begin{figure}
     \centering
     {\tiny\bf
     \begin{align*}     {\color{blue}\xymatrix@R=0.5cm@C=0.5cm
 {
 &1\ar[rd]^{c}\ar[ld]_{\sqrt 2}&\\
 2&&3\ar[ll]^{\sqrt 2} 
 }}
 &&
 {\color{blue}\xymatrix@R=0.5cm@C=0.5cm
 {
 &1\ar[rd]^{c}\ar[ld]_{\sqrt 2}&\\
 2&&3\ar[ll]^{\sqrt 2} \\
 }}
 &
 &
 {\color{blue}\xymatrix@R=0.5cm@C=0.5cm
 {
 &1\ar[rd]^{c}\ar[ld]_{\sqrt 2}&\\
 2&&3\ar[ll]^{\sqrt 2} \\
 }}&&
     \\
 \xymatrix@R=0.5cm@C=0.8cm
 {
 &\ar@{-}[d]^{\mu_2}&\\
 &&\\
 }
 &&
 \xymatrix@R=0.5cm@C=0.8cm
 {
 &&\\
 &\ar@{~>}[u]^{ \mu_{323}}&\\
 }
 &
 &
 \xymatrix@R=0.5cm@C=0.8cm
 {
 &\ar@{-}[d]^{\mu_1}&\\
 &&\\
 }
 &&
 \\
     t_0:
 {\color{blue}\xymatrix@R=0.5cm@C=0.5cm
 {
 &1\ar[rd]^{c}&\\
 2\ar[ru]^{\sqrt 2}\ar[rr]_{\sqrt 2}&&3\\
 }}
 &
 \xymatrix@R=0.5cm@C=0.6cm
 {
 &\\
 \ar@{~>}[r]^{\mu_{121}}&\\
 }
 &
 \xymatrix@R=0.5cm@C=0.5cm
 {
 &1\ar[ld]_{\sqrt 2}&\\
 2\ar[rr]_{\sqrt 2}&&3\ar[lu]_{c+2}\\
 }
 &
 \xymatrix@R=0.5cm@C=0.6cm
 {
 &\\
 \ar@{-}[r]^{\mu_2}&\\
 }
 &
 \xymatrix@R=0.5cm@C=0.5cm
 {
 &1&\\
 2\ar[ru]^{\sqrt 2}&&3\ar[ll]^{\sqrt 2}\ar[lu]_c\\
 }
 &
 \xymatrix@R=0.5cm@C=0.6cm
 {
 &\\
 \ar@{-}[r]^{\mu_3}&\\
 }
 &
 {\color{blue}\xymatrix@R=0.5cm@C=0.5cm
 {
 &1\ar[rd]^{c}&\\
 2\ar[ru]^{\sqrt 2}\ar[rr]_{\sqrt 2}&&3\\
 }}
 \\
 \xymatrix@R=0.5cm@C=0.8cm
 {
 &\ar@{-}[d]^{\mu_3}&\\
 &&\\
 }&&&&&&
 \\
 \xymatrix@R=0.5cm@C=0.5cm
 {
 &1&\\
 2\ar[ru]^{\sqrt 2}&&3\ar[lu]_{c}\ar[ll]^{\sqrt 2}\\
 }
 &
 \xymatrix@R=0.5cm@C=0.6cm
 {
 &\\
 \ar@{-}[r]^{\mu_2}&\\
 }
 &
 \xymatrix@R=0.5cm@C=0.5cm
 {
 &1\ar[ld]_{\sqrt 2}&\\
 2\ar[rr]_{\sqrt 2}&&3\ar[lu]_{c+2}\\
 }
 &
 \xymatrix@R=0.5cm@C=0.6cm
 {
 &\\
 \ar@{~>}[r]^{\mu_{121}}&\\
 }
 &
 {\color{blue}\xymatrix@R=0.5cm@C=0.5cm
 {
 &1\ar[rd]^{c}&\\
 2\ar[ru]^{\sqrt 2}\ar[rr]_{\sqrt 2}&&3
 }}&&
 \\
 \xymatrix@R=0.5cm@C=0.8cm
 {
 &\ar@{-}[d]_{\mu_1}& \\
 &&\\
 }
 &
 &
 \xymatrix@R=0.5cm@C=0.8cm
 {
 &\ar@{~>}[d]_{\mu_{323}}& \\
 &&\\
 }&&&&\\
 {\color{blue}\xymatrix@R=0.5cm@C=0.5cm
 {
 &1\ar[rd]^{c}\ar[ld]_{\sqrt 2}&\\
 2&&3\ar[ll]^{\sqrt 2} 
 }}
 &
 &
 {\color{blue}\xymatrix@R=0.5cm@C=0.5cm
 {
 &1\ar[rd]^{c}\ar[ld]_{\sqrt 2}&\\
 2&&3\ar[ll]^{\sqrt 2} 
 }}&&&&
     \end{align*}
     }
     \caption{Mutation subgraph of $\Gamma_0$ of $(\sqrt{2},\sqrt{2},c\ge 2)$.}
     \label{fig:sqrt2sqrt2c-123}
 \end{figure}

Mutation subgraphs of $\Gamma_0$ are given in Figure \ref{fig:sqrt2sqrt2c-123}.

Firstly, assume that $d_1=d_3$. By Theorem \ref{thm:generator-set}, we conclude that $\operatorname{Aut}(\mathcal{A})$ is generated by
\begin{align*}
     g_{2}^{(1 3),-}, 
     g_{13}^{(1 3), -},
    g_{12121}^{(1 3), -},
    g_{32121}^{\mathbf{id},+},
    g_{323121}^{(1 3), -},
    g_{32323}^{(1 3), -},
    g_{12123}^{\mathbf{id},+}
    \in H_1(t_0).
\end{align*}
By Lemma \ref{lem:relation-rank-2-finite} and direction calculation, we obtain 
\[
    g_{12121}^{(1 3), -}\circ g_2^{(1 3), -}= g_2^{(1 3), -}\circ g_{32323}^{(1 3), -}=(g_{2}^{(1 3), -})^2=(g_{13}^{(1 3), -})^2=\mathbf{id},
    \]
    \[
    g_{32121}^{\mathbf{id},+}\circ g_{13}^{(1 3), -}=g_{12121}^{(1 3), -}, \quad
    g_{13}^{(1 3), -}\circ g_{12123}^{\mathbf{id},+}=g_{32323}^{(1 3), -}, \quad
    g_{12121}^{(1 3), -}\circ g_{12123}^{\mathbf{id},+}=g_{323121}^{(1 3), -}.
    \]
It follows that $\operatorname{Aut}(\mathcal{A})$ is generated by $g_{2}^{(1 3), -}$ and $g_{13}^{(1 3), -}$. Similar to the previous cases, one can show that \[\operatorname{Aut}(\mathcal{A})=\langle g_2^{(1 3), -},g_{13}^{(1 3), -}| (g_2^{(1 3), -})^2=(g_{13}^{(1 3), -})^2=\mathbf{id}\rangle \cong D_\infty.\]

Now assume that $d_1\ne d_3$. Every cluster automorphism is also a cluster automorphism for the case $d_1= d_3$. It is routine to check that  $\operatorname{Aut}(\mathcal{A})$ is the subgroup of $\langle g_2^{(1 3), -},g_{13}^{(1 3), -}\mid (g_2^{(1 3), -})^2=(g_{13}^{(1 3), -})^2=\mathbf{id}\rangle$ generated by $g_2^{(1 3), -}\circ g_{13}^{(1 3), -}=g_{312}^{\bf{id},+}$. Hence, $\operatorname{Aut}(\mathcal{A})\cong \mathbb{Z}$.

\subsubsection{$a=\sqrt{2},b=\sqrt{3},c\ge 2$.}
The mutation graph of $\Gamma_0$ are given in Figure \ref{fig:sqrt2sqrt3c-1} and \ref{fig:sqrt2sqrt3c-23}.
It is obvious that $t_\diamond:=\mu_{12121}(t_0)$ satisfies the condition $(\diamond)$.  According to Proposition \ref{prop:generator-set-assumption} and Remark \ref{rem:sharpen-K-2}, $\operatorname{Aut}(\mathcal{A})$ is generated by
\begin{align*}
     & g_{312}^{\mathbf{id},+}, g_{213}^{\mathbf{id},+}
     &\in H_{1,0}^{t_\diamond}(t_0),\\
   & g_{232312121}^{\mathbf{id},+},
    g_{3132312121}^{\mathbf{id},+},
    g_{21212121}^{\mathbf{id},+},
    g_{3212121}^{\mathbf{id},+},
    g_{121213232}^{\mathbf{id},+},
    g_{323232}^{\mathbf{id},+},
    g_{2123232}^{\mathbf{id},+},
    g_{2323212}^{\mathbf{id},+}
    &\in H_{1,1}^{t_\diamond}(t_0),\\
   & g_{31323212}^{\mathbf{id},+},
    g_{12121212}^{\mathbf{id},+},
    g_{1212132313}^{\mathbf{id},+},
    g_{3232313}^{\mathbf{id},+},
    g_{21232313}^{\mathbf{id},+},
    g_{232323}^{\mathbf{id},+},
    g_{3132323}^{\mathbf{id},+},
    g_{1212123}^{\mathbf{id},+}
   & \in H_{1,1}^{t_\diamond}(t_0),\\
   & g_{323231212121}^{\mathbf{id},+},
    g_{23232132312121}^{\mathbf{id},+}
   & \in K_{1}^{t_\diamond}(t_0).
\end{align*}
Similar to the previous cases, one can show that all the generators can be generated by $g_{312}^{\mathbf{id},+}$. It follows that $\operatorname{Aut}(\mathcal{A})\cong \mathbb{Z}$ by Lemma \ref{lem:order-of-tau}.



\begin{figure}
    \centering
{\tiny\bf
    \begin{align*}
t_0: {\color{blue} \xymatrix@R=0.5cm@C=0.5cm
{
&1\ar[rd]^{c}&\\
2\ar[ru]^{\sqrt 3}\ar[rr]_{\sqrt 2}&&3\\
} } 
&&
{\color{blue} 
\xymatrix@R=0.5cm@C=0.5cm
{
&1\ar[rd]^{c}&\\
2\ar[ru]^{\sqrt 3}\ar[rr]_{\sqrt 2}&&3\\
}}
&
&
{\color{blue}
\xymatrix@R=0.5cm@C=0.5cm
{
&1\ar[rd]^{c}&\\
2\ar[ru]^{\sqrt 3}\ar[rr]_{\sqrt 2}&&3\\
}}
&&
\\
\xymatrix@R=0.5cm@C=0.8cm
{
&\ar@{~>}[d]_{\mu_{12121}}& \\
&&\\
}
&&
\xymatrix@R=0.5cm@C=0.8cm
{
&\ar@{-}[d]_{\mu_{3}}& \\
&&\\
}
&&
\xymatrix@R=0.5cm@C=0.8cm
{
&\ar@{-}[d]_{\mu_{2}}& \\
&&\\
}
&&
\\
{\color{red}\xymatrix@R=0.5cm@C=0.5cm
{
&1\ar[ld]_{\sqrt 3}&\\
2\ar[rr]_{\sqrt 2}&&3\ar[lu]_{c+\sqrt 6}\\
}}
&
\xymatrix@R=0.5cm@C=0.6cm
{
&\\
\ar@{-}[r]^{\mu_2}&\\
}
&
\xymatrix@R=0.5cm@C=0.5cm
{
&1&\\
2\ar[ru]^{\sqrt 3}&&3\ar[lu]_{c}\ar[ll]^{\sqrt 2}\\
}
&
\xymatrix@R=0.5cm@C=0.6cm
{
&\\
\ar@{-}[r]^{\mu_1}&\\
}
&
\xymatrix@R=0.5cm@C=0.5cm
{
&1\ar[rd]^{c}\ar[ld]_{\sqrt 3}&\\
2&&3\ar[ll]^{\sqrt 2}\\
}
&
\xymatrix@R=0.5cm@C=0.6cm
{
&\\
\ar@{~>}[r]^{\mu_{323}}&\\
}
&
{\color{orange}
\xymatrix@R=0.5cm@C=0.5cm
{
&1\ar[ld]_{\sqrt 3}&\\
2\ar[rr]_{\sqrt 2}&&3\ar[lu]_{c+\sqrt 6}\\
}}
\\
\xymatrix@R=0.5cm@C=0.8cm
{
&\ar@{~>}[d]_{\mu_{323}}& \\
&&\\
}&&&&&&
\\
\xymatrix@R=0.5cm@C=0.5cm
{
&1\ar[rd]^{c}\ar[ld]_{\sqrt 3}&\\
2&&3\ar[ll]^{\sqrt 2}\\
}
&
\xymatrix@R=0.5cm@C=0.6cm
{
&\\
\ar@{-}[r]^{\mu_1}&\\
}
&
\xymatrix@R=0.5cm@C=0.5cm
{
&1&\\
2\ar[ru]^{\sqrt 3}&&3\ar[lu]_{c}\ar[ll]^{\sqrt 2}\\
}
&
\xymatrix@R=0.5cm@C=0.6cm
{
&\\
\ar@{-}[r]^{\mu_2}&\\
}
&
{\color{orange}
\xymatrix@R=0.5cm@C=0.5cm
{
&1\ar[ld]_{\sqrt 3}&\\
2\ar[rr]_{\sqrt 2}&&3\ar[lu]_{c+\sqrt 6}\\
}}&&
\\
\xymatrix@R=0.5cm@C=0.8cm
{
&\ar@{-}[d]_{\mu_{2}}& \\
&&\\
}
&
&
\xymatrix@R=0.5cm@C=0.8cm
{
&\ar@{-}[d]_{\mu_{3}}& \\
&&\\
}&&&&
\\
{\color{blue}
\xymatrix@R=0.5cm@C=0.5cm
{
&1\ar[rd]^{c}&\\
2\ar[ru]^{\sqrt 3}\ar[rr]_{\sqrt 2}&&3\\
}}
&
&
{
\color{blue}\xymatrix@R=0.5cm@C=0.5cm
{
&1\ar[rd]^{c}&\\
2\ar[ru]^{\sqrt 3}\ar[rr]_{\sqrt 2}&&3\\
}}
&&&&
\\
    \end{align*}
    }
    \caption{Mutation subgraph of $\Gamma_0$ of $(\sqrt{2},\sqrt{3},c\ge 2)$ in direction $1$.}
      \label{fig:sqrt2sqrt3c-1}
\end{figure}

\begin{figure}
    \centering
{\tiny\bf
    \begin{align*}
    &&
    {\color{orange}
\xymatrix@R=0.5cm@C=0.5cm
{
&1\ar[ld]_{\sqrt 3}&\\
2\ar[rr]_{\sqrt 2}&&3\ar[lu]_{c+\sqrt 6}\\
}}
&
&
{\color{orange}
\xymatrix@R=0.5cm@C=0.5cm
{
&1\ar[ld]_{\sqrt 3}&\\
2\ar[rr]_{\sqrt 2}&&3\ar[lu]_{c+\sqrt 6}\\
}}
&&
    \\
&&
\xymatrix@R=0.5cm@C=0.8cm
{
&& \\
&\ar@{~>}[u]_{\mu_{323}}&\\
}
&
&
\xymatrix@R=0.5cm@C=0.8cm
{
&\ar@{-}[d]_{\mu_{2}}& \\
&&\\
}
&&
    \\
 t_0:{\color{blue}\xymatrix@R=0.5cm@C=0.5cm
{
&1\ar[rd]^{c}&\\
2\ar[ru]^{\sqrt 3}\ar[rr]_{\sqrt 2}&&3\\
}}
&
\xymatrix@R=0.5cm@C=0.6cm
{
&\\
\ar@{-}[r]^{\mu_2}&\\
}
&
\xymatrix@R=0.5cm@C=0.5cm
{
&1\ar[rd]^{c}\ar[ld]_{\sqrt 3}&\\
2&&3\ar[ll]^{\sqrt 2} \\
}
&
\xymatrix@R=0.5cm@C=0.6cm
{
&\\
\ar@{-}[r]^{\mu_1}&\\
}
&
\xymatrix@R=0.5cm@C=0.5cm
{
&1&\\
2\ar[ru]^{\sqrt 3}&&3\ar[lu]_{c}\ar[ll]^{\sqrt 2}\\
}
&
\xymatrix@R=0.5cm@C=0.6cm
{
&\\
\ar@{-}[r]^{\mu_3}&\\
}
&
{\color{blue}
\xymatrix@R=0.5cm@C=0.5cm
{
&1\ar[rd]^{c}&\\
2\ar[ru]^{\sqrt 3}\ar[rr]_{\sqrt 2}&&3\\
}}
\\
\xymatrix@R=0.5cm@C=0.8cm
{
&\ar@{-}[d]_{\mu_{3}}& \\
&&\\
}&&&&&&
\\
\xymatrix@R=0.5cm@C=0.5cm
{
&1&\\
2\ar[ru]^{\sqrt 3}&&3\ar[lu]_{c}\ar[ll]^{\sqrt 2}\\
}
&
\xymatrix@R=0.5cm@C=0.6cm
{
&\\
\ar@{-}[r]^{\mu_1}&\\
}
&
\xymatrix@R=0.5cm@C=0.5cm
{
&1\ar[rd]^{c}\ar[ld]_{\sqrt 3}&\\
2&&3\ar[ll]^{\sqrt 2}\\
}
&
\xymatrix@R=0.5cm@C=0.6cm
{
&\\
\ar@{-}[r]^{\mu_2}&\\
}
&
{\color{blue}
\xymatrix@R=0.5cm@C=0.5cm
{
&1\ar[rd]^{c}&\\
2\ar[ru]^{\sqrt 3}\ar[rr]_{\sqrt 2}&&3\\
}}
&
&
\\
\xymatrix@R=0.5cm@C=0.8cm
{
&\ar@{-}[d]_{\mu_{2}}& \\
&&\\
}
&
&
\xymatrix@R=0.5cm@C=0.8cm
{
&\ar@{~>}[d]_{\mu_{323}}& \\
&&\\
}&&&&
\\
{\color{orange}
\xymatrix@R=0.5cm@C=0.5cm
{
&1\ar[ld]_{\sqrt 3}&\\
2\ar[rr]_{\sqrt 2}&&3\ar[lu]_{c+\sqrt 6}\\
}}
&
&
{\color{orange}
\xymatrix@R=0.5cm@C=0.5cm
{
&1\ar[ld]_{\sqrt 3}&\\
2\ar[rr]_{\sqrt 2}&&3\ar[lu]_{c+\sqrt 6}\\
}}
&&&&
\\
    \end{align*}
    }
    \caption{Mutation subgraph of $\Gamma_0$ of $(\sqrt{2},\sqrt{3},c\ge 2)$ in directions $2$ and $3$.}
      \label{fig:sqrt2sqrt3c-23}
\end{figure}

\subsubsection{$a=\sqrt{2},b\geq 2,c\geq 2$.}\label{ss:sqrt2-b-c}
It is clear that $d_2\neq d_3$ and $t_\diamond:=\mu_2(t_0)$ satisfies the condition $(\diamond)$. Mutation subgraph of $\Gamma_0$ is given in Figure \ref{fig:sqrt2bc-123}.


\begin{figure}
    \centering
{\tiny\bf
\begin{align*}
\xymatrix@R=0.5cm@C=0.5cm
{
&&\\
&\times&\\
}
&&&&
{\color{orange}\xymatrix@R=0.5cm@C=0.5cm
{
&1\ar[ld]_{b}\ar[rd]^c&\\
2&&3\ar[ll]^{\sqrt 2} \\
}}&&
\\
\xymatrix@R=0.5cm@C=0.8cm
{
&\ar@{-}[d]^{\mu_{1}}&\\
&&\\
}&&&&
\xymatrix@R=0.5cm@C=0.8cm
{
&&\\
&\ar@{~>}[u]^{\mu_{3232}}&\\
}
&&
\\
t_0:{\color{blue}\xymatrix@R=0.5cm@C=0.5cm
{
&1\ar[rd]^{c}&\\
2\ar[ru]^{b}\ar[rr]_{\sqrt 2}&&3\\
}}
&
\xymatrix@R=0.5cm@C=0.6cm
{
&\\
\ar@{-}[r]^{\mu_2}&\\
}
&
{\color{red}\xymatrix@R=0.5cm@C=0.5cm
{
&1\ar[ld]_{b}\ar[rd]^c&\\
2&&3\ar[ll]^{\sqrt 2} \\
}}
&
\xymatrix@R=0.5cm@C=0.6cm
{
&\\
\ar@{-}[r]^{\mu_1}&\\
}
&
\xymatrix@R=0.5cm@C=0.5cm
{
&1&\\
2\ar[ru]^{b}&&3\ar[ll]^{\sqrt 2}\ar[lu]_c\\
}
&
\xymatrix@R=0.5cm@C=0.6cm
{
&\\
\ar@{-}[r]^{\mu_3}&\\
}
&
{\color{blue}\xymatrix@R=0.5cm@C=0.5cm
{
&1\ar[rd]^{c}&\\
2\ar[ru]^{b}\ar[rr]_{\sqrt 2}&&3\\
}}
\\
&
&
\xymatrix@R=0.5cm@C=0.8cm
{
&\ar@{~>}[d]^{\mu_{2323}}&\\
&&\\
}
&
&
&
&
\\
{\color{orange}\xymatrix@R=0.5cm@C=0.5cm
{
&1\ar[ld]_{b}\ar[rd]^c&\\
2&&3\ar[ll]^{\sqrt 2} \\
}}
&
\xymatrix@R=0.5cm@C=0.6cm
{
&\\
\ar@{-}[r]^{\mu_1}&\\
}
&
\xymatrix@R=0.5cm@C=0.5cm
{
&1&\\
2\ar[ru]^b&&3\ar[ll]^{\sqrt 2}\ar[lu]_{c}\\
}
&
\xymatrix@R=0.5cm@C=0.6cm
{
&\\
\ar@{-}[r]^{\mu_3}&\\
}
&
{\color{blue}\xymatrix@R=0.5cm@C=0.5cm
{
&1\ar[rd]^{c}&\\
2\ar[ru]^{b}\ar[rr]_{\sqrt 2}&&3\\
}}
&
&
\end{align*}
  
\begin{align*}
t_0:{\color{blue}\xymatrix@R=0.5cm@C=0.5cm
{
&1\ar[rd]^{c}&\\
2\ar[ru]^{b}\ar[rr]_{\sqrt 2}&&3\\
}}
&
\xymatrix@R=0.5cm@C=0.6cm
{
&\\
\ar@{-}[r]^{\mu_3}&\\
}
&
\xymatrix@R=0.5cm@C=0.5cm
{
&1&\\
2\ar[ru]^{b}&&3\ar[lu]_{c}\ar[ll]^{\sqrt 2}\\
}
&
\xymatrix@R=0.5cm@C=0.6cm
{
&\\
\ar@{-}[r]^{\mu_1}&\\
}
&
{\color{orange}\xymatrix@R=0.5cm@C=0.5cm
{
&1\ar[rd]^{c}\ar[ld]_{b}&\\
2&&3\ar[ll]^{\sqrt 2} \\
}}
&
\xymatrix@R=0.5cm@C=0.6cm
{
&\\
\ar@{-}[r]^{\mu_2}&\\
}
&
{\color{blue}\xymatrix@R=0.5cm@C=0.5cm
{
&1\ar[rd]^{c}&\\
2\ar[ru]^{b}\ar[rr]_{\sqrt 2}&&3\\
}}
\\
&
\xymatrix@R=0.5cm@C=0.6cm
{
&\ar@{~>}[ld]^{\mu_{3232}}\\
&\\
}
&
&
&
\xymatrix@R=0.5cm@C=0.8cm
{
&\ar@{~>}[d]^{\mu_{2323}}&\\
&&\\
}
&
&
\\
{\color{orange}\xymatrix@R=0.5cm@C=0.5cm
{
&1\ar[ld]_{b}\ar[rd]^c&\\
2&&3\ar[ll]^{\sqrt 2} \\
}}
&
\xymatrix@R=0.5cm@C=0.6cm
{
&\\
&\\
}
&
{\color{orange}\xymatrix@R=0.5cm@C=0.5cm
{
&1\ar[ld]_{b}\ar[rd]^c&\\
2&&3\ar[ll]^{\sqrt 2} \\
}}
&
\xymatrix@R=0.5cm@C=0.6cm
{
&\\
\ar@{-}[r]^{\mu_1}&\\
}
&
\xymatrix@R=0.5cm@C=0.5cm
{
&1&\\
2\ar[ru]^b&&3\ar[ll]^{\sqrt 2}\ar[lu]_{c}\\
}
&
\xymatrix@R=0.5cm@C=0.6cm
{
&\\
\ar@{-}[r]^{\mu_3}&\\
}
&
{\color{blue}\xymatrix@R=0.5cm@C=0.5cm
{
&1\ar[rd]^{c}&\\
2\ar[ru]^{b}\ar[rr]_{\sqrt 2}&&3\\
}}
\end{align*}
}
    \caption{Mutation subgraph of $\Gamma_0$  of $(\sqrt{2},b\ge 2,c\ge 2)$.}
     \label{fig:sqrt2bc-123}
\end{figure}
According to Proposition \ref{prop:generator-set-assumption}, $\operatorname{Aut}(\mathcal{A})$ is generated by 
\begin{align*}
     & g_{312}^{\mathbf{id},+}, g_{213}^{\mathbf{id},+},
    g_{323232}^{\mathbf{id},+},
    g_{232323}^{\mathbf{id},+}, g_{3132323}^{\bf{id},+},
    g_{3232313}^{\mathbf{id},+}
    &\in H_{1,1}^{t_\diamond}(t_0),\\
    & g_{2323212}^{\mathbf{id},+},
    g_{2123232}^{\mathbf{id},+}
   & \in K_{1}^{t_\diamond}(t_0).
\end{align*}
Again, one can show that 
\[
     g_{3232313}^{\mathbf{id},+}=g_{213}^{\mathbf{id},+}\circ g_{323232}^{\mathbf{id},+}, g_{2123232}^{\mathbf{id},+}=g_{323232}^{\mathbf{id},+}\circ g_{213}^{\mathbf{id},+}, g_{3132323}^{\bf{id},+}\circ g_{213}^{\bf{id},+}=g_{232323}^{\bf{id},+},
\]
and 
\[
g_{323232}^{\mathbf{id},+}=g_{232323}^{\mathbf{id},+}=\mathbf{id},g_{213}^{\mathbf{id},+}\circ g_{2323212}^{\mathbf{id},+}=\mathbf{id},g_{213}^{\mathbf{id},+}\circ g_{312}^{\mathbf{id},+}=\mathbf{id}.
\]
Consequently,  $\operatorname{Aut}(\mathcal{A})$ is generated by $g_{312}^{\mathbf{id},+}$  and $\operatorname{Aut}(\mathcal{A})\cong \mathbb{Z}$.


\subsection{Case 4: $a=\sqrt{3}$.}

\subsubsection{$a=\sqrt{3},b=\sqrt{3},c\geq 2$.}
Mutation subgraph of $\Gamma_0$ is given in Figure \ref{fig:sqrt3sqrt3c-123}.
Note that this figure is almost identical to the Figure \ref{fig:sqrt2sqrt2c-123}. Similar to Section \ref{ss:a=sqrt2-b=sqrt2-cge2}, we obtain
\[
   \operatorname{Aut}(\mathcal{A})\cong
   \begin{cases}
       D_\infty & d_1=d_3;\\
       \mathbb{Z} &d_1\ne d_3.
   \end{cases}
\]

\begin{figure}[t]
    \centering
{\tiny\bf
    \begin{align*}
{\color{blue} \xymatrix@R=0.5cm@C=0.5cm
{
&1\ar[rd]^{c}\ar[ld]_{\sqrt 3}&\\
2&&3\ar[ll]^{\sqrt 3} \\
}}
&&
{\color{blue}\xymatrix@R=0.5cm@C=0.5cm
{
&1\ar[rd]^{c}\ar[ld]_{\sqrt 3}&\\
2&&3\ar[ll]^{\sqrt 3}\\
}}
&
&
{\color{blue}\xymatrix@R=0.5cm@C=0.5cm
{
&1\ar[rd]^{c}&\\
2\ar[ru]^{\sqrt 3}\ar[rr]_{\sqrt 3}&&3\\
}}
&&
 \\
    \xymatrix@R=0.5cm@C=0.8cm
{
&\ar@{-}[d]^{\mu_2}&\\
&&\\
}
    &&
\xymatrix@R=0.5cm@C=0.8cm
{
&& \\
&\ar@{~>}[u]_{\mu_{32323}}&\\
}
&
&
\xymatrix@R=0.5cm@C=0.8cm
{
&\ar@{-}[d]^{\mu_3}&\\
&&\\
}&&
    \\
  t_0:  {\color{blue}\xymatrix@R=0.5cm@C=0.5cm
{
&1\ar[rd]^{c}&\\
2\ar[ru]^{\sqrt 3}\ar[rr]_{\sqrt 3}&&3\\
}}
&
\xymatrix@R=0.5cm@C=0.6cm
{
&\\
\ar@{~>}[r]^{\mu_{12121}}&\\
}
&
\xymatrix@R=0.5cm@C=0.5cm
{
&1\ar[ld]_{\sqrt 3}&\\
2\ar[rr]_{\sqrt 3}&&3\ar[lu]_{c+3}\\
}
&
\xymatrix@R=0.5cm@C=0.6cm
{
&\\
\ar@{-}[r]^{\mu_2}&\\
}
&
\xymatrix@R=0.5cm@C=0.5cm
{
&1&\\
2\ar[ru]^{\sqrt 3}&&3\ar[lu]_{c}\ar[ll]^{\sqrt 3}\\
}
&
\xymatrix@R=0.5cm@C=0.6cm
{
&\\
\ar@{-}[r]^{\mu_1}&\\
}
&
{\color{blue}\xymatrix@R=0.5cm@C=0.5cm
{
&1\ar[rd]^{c}\ar[ld]_{\sqrt 3}&\\
2&&3\ar[ll]^{\sqrt 3}\\
}}
\\
\xymatrix@R=0.5cm@C=0.8cm
{
&\ar@{-}[d]^{\mu_3}&\\
&&\\
}
&&&&&&
\\
\xymatrix@R=0.5cm@C=0.5cm
{
&1&\\
2\ar[ru]^{\sqrt 3}&&3\ar[lu]_{c}\ar[ll]^{\sqrt 3}\\
}
&
\xymatrix@R=0.5cm@C=0.6cm
{
&\\
\ar@{-}[r]^{\mu_2}&\\
}
&
\xymatrix@R=0.5cm@C=0.5cm
{
&1\ar[ld]_{\sqrt 3}&\\
2\ar[rr]_{\sqrt 3}&&3\ar[lu]_{c+3}\\
}
&
\xymatrix@R=0.5cm@C=0.6cm
{
&\\
\ar@{~>}[r]^{\mu_{12121}}&\\
}
&
{\color{blue}\xymatrix@R=0.5cm@C=0.5cm
{
&1\ar[rd]^{c}&\\
2\ar[ru]^{\sqrt 3}\ar[rr]_{\sqrt 3}&&3\\
}}
&
&
\\
\xymatrix@R=0.5cm@C=0.8cm
{
&\ar@{-}[d]^{\mu_1}&\\
&&\\
}
&
&
\xymatrix@R=0.5cm@C=0.8cm
{
&\ar@{~>}[d]^{\mu_{32323}}&\\
&&\\
}&&&&
\\
{\color{blue}\xymatrix@R=0.5cm@C=0.5cm
{
&1\ar[rd]^{c}\ar[ld]_{\sqrt 3}&\\
2&&3\ar[ll]^{\sqrt 3} \\
}}
&
&
{\color{blue}\xymatrix@R=0.5cm@C=0.5cm
{
&1\ar[rd]^{c}\ar[ld]_{\sqrt 3}&\\
2&&3\ar[ll]^{\sqrt 3} \\
}}
&&&&\\
    \end{align*}
    }
    \caption{Mutation subgraph of $\Gamma_0$ of $(\sqrt{3},\sqrt{3},c\ge 2)$.}
   \label{fig:sqrt3sqrt3c-123}
\end{figure}

\subsubsection{$a=\sqrt{3},b\geq 2,c\geq 2$.}
Mutation subgraph of $\Gamma_0$ is given in Figure \ref{fig:sqrt3bc-123}, which is almost identical to the Figure \ref{fig:sqrt2bc-123}. Similar to Section \ref{ss:sqrt2-b-c}, one can show that $\operatorname{Aut}(\mathcal{A})$ is generated by $g_{312}^{\bf{id},+}$, thus we have
\[
   \operatorname{Aut}(\mathcal{A})\cong \mathbb{Z}. 
\]




\begin{figure}
    \centering
{\tiny\bf
\begin{align*}
\xymatrix@R=0.5cm@C=0.5cm
{
&&\\
&\times&\\
}
&&&&
{\color{orange}\xymatrix@R=0.5cm@C=0.5cm
{
&1\ar[ld]_{b}\ar[rd]^c&\\
2&&3\ar[ll]^{\sqrt 3} \\
}}&&
\\
\xymatrix@R=0.5cm@C=0.8cm
{
&\ar@{-}[d]^{\mu_{1}}&\\
&&\\
}&&&&
\xymatrix@R=0.5cm@C=0.8cm
{
&&\\
&\ar@{~>}[u]^{\mu_{323232}}&\\
}
&&
\\
t_0:{\color{blue}\xymatrix@R=0.5cm@C=0.5cm
{
&1\ar[rd]^{c}&\\
2\ar[ru]^{b}\ar[rr]_{\sqrt 3}&&3\\
}}
&
\xymatrix@R=0.5cm@C=0.6cm
{
&\\
\ar@{-}[r]^{\mu_2}&\\
}
&
{\color{red}\xymatrix@R=0.5cm@C=0.5cm
{
&1\ar[ld]_{b}\ar[rd]^c&\\
2&&3\ar[ll]^{\sqrt 3} \\
}}
&
\xymatrix@R=0.5cm@C=0.6cm
{
&\\
\ar@{-}[r]^{\mu_1}&\\
}
&
\xymatrix@R=0.5cm@C=0.5cm
{
&1&\\
2\ar[ru]^{b}&&3\ar[ll]^{\sqrt 3}\ar[lu]_c\\
}
&
\xymatrix@R=0.5cm@C=0.6cm
{
&\\
\ar@{-}[r]^{\mu_3}&\\
}
&
{\color{blue}\xymatrix@R=0.5cm@C=0.5cm
{
&1\ar[rd]^{c}&\\
2\ar[ru]^{b}\ar[rr]_{\sqrt 3}&&3\\
}}
\\
&
&
\xymatrix@R=0.5cm@C=0.8cm
{
&\ar@{~>}[d]^{\mu_{232323}}&\\
&&\\
}
&
&
&
&
\\
{\color{orange}\xymatrix@R=0.5cm@C=0.5cm
{
&1\ar[ld]_{b}\ar[rd]^c&\\
2&&3\ar[ll]^{\sqrt 3} \\
}}
&
\xymatrix@R=0.5cm@C=0.6cm
{
&\\
\ar@{-}[r]^{\mu_1}&\\
}
&
\xymatrix@R=0.5cm@C=0.5cm
{
&1&\\
2\ar[ru]^b&&3\ar[ll]^{\sqrt 3}\ar[lu]_{c}\\
}
&
\xymatrix@R=0.5cm@C=0.6cm
{
&\\
\ar@{-}[r]^{\mu_3}&\\
}
&
{\color{blue}\xymatrix@R=0.5cm@C=0.5cm
{
&1\ar[rd]^{c}&\\
2\ar[ru]^{b}\ar[rr]_{\sqrt 3}&&3\\
}}
&
&
\end{align*}
  
\begin{align*}
t_0:{\color{blue}\xymatrix@R=0.5cm@C=0.5cm
{
&1\ar[rd]^{c}&\\
2\ar[ru]^{b}\ar[rr]_{\sqrt 3}&&3\\
}}
&
\xymatrix@R=0.5cm@C=0.6cm
{
&\\
\ar@{-}[r]^{\mu_3}&\\
}
&
\xymatrix@R=0.5cm@C=0.5cm
{
&1&\\
2\ar[ru]^{b}&&3\ar[lu]_{c}\ar[ll]^{\sqrt 3}\\
}
&
\xymatrix@R=0.5cm@C=0.6cm
{
&\\
\ar@{-}[r]^{\mu_1}&\\
}
&
{\color{orange}\xymatrix@R=0.5cm@C=0.5cm
{
&1\ar[rd]^{c}\ar[ld]_{b}&\\
2&&3\ar[ll]^{\sqrt 3} \\
}}
&
\xymatrix@R=0.5cm@C=0.6cm
{
&\\
\ar@{-}[r]^{\mu_2}&\\
}
&
{\color{blue}\xymatrix@R=0.5cm@C=0.5cm
{
&1\ar[rd]^{c}&\\
2\ar[ru]^{b}\ar[rr]_{\sqrt 3}&&3\\
}}
\\
&
\xymatrix@R=0.5cm@C=0.6cm
{
&\ar@{~>}[ld]^{\mu_{323232}}\\
&\\
}
&
&
&
\xymatrix@R=0.5cm@C=0.8cm
{
&\ar@{~>}[d]^{\mu_{232323}}&\\
&&\\
}
&
&
\\
{\color{orange}\xymatrix@R=0.5cm@C=0.5cm
{
&1\ar[ld]_{b}\ar[rd]^c&\\
2&&3\ar[ll]^{\sqrt 3} \\
}}
&
\xymatrix@R=0.5cm@C=0.6cm
{
&\\
&\\
}
&
{\color{orange}\xymatrix@R=0.5cm@C=0.5cm
{
&1\ar[ld]_{b}\ar[rd]^c&\\
2&&3\ar[ll]^{\sqrt 3} \\
}}
&
\xymatrix@R=0.5cm@C=0.6cm
{
&\\
\ar@{-}[r]^{\mu_1}&\\
}
&
\xymatrix@R=0.5cm@C=0.5cm
{
&1&\\
2\ar[ru]^b&&3\ar[ll]^{\sqrt 3}\ar[lu]_{c}\\
}
&
\xymatrix@R=0.5cm@C=0.6cm
{
&\\
\ar@{-}[r]^{\mu_3}&\\
}
&
{\color{blue}\xymatrix@R=0.5cm@C=0.5cm
{
&1\ar[rd]^{c}&\\
2\ar[ru]^{b}\ar[rr]_{\sqrt 3}&&3\\
}}
\end{align*}
}
    \caption{Mutation subgraph of $\Gamma_0$  of $(\sqrt{3},b\ge 2,c\ge 2)$.}
   \label{fig:sqrt3bc-123}
\end{figure}

\subsection{Case 5. $a\geq 2$.}
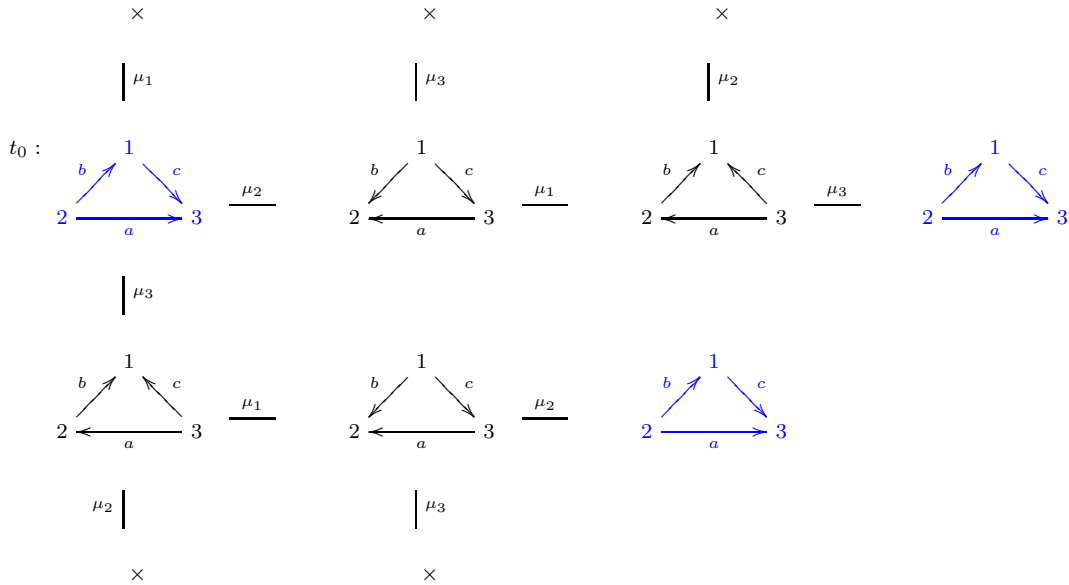
\begin{figure}
    \centering
{\tiny\bf
    \begin{align*}
    \xymatrix@R=0.5cm@C=0.5cm
{
&&\\
&\times& \\
}
&
&
\xymatrix@R=0.5cm@C=0.5cm
{
&&\\
&\times& \\
}
&
&
\xymatrix@R=0.5cm@C=0.5cm
{
&&\\
&\times& \\
}&&
\\
\xymatrix@R=0.5cm@C=0.8cm
{
&\ar@{-}[d]^{\mu_1}&\\
&&\\
}
&
&
\xymatrix@R=0.5cm@C=0.8cm
{
&\ar@{-}[d]^{\mu_3}&\\
&&\\
}
&
&
\xymatrix@R=0.5cm@C=0.8cm
{
&\ar@{-}[d]^{\mu_2}&\\
&&\\
}
&
&
\\
t_0: {\color{blue}\xymatrix@R=0.5cm@C=0.5cm
{
&1\ar[rd]^{c}&\\
2\ar[ru]^{b}\ar[rr]_{a}&&3\\
}}
&
\xymatrix@R=0.5cm@C=0.6cm
{
&\\
\ar@{-}[r]^{\mu_2}&\\
}
&
\xymatrix@R=0.5cm@C=0.5cm
{
&1\ar[ld]_{b}\ar[rd]^c&\\
2&&3\ar[ll]^{a} \\
}
&
\xymatrix@R=0.5cm@C=0.6cm
{
&\\
\ar@{-}[r]^{\mu_1}&\\
}
&
\xymatrix@R=0.5cm@C=0.5cm
{
&1&\\
2\ar[ru]^{b}&&3\ar[ll]^{a}\ar[lu]_c\\
}
&
\xymatrix@R=0.5cm@C=0.6cm
{
&\\
\ar@{-}[r]^{\mu_3}&\\
}
&
{\color{blue}\xymatrix@R=0.5cm@C=0.5cm
{
&1\ar[rd]^{c}&\\
2\ar[ru]^{b}\ar[rr]_a&&3\\
}}
\\
\xymatrix@R=0.5cm@C=0.8cm
{
&\ar@{-}[d]^{\mu_3}&\\
&&\\
}&&&&&&\\
\xymatrix@R=0.5cm@C=0.5cm
{
&1&\\
2\ar[ru]^{b}&&3\ar[lu]_{c}\ar[ll]^{a}\\
}
&
\xymatrix@R=0.5cm@C=0.6cm
{
&\\
\ar@{-}[r]^{\mu_1}&\\
}
&
\xymatrix@R=0.5cm@C=0.5cm
{
&1\ar[rd]^{c}\ar[ld]_{b}&\\
2&&3\ar[ll]^{a} \\
}
&
\xymatrix@R=0.5cm@C=0.6cm
{
&\\
\ar@{-}[r]^{\mu_2}&\\
}
&
{\color{blue}\xymatrix@R=0.5cm@C=0.5cm
{
&1\ar[rd]^{c}&\\
2\ar[ru]^{b}\ar[rr]_{a}&&3\\
}}
\\
\xymatrix@R=0.5cm@C=0.8cm
{
&\ar@{-}[d]_{\mu_2}& \\
&&\\
}
&
&
\xymatrix@R=0.5cm@C=0.8cm
{
&\ar@{-}[d]^{\mu_3}&\\
&&\\
}
&&&&\\
\xymatrix@R=0.5cm@C=0.5cm
{
&\times&\\
&& \\
}
&
&
\xymatrix@R=0.5cm@C=0.5cm
{
&\times&\\
&& \\
}&&&&
    \end{align*}
    }
   \caption{Mutation subgraph of $\Gamma_0$ of $(a\ge 2,b\ge 2,c\ge 2)$.}
    \label{fig:abc-123}
\end{figure}

Mutation subgraph of $\Gamma_0$ is given by Figure \ref{fig:abc-123}. 
We clearly have $|G_0(t_0)|\neq 1$ if and only if $b=c$ and $d_2=d_3$. In this case, $\psi_{(23)}^-$ is a generator of $G_0(t_0)$.
By Theorem \ref{thm:generator-set}, it is routine to check that $\operatorname{Aut}(\mathcal{A})$ is generated by
\[
    \begin{cases}
        \psi_{(2 3)}^{-}, g_{2}^{(1  3),-}, g_{3}^{(1  2),-} &\text{if $a=b=c, d_1=d_2=d_3$};\\
        g_3^{(1  2), -}, g_{12}^{(1  2),-}& \text{if $a=b=c, d_1=d_2\neq d_3$};\\
        g_{2}^{(1  3), -},g_{13}^{(1  3),-} & \text{if $a=b=c,d_1=d_3\neq d_2$ or $a=b<c, d_1=d_3$};\\
        \psi_{(2  3)}^-, 
        g_{312}^{\bf{id},+},
        g_{213}^{\bf{id},+} &\text{if $a=b=c, d_2=d_3\neq d_1$ or $a<b=c,d_2=d_3$};\\
        g_{312}^{\bf{id},+},
        g_{213}^{\bf{id},+} &\text{else.}        
    \end{cases}
\]
and
\[
    \operatorname{Aut}(\mathcal{A})\cong \begin{cases}
    D_\infty&\text{if $a=b=c, d_1=d_2=d_3$};\\
        D_\infty& \text{if $a=b=c, d_1=d_2\neq d_3$};\\
        D_\infty & \text{if $a=b=c,d_1=d_3\neq d_2$ or $a=b<c, d_1=d_3$};\\
        D_\infty &\text{if $a=b=c, d_2=d_3\neq d_1$ or $a<b=c,d_2=d_3$};\\
        \mathbb{Z} &\text{else.}   
    \end{cases}
\]

\bibliographystyle{alpha}
\bibliography{ref}

\begin{thebibliography}{ABBS08}

\bibitem[ABBS08]{ABBS08}
Ibrahim Assem, Martin Blais, Thomas Br\"ustle, and Audrey Samson.
\newblock Mutation classes of skew-symmetric {$3\times 3$}-matrices.
\newblock {\em Comm. Algebra}, 36(4):1209--1220, 2008.

\bibitem[ASS12]{ASS2012}
Ibrahim Assem, Ralf Schiffler, and Vasilisa Shramchenko.
\newblock Cluster automorphisms.
\newblock {\em Proc. Lond. Math. Soc. (3)}, 104(6):1271--1302, 2012.

\bibitem[ASS14]{ASS14}
Ibrahim Assem, Vasilisa Shramchenko, and Ralf Schiffler.
\newblock Cluster automorphisms and compatibility of cluster variables.
\newblock {\em Glasg. Math. J.}, 56(3):705--720, 2014.

\bibitem[BD15]{BD15}
J\'er\'emy Blanc and Igor Dolgachev.
\newblock Automorphisms of cluster algebras of rank 2.
\newblock {\em Transform. Groups}, 20(1):1--20, 2015.

\bibitem[BM16]{BM16}
V\'eronique Bazier-Matte.
\newblock Unistructurality of cluster algebras of type {$\tilde{\Bbb A}$}.
\newblock {\em J. Algebra}, 464:297--315, 2016.

\bibitem[BMP20]{BMP20}
V\'eronique Bazier-Matte and Pierre-Guy Plamondon.
\newblock Unistructurality of cluster algebras from unpunctured surfaces.
\newblock {\em Proc. Amer. Math. Soc.}, 148(6):2397--2409, 2020.

\bibitem[BS15]{BS15}
Tom Bridgeland and Ivan Smith.
\newblock Quadratic differentials as stability conditions.
\newblock {\em Publ. Math. Inst. Hautes \'Etudes Sci.}, 121:155--278, 2015.

\bibitem[CHL24]{CHL24}
Wen Chang, Min Huang, and Jian-Rong Li.
\newblock Quasi-homomorphisms of quantum cluster algebras.
\newblock {\em J. Algebra}, 638:506--534, 2024.

\bibitem[CL20a]{CL20}
Peigen Cao and Fang Li.
\newblock The enough {{\(g\)}}-pairs property and denominator vectors of
  cluster algebras.
\newblock {\em Math. Ann.}, 377(3-4):1547--1572, 2020.

\bibitem[CL20b]{CL2020}
Peigen Cao and Fang Li.
\newblock Unistructurality of cluster algebras.
\newblock {\em Compos. Math.}, 156(5):946--958, 2020.

\bibitem[CLLP19]{CLLP19}
Peigen Cao, Fang Li, Siyang Liu, and Jie Pan.
\newblock A conjecture on cluster automorphisms of cluster algebras.
\newblock {\em Electron. Res. Arch.}, 27:1--6, 2019.

\bibitem[CS19]{CS19}
Wen Chang and Ralf Schiffler.
\newblock Cluster automorphisms and quasi-automorphisms.
\newblock {\em Adv. in Appl. Math.}, 110:342--374, 2019.

\bibitem[CZ16a]{CZ16a}
Wen Chang and Bin Zhu.
\newblock Cluster automorphism groups of cluster algebras of finite type.
\newblock {\em J. Algebra}, 447:490--515, 2016.

\bibitem[CZ16b]{CZ16b}
Wen Chang and Bin Zhu.
\newblock Cluster automorphism groups of cluster algebras with coefficients.
\newblock {\em Sci. China Math.}, 59(10):1919--1936, 2016.

\bibitem[CZ20]{CZ20}
Wen Chang and Bin Zhu.
\newblock Cluster automorphism groups and automorphism groups of exchange
  graphs.
\newblock {\em Pacific J. Math.}, 307(2):283--302, 2020.

\bibitem[DL23]{DL23}
Jinlei Dong and Fang Li.
\newblock Presentations of mapping class groups and an application to cluster
  algebras from surfaces.
\newblock {\em arXiv:2307.15227}, 2023.

\bibitem[Fra16]{Fr16}
Chris Fraser.
\newblock Quasi-homomorphisms of cluster algebras.
\newblock {\em Adv. in Appl. Math.}, 81:40--77, 2016.

\bibitem[FZ02]{FZ02}
Sergey Fomin and Andrei Zelevinsky.
\newblock Cluster algebras. {I}. {F}oundations.
\newblock {\em J. Amer. Math. Soc.}, 15(2):497--529, 2002.

\bibitem[FZ03]{FZ03}
Sergey Fomin and Andrei Zelevinsky.
\newblock Cluster algebras. {II}. {F}inite type classification.
\newblock {\em Invent. Math.}, 154(1):63--121, 2003.

\bibitem[Gu11]{Gu11}
Weiwen Gu.
\newblock A decomposition algorithm for the oriented adjacency graph of the
  triangulations of a bordered surface with marked points.
\newblock {\em Electron. J. Combin.}, 18(1):Paper 91, 45, 2011.

\bibitem[LL21]{LL21}
Siyang Liu and Fang Li.
\newblock Periodicities in cluster algebras and cluster automorphism groups.
\newblock {\em Algebra Colloq.}, 28(4):601--624, 2021.

\bibitem[Sev13]{Sev13}
Ahmet~I. Seven.
\newblock Mutation classes of skew-symmetrizable {$3\times 3$} matrices.
\newblock {\em Proc. Amer. Math. Soc.}, 141(5):1493--1504, 2013.

\end{thebibliography}

\end{document}